# Large deviations and stochastic calculus for large random matrices[*]

## Alice Guionnet

*UMPA, Ecole Normale Supérieure de Lyon, 46, allée d'Italie, 69364 Lyon Cedex 07, France*
*e-mail:* `Alice.GUIONNET@umpa.ens-lyon.fr`

**Abstract:** Large random matrices appear in different fields of mathematics and physics such as combinatorics, probability theory, statistics, operator theory, number theory, quantum field theory, string theory etc... In the last ten years, they attracted lots of interests, in particular due to a serie of mathematical breakthroughs allowing for instance a better understanding of local properties of their spectrum, answering universality questions, connecting these issues with growth processes etc. In this survey, we shall discuss the problem of the large deviations of the empirical measure of Gaussian random matrices, and more generally of the trace of words of independent Gaussian random matrices. We shall describe how such issues are motivated either in physics/combinatorics by the study of the so-called matrix models or in free probability by the definition of a non-commutative entropy. We shall show how classical large deviations techniques can be used in this context.

These lecture notes are supposed to be accessible to non probabilists and non free-probabilists.

---

[*]These notes were prepared for a 6 hours course at the XXIX Conference on Stochastic Processes and Applications and slightly completed for publication. This work was partially supported by the accord France-Bresil.

# Contents

# Chapter 1

# Introduction

Large random matrices have been studied since the thirties when Wishart [132] considered them to analyze some statistics problems. Since then, random matrices appeared in various fields of mathematics. Let us briefly summarize some of them and the mathematical questions they raised.

1. **Large random matrices and statistics** : In 1928, Wishart considered matrices of the form $\mathbf{Y}^{N,M} = \mathbf{X}^{N,M}(\mathbf{X}^{N,M})^*$ with an $N \times M$ matrix $\mathbf{X}^{N,M}$ with random entries. Typically, the matrix $\mathbf{X}^{N,M}$ is made of independent equidistributed vectors $\{X^1, \cdots, X^N\}$ in $\mathbb{C}^M$ with covariance matrix $\Sigma$, $(\Sigma)_{ij} = \mathbb{E}[X_i^1 X_j^1]$ for $1 \leq i,j \leq M$. Such random vectors naturally appear in multivariate analysis context where $\mathbf{X}^{N,M}$ is a data matrix, the column vectors of which represent an observation of a vector in $\mathbb{C}^M$. In such a setup, one would like to find the effective dimension of the system, that is the smallest dimension with which one can encode all the variations of the data. Such a principal components analysis is based on the study of the eigenvalues and eigenvectors of the covariance matrix $\mathbf{X}^{N,M}(\mathbf{X}^{N,M})^*$. When one assumes that the column vectors have i.i.d Gaussian entries, $\mathbf{Y}^{N,M}$ is called a standard Gaussian Wishart matrix. In statistics, it used to be reasonable to assume that $N/M$ was large. However, the case where $N/M$ is of order one is nowadays commonly considered, which corresponds to the cases where either the number of observations is rather small or when the dimension of the observation is very large. Such cases appear for instance in problems related with telecommunications and more precisely the analysis of cellular phones data, where a very large number of customers have to be treated simultaneously (see [70, 116, 120] and references therein). Other examples are provided in [80].
    In this setting, the main questions concern local properties of the spectrum (such as the study of the large $N, M$ behavior of the spectral radius of $\mathbf{Y}^{N,M}$, see [80], the asymptotic behaviour of the $k$ largest eigenvalues etc.), or the form of the eigenvectors of $\mathbf{Y}^{N,M}$ (see [120] and references therein).

2. **Large random matrices and quantum mechanics** : Wigner, in 1951 [131], suggested to approximate the Hamiltonians of highly excited nuclei by large random matrices. The basic idea is that there are so many phenomena going on in such systems that they can not be analyzed exactly and only a statistical approach becomes reasonable. The random matrices should be chosen as randomly as possible within the known physical restrictions of the model. For instance, he considered what we shall later on call Wigner's matrices, that is Hermitian (since the Hamiltonian has to be Hermitian) matrices with i.i.d entries (modulo the symmetry constraint). In the case where the system is invariant by time inversion, one can consider real symmetric matrices etc... As Dyson pointed out, the general idea is to chose the most random model within the imposed symmetries and to check if the theoretical predictions agree with the experiment, a disagreement pointing out that an important symmetry of the problem has been neglected. It turned out that experiments agreed exceptionally well with these models; for instance, it was shown that the energy states of the atom of hydrogen submitted to a strong magnetic field can be compared with the eigenvalues of an Hermitian matrix with i.i.d Gaussian entries. The book [59] summarizes a few similar experiments as well as the history of random matrices in quantum mechanics.
   In quantum mechanics, the eigenvalues of the Hamiltonian represent the energy states of the system. It is therefore important to study, following Wigner, the spectral distribution of the random matrix under study, but even more important, is its spacing distribution which represents the energy gaps and its extremal eigenvalues which are related with the ground states. Such questions were addressed in the reference book of M.L. Mehta [93], but got even more popular in mathematics since the work of C. Tracy et H. Widom [117] . It is also important to make sure that the results obtained do not depend on the details of the large random matrix models such as the law of the entries; this important field of investigation is often referred to as universality. An important effort of investigation was made in the last ten years in this direction for instance in [23], [54], [76],[89], [110], [112], [118], [102] ...
3. **Large random matrices and Riemann Zeta function :** The Riemann Zeta function is given by

$$\zeta(s) = \sum_{n=1}^{\infty} n^{-s}$$

   with $\text{Re}(s) > 1$ and can be analytically continued to the complex plane. The study of the zeroes of this function in the strip $0 \leq \text{Re}(s) < 1$ furnishes one of the most famous open problems. It is well known that $\zeta$ has trivial zeroes at $-2, -4, -6....$ and that its zeroes are distributed symmetrically with respect to the line $\text{Re}(s) = 2^{-1}$. The Riemann conjecture is that all the non trivial zeroes are located on this line. It was suggested by Hilbert and Polya that these zeroes might be related to the eigenvalues of a Hermitian operator, which would immediately imply that they are aligned. To investigate this idea, H. Montgomery (1972), assuming the Riemann conjecture, studied the number of zeroes of the zeta function in $\text{Re}(s) = 2^{-1}$ up to a distance $T$ of the real axis. His result suggests a striking similarity with corresponding statistics of the distribution of the eigenvalues of random Hermitian or unitary matrices when $T$ is large. Since then, an extensive literature was devoted to understand this relation. Let us only point out that the statistical evidence of this link can only be tested thanks to enormous numerical work, in particular due to A. Odlyzko [99, 100] who could determine a few hundred of millions of zeroes of Riemann zeta function around the $10^{20}$-th zeroes on the line $\text{Re}(s) = 2^{-1}$.

In somewhat the same direction, there is numerical evidence that the eigenvalues distribution of large Wigner matrices also describes the large eigenvalues of the Laplacian in some bounded domain such as the cardioid. This is related to quantum chaos since these eigenvalues describe the long time behavior of the classical ray dynamics in this domain (i.e. the billiard dynamics).

4. **Large random matrices and free probability** Free probability is a probability theory in a non-commutative framework. Probability measures are replaced by tracial states on von Neumann algebras. Free probability also contains the central notion of freeness which can be seen as a non-commutative analogue of the notion of independence. At the algebraic level, it can be related with the usual notion of freeness. This is why free probability could be well suited to solve important questions in von Neumann algebras, such as the question of isomorphism between free group factors. Eventhough this goal is not yet achieved, let us quote a few results on von Neumann algebras which were proved thanks to free probability machinery [56],[57], [124].

   In the 1990's, Voiculescu [121] proved that large random matrices are asymptotically free as their size go to infinity. Hence, large random matrices became a source for constructing many non-commutative laws, with nice properties with respect to freeness. Thus, free probability can be considered as the natural asymptotic large random matrices framework. Conversely, if one believes that any tracial state could be approximated by the empirical distribution of large matrices (which we shall define more precisely later), which would answer in the affirmative a well known question of A. Connes, then any tracial state could be obtained as such a limit.

   In this context, one often studies the asymptotic behavior of traces of polynomial functions of several random matrices with size going to infinity, trying to deduce from this limit either intuition or results concerning tracial states. For instance, free probability and large random matrices can be used to construct counter examples to some operator algebra questions.

5. **Combinatorics, enumeration of maps and matrix models**
   It is well known that the evaluation of the expectation of traces of random matrices possesses a combinatorial nature. For instance, if one considers a $N \times N$ symmetric or Hermitian matrix $\mathbf{X}_N$ with i.i.d centered entries with covariance $N^{-1}$, it is well known that $E[N^{-1}\text{Tr}(\mathbf{X}_N^p)]$ converges toward 0 if $p$ is odd and toward the Catalan number $C_{\frac{p}{2}}$ if $p$ is even. $C_p$ is the number of non crossing partitions of $\{1, \cdots, 2p\}$ and arises very often in combinatorics. This idea was pushed forward by J. Harer and D. Zagier [68] who computed exactly moments of the trace of $\mathbf{X}_N^p$ to enumerate maps with given number of vertices and genus. This combinatorial aspect of large random matrices was developed in the free probability context by R. Speicher [113].

This strategy was considerably generalized by 't Hooft who saw that matrix integrals such as

$$Z_N(P) = E[e^{N\text{Tr}(P(\mathbf{X}_N^1, \cdots, \mathbf{X}_N^k))}]$$

with a polynomial function $P$ and independent copies $\mathbf{X}_N^i$ of $\mathbf{X}_N$, can be seen as generating functions for the enumeration of maps of various types. The formal proof follows from Feynman diagrams expansion. This relation is nicely summarized in an article by A. Zvonkin [136] and we shall describe it more precisely in Chapter 5. One-matrix integrals can be used to enumerate various maps of arbitrary genus (maps with a given genus $g$ appearing as the $N^{-2g}$ correction terms in the expansion of $Z_N(P)$), and several matrix integrals can serve to consider the case where the vertices of these maps are colored, i.e. can take different states. For example, two-matrix integrals can therefore serve to define an Ising model on random graphs.

Matrix models were also used in physics to construct string theory models. Since string theory concerns maps with arbitrary genus, matrix models have to be considered at criticality and with temperature parameters well tuned with the dimension in order to have any relevance in this domain. It seems that this subject had a great revival in the last few years, but it seems still far from mathematical (or at least my) understanding.

Haar distributed Unitary matrices also can be used to enumerate combinatorial objects due to their relation with representations of the symmetric group (cf. [34] for instance). Nice applications to the enumeration of magic squares can be found in [38].

In this domain, one tries to estimate integrals such as $Z_N(P)$, and in particular tries to obtain the full expansion of $\log Z_N(P)$ in terms of the dimension $N$. This could be done rigorously so far only for one matrix models by use of Riemann-Hilbert problem techniques by J. Mc Laughlin et N. Ercolani [46]. First order asymptotics for a few several-matrix models could be obtained by orthogonal polynomial methods by M. L. Mehta [93, 90, 32] and by large deviations techniques in [61]. The physics literature on the subject is much more consistent as can be seen on the arxiv (see work by V. Kazakov, I. Kostov, M. Staudacher, B. Eynard, P. Zinn Justin etc.).

## 6. Large random matrices, random partitions and determinantal laws

It is well know [93] that Gaussian matrices have a determinantal form, i.e. the law of the eigenvalues $(\lambda_1, \cdots, \lambda_N)$ of a Wigner matrix with complex Gaussian entries (also called the **GUE**) is given by

$$dP(\lambda_1, \cdots, \lambda_N) = Z_N^{-1} \Delta(\lambda)^2 e^{-\frac{N}{4} \sum_{i=1}^N \lambda_i^2} \prod d\lambda_i$$

with $Z_N$ the normalizing constant and

$$\Delta(\lambda) = \prod_{i<j}(\lambda_i - \lambda_j) = \det \begin{pmatrix} 1 & \lambda_1 & \lambda_1^2 & \cdots & \lambda_1^{N-1} \\ 1 & \lambda_2 & \lambda_2^2 & \cdots & \lambda_2^{N-1} \\ . & . & . & . & . \\ 1 & \lambda_N & \lambda_N^2 & \cdots & \lambda_N^{N-1} \end{pmatrix}$$

Because $\Delta$ is a determinant, specific techniques can be used to study for instance the law of the top eigenvalue or the spacing distribution in the bulk or next to the top (cf. [117]). Such laws appear actually in different contexts such as random partitions as illustrated in the work of K. Johansson [77] or tilling problems [78]. For more general remarks on the relation between random matrices and random partitions, see [101].

In fact, determinantal laws appear naturally when non-intersecting paths are involved. Indeed, following [83], if $k_T$ is the transition probability of a homogeneous continuous Markov process, and $P_T^N$ the distribution of $N$ independent copies $X_t^N = (x_1(t), \cdots, x_N(t))$ of this process, then for any $X = (x_1, \cdots, x_N), x_1 < x_2 < \cdots < x_N, Y = (y_1, \cdots, y_N), y_1 < y_2 < \cdots < y_N$, the reflection principle shows that

$$P(X_N(0) = X, X_N(T) = Y | \forall t \geq 0, x_1(t) \leq x_2(t) \leq \cdots x_N(t))$$

$$= C(x) \det \left( (k_T(x_i, y_j))_{\substack{1 \leq i \leq N \\ 1 \leq j \leq N}} \right) \quad (1.0.1)$$

with

$$C(x)^{-1} = \int \det \left( (k_T(x_i, y_j))_{\substack{1 \leq i \leq N \\ 1 \leq j \leq N}} \right) dy.$$

This might provide an additional motivation to study determinantal laws. Even more striking is the occurrence of large Gaussian matrices laws for the problem of the longest increasing subsequence [8], directed polymers and the totally asymmetric simple exclusion process [75]. These relations are based on bijections with pairs of Young tableaux.

In fact, the law of the hitting time of the totally asymmetric simple exclusion process (TASEP) starting from Heaviside initial condition can be related with the law of the largest eigenvalue of a Wishart matrix. Let us remind the reader that the (TASEP) is a process with values in $\{0,1\}^{\mathbb{Z}}$, 0 representing the fact that the site is empty and 1 that it is occupied, the dynamics of which are described as follows. Each site of $\mathbb{Z}$ is equipped with a clock which rings at times with exponential law. When the clock rings at site $i$, nothing happens if there is no particle at $i$ or if there is one at $i+1$. Otherwise, the particle jumps from $i$ to $i+1$. Once this clock rang, it is replaced by a brand new independent clock. K. Johansson [75] considered these dynamics starting from the initial condition where there is no particles on $\mathbb{Z}^+$ but one particle on each site of $\mathbb{Z}^-$. The paths of the particles do not intersect by construction and therefore one can expect the law of the configurations to be determinantal. The main question to understand is to know where the particle which was at site $-N$, $N \in \mathbb{N}$, at time zero will be at time $T$. In other words, one wants to study the time $H(N, M)$ that the particle which was initially at $-N$ needs to get to $M - N$. K. Johansson [75] has shown that $H(M, N)$ has the same law as of the largest eigenvalue of a Gaussian complex Wishart matrix $X^{N+1,M}(X^{N+1,M})^*$ where $X^{N+1,M}$ is a $(N + 1) \times M$ matrix with i.i.d complex Gaussian entries with covariance $2^{-1}$. This remark allowed him to complete the law of large numbers result of Rost [106] by the study of the fluctuations of order $N^{\frac{1}{3}}$.

This paper opens the field of investigation of diverse growth processes (cf. Forrester [53]), to the problem of generalizing this result to different initial conditions or to other problems such as tilling models [78]. In this last context, one of the main results is the description of the fluctuation of the boundary of the tilling in terms of the Airy process (cf. M. Prahofer and H. Spohn [114] and K. Johansson [79]).

In this set of problems, one usually meets the problem of analyzing the largest eigenvalue of a large matrix, which is a highly non trivial analysis since the eigenvalues interact by a Coulomb gas potential.

In short, large random matrices became extremely fashionable during the last ten years. It is somewhat a pity that there is no good introductory book to the field. Having seen the six aspects of the topic I tried to describe above and imagining all those I forgot, the task looks like a challenge.

These notes are devoted to a very particular aspect of the study of large random matrices, namely, the study of the deviations of the law of large random matrices macroscopic quantities such as their spectral measures. It is only connected to points 4 and 5 listed above. Since large deviations results are refinements of law of large numbers theorems, let us briefly summarize these last results here.

It has been known since Wigner that the spectral measure of Wigner matrices converges toward the semicircle law almost surely. More precisely, let us consider a Wigner matrix, that is a $N \times N$ selfadjoint matrix $\mathbf{X}^N$ with independent (modulo the symmetry constraint) equidistributed centered entries with covariance $N^{-1}$. Let $(\lambda_1, \cdots, \lambda_N)$ be the eigenvalues of $\mathbf{X}^N$. Then, it was shown by Wigner [131], under appropriate assumptions on the moments of the entries, that the spectral measure $\hat{\mu}^N = N^{-1} \sum \delta_{\lambda_i}$ converges almost surely toward the

semi-circle distribution

$$\sigma(dx) = C\sqrt{4-x^2}1_{|x|\leq 2}dx.$$

This result was originally proved by estimating the moments $\{N^{-1}\text{Tr}((\mathbf{X}^N)^p), p \in \mathbb{N}\}$, which is a common strategy to study the spectral measure of self-adjoint random matrices.

This convergence can also be proved by considering the Stieljes transform of the spectral measure following Z. Bai [4], which demands less hypothesis on the moments of the entries of $\mathbf{X}^N$. In the case of Gaussian entries, this result can be easily deduced from the large deviation principle of Section 3. The convergence of the spectral measure was generalized to Wishart matrices (matrices of the form $\mathbf{X}^N R^N (\mathbf{X}^N)^*$ with a matrix $\mathbf{X}^N$ with independent entries and a diagonal matrix $R^N$) by Pastur and Marchenko [103]. Another interesting question is to wonder, if you are given two arbitrary large matrices $(A, B)$ with given spectrum, how the spectrum of the sum of these two matrices behave. Of course, this depends a lot on their eigenvectors. If one assumes that $A$ and $B$ have the same eigenvectors and i.i.d eigenvalues with law $\mu$ and $\nu$ respectively, the law of the eigenvalues of $A + B$ is the standard convolution $\mu * \nu$. On the contrary, if the eigenvectors of $A$ and $B$ are a priori not related, it is natural to consider $A + UBU^*$ with $U$ following the Haar measure on the unitary group. It was proved by D. Voiculescu [122] that the spectral measure of this sum converges toward the free convolution $\mu_A \boxplus \mu_B$ if the spectral measure of $A$ (resp. $B$) converges toward $\mu_A$ (resp . $\mu_B$) as the size of the matrices goes to infinity. More generally, if one considers the normalized trace of a word in two independent Wigner matrices then Voiculescu [122] proved that it converges in expectation (but actually also almost surely) toward a limit which is described by the trace of this word taken at two free semi-circular variables. We shall describe the notion of freeness in Chapter 6.

The question of the fluctuations of the spectral measure of random matrices was initiated in 1982 by D. Jonsson [81] for Wishart matrices by using moments method. This approach was applied and improved by A. Soshnikov an Y. Sinai [110] who considered Wigner matrices with non Gaussian entries but sufficient bounds on their moments and who obtained precise estimates on the moments $\{N^{-1}\text{Tr}((\mathbf{X}^N)^p), p \in \mathbb{N}\}$. Such results were generalized to the non-commutative setting where one considers polynomial functions of several independent random matrices by T. Cabanal Duvillard [28] and myself [60]. Recently, J. Mingo and R. Speicher [96] gave a combinatorial interpretation of the limiting covariance via a notion of second order freeness which places the problem of fluctuations to its natural non-commutative framework. They applied it with P. Sniady [97] to unitary matrices, generalizing to a non-commutative framework the results of P. Diaconis and M. Shahshahani [37] showing that traces of moments of unitary matrices converge towards Gaussian variables. In [60], I used the non-commutative framework to study fluctuations of the spectral measure of Gaussian band matrices, following an idea of D. Shlyakhtenko [109]. On the other hand, A. Khorunzhy, B. Khoruzhenko and L. Pastur [89] and more recently Z. Bai and J.F Yao [6] developed Stieljes transforms technology to study the central limit theorems for entries with eventually only the four first moments bounded. Such techniques apply at best to prove central limit theorem for nice analytic functions of the matrix under study. K. Johansson [73] considered Gaussian entries in order to take advantage that in this case, the eigenvalues have a simple joint law, given by a Coulomb gas type interaction. In this case, he could describe the optimal set of functions for which a central limit theorem can hold. Note here that in [60], the covariance is described in terms of a positive symmetric operator and therefore such an optimal set should be described as the domain of this operator. However, because this operator acts on non-commutative functions, its domain remains rather mysterious. A general combinatorial approach for understanding the fluctuations of band matrices with entries satisfying for instance Poincaré inequalities and rather general test functions has recently been undertaken by G. Anderson and O. Zeitouni [2].

In these notes, we shall study the error to the typical behavior in terms of large deviations in the cases listed above, with the restriction to Gaussian entries. They rely on a series of papers I have written on this subject with different coauthors [10, 18, 29, 30, 42, 60, 62, 64, 65] and try to give a complete accessible overview of this work to uninitiated readers. Some statements are improved or corrected and global introductions to free probability and hydrodynamics/large deviations techniques are given. While full proofs are given in Chapter 3 and rather detailed in Chapter 4, Chapter 7 only outlines how to adapt the ideas of Chapter 4 to the non-commutative setting. Chapter 5 uses the results of Chapter 1 and Chapter 4 to study matrix models. These notes are supposed to be accessible to non probabilists, if they assume some facts concerning Itô's calculus.

First, we shall consider the case of Wigner Gaussian matrices (see Chapter 3). The case of non Gaussian entries is still an open problem. We generalize our approach to non centered Gaussian entries in Chapter 4, which corresponds to the deviations of the law of the spectral measure of $A + \mathbf{X}$ with a deterministic diagonal matrix $A$ and a Wigner matrix $\mathbf{X}$. This result in turn gives the first order asymptotics of spherical integrals. The asymptotics of spherical integrals allows us to estimate matrix integrals in the case of quadratic (also called $AB$) interaction. Such a study puts on a firm ground some physics papers of Matytsin for instance. It is related with the enumeration of colored planar maps. We finally present the natural generalization of these results to several matrices, which deals with the so-called free micro-states entropy.

# Frequently used notations

For $N \in \mathbb{N}$, $\mathcal{M}_N$ will denote the set of $N \times N$ matrices with complex entries, $\mathcal{H}_N$ (resp. $\mathcal{S}_N$) will denote the set of $N \times N$ Hermitian (resp. symmetric) matrices. $U(N)$ (resp. $O(N)$, resp $S(N)$) will denote the unitary (resp. orthogonal, resp. symplectic) group. We denote Tr the trace on $\mathcal{M}_N$, $\text{Tr}(A) = \sum_{i=1}^{N} A_{ii}$ and tr the normalized trace $\text{tr}(A) = N^{-1}\text{Tr}(A)$.

To denote an ordered product of non-commutative variables $X_1, \cdots X_n$ (such as matrices), we write in short

$$X_1 X_2 \cdots X_n = \prod_{1 \leq i \leq n}^{\rightarrow} X_i.$$

$\mathbb{C}[X_1, \cdots, X_n]$ (resp. $\mathbb{C}\langle X_1, \cdots, X_n \rangle$) denotes the space of commutative (resp. non-commutative) polynomials in $n$ variables for which $\prod_{1 \leq i \leq n}^{\rightarrow} X_{j_i} = \prod_{1 \leq i \leq n}^{\rightarrow} X_{\sigma(j_i)}$ (resp. $\prod_{1 \leq i \leq n}^{\rightarrow} X_{j_i} \neq \prod_{1 \leq i \leq n}^{\rightarrow} X_{\sigma(j_i)}$) for all choices of indices $\{j_i, 1 \leq i \leq n, n \in \mathbb{N}\}$ (resp. eventually for a choice of $\{j_i, 1 \leq i \leq n, n \in \mathbb{N}\}$) and for all permutation $\sigma$ (resp. eventually for some permutation $\sigma$).

For a Polish space $X$, $\mathcal{P}(X)$ shall denote the set of probability measures on $X$. $\mathcal{P}(X)$ will be equipped with the usual weak topology, ie a sequence $\mu_n \in \mathcal{P}(X)$ converges toward $\mu$ iff for any bounded continuous function $f$ on $X$, $\mu_n(f)$ converges toward $\mu(f)$. Here, we denote in short

$$\mu(f) = \int f(x) d\mu(x).$$

For two Polish spaces $X, Y$ and a measurable function $\phi : X \rightarrow Y$, for any $\mu \in \mathcal{P}(X)$ we denote $\phi_\# \mu \in \mathcal{P}(Y)$ the push forward of $\mu$ by $\phi$, that is the probability measure on $Y$ such that for any bounded continuous $f : Y \rightarrow \mathbb{R}$,

$$\phi_\# \mu(f) = \int f(\phi(x)) d\mu(x).$$

For a given selfadjoint $N \times N$ matrix $\mathbf{A}$, we denote $(\lambda_1(\mathbf{A}), \cdots, \lambda_N(\mathbf{A}))$ its $N$ (real) eigenvalues and by $\hat{\mu}_\mathbf{A}^N$ its spectral measure

$$\hat{\mu}_\mathbf{A}^N = \frac{1}{N} \sum_{i=1}^{N} \delta_{\lambda_i(\mathbf{A})} \in \mathcal{P}(\mathbb{R}).$$

For two Polish spaces $X, Y$ we denote by $\mathcal{C}_b^0(X, Y)$ (or $\mathcal{C}(X, Y)$ when no ambiguity is possible) the space of bounded continuous functions from $X$ to $Y$. For instance, we shall denote $\mathcal{C}([0, 1], \mathcal{P}(\mathbb{R}))$ the set of continuous processes on $[0, 1]$ with values in the set $\mathcal{P}(\mathbb{R})$ of probability measures on $\mathbb{R}$, endowed with its usual weak topology. For a measurable set $\Omega$ of $\mathbb{R} \times [0, 1]$, $\mathcal{C}_b^{p,q}(\Omega)$ denotes the set of real-valued functions on $\Omega$ which are p times continuously differentiable with respect to the (first) space variable and $q$ times continuously differentiable with respect to the (second) time variable with bounded derivatives. $\mathcal{C}_c^{p,q}(\Omega)$ will denote the functions of $\mathcal{C}_b^{p,q}(\Omega)$ with compact support in the interior of the measurable set $\Omega$. For a probability measure $\mu$ on a Polish space $X$, $L_p(d\mu)$ denotes the space of measurable functions with finite $p^{th}$ moment under $\mu$. We shall say that an equality holds in the sense of distribution on a measurable set $\Omega$ if it holds, once integrated with respect to any $\mathcal{C}_c^{\infty,\infty}(\Omega)$ functions.

# Chapter 2

# Basic notions of large deviations

Since these notes are devoted to the proof of large deviations principles, let us remind the reader what is a large deviation principle and the few main ideas which are commonly used to prove it. We refer the reader to [41] and [43] for further developments. In what follows, $X$ will be a Polish space (that is a complete separable metric space). We then have

**Definition 2.1.**  • $I : X \to \mathbb{R}^+ \cup \{+\infty\}$ is a rate function, iff it is lower semi-continuous, i.e. its level sets $\{x \in X : I(x) \leq M\}$ are closed for any $M \geq 0$. It is a good rate function if its level sets $\{x \in X : I(x) \leq M\}$ are compact for any $M \geq 0$.
• A sequence $(\mu_N)_{N \in \mathbb{N}}$ of probability measures on $X$ satisfies a large deviation principle with speed (or in the scale) $a_N$ (going to infinity with $N$) and rate function $I$ iff
a) For any closed subset $F$ of $X$,

$$\limsup_{N \to \infty} \frac{1}{a_N} \log \mu_N(F) \leq -\inf_F I.$$

b) For any open subset $O$ of $X$,

$$\liminf_{N \to \infty} \frac{1}{a_N} \log \mu_N(O) \geq -\inf_O I.$$

The proof of a large deviation principle often proceeds first by the proof of a weak large deviation principle (which is defined as in definition (2.1) except that the upper bound is only required to hold for compact sets) and the so-called exponential tightness property

**Definition 2.2.** A sequence $(\mu_N)_{N \in \mathbb{N}}$ of probability measures on $X$ is exponentially tight iff there exists a sequence $(K_L)_{L \in \mathbb{N}}$ of compact sets such that

$$\limsup_{L \to \infty} \limsup_{N \to \infty} \frac{1}{a_N} \log \mu_N(K_L^c) = -\infty.$$

A weak large deviation principle is itself equivalent to the estimation of the probability of deviations towards small balls

**Theorem 2.3 ([41], Theorem 4.1.11).** *Let $\mathcal{A}$ be a base of the topology of $X$. For every $A \in \mathcal{A}$, define*

$$\mathcal{L}_A = -\liminf_{N\to\infty} \frac{1}{a_N} \log \mu_N(A)$$

*and*

$$I(x) = \sup_{A\in\mathcal{A}: x\in A} \mathcal{L}_A.$$

*Suppose that for all $x \in X$,*

$$I(x) = \sup_{A\in\mathcal{A}: x\in A} \left[ -\limsup_{N\to\infty} \frac{1}{a_N} \log \mu_N(A) \right]$$

*Then, $\mu_N$ satisfies a weak large deviation principle with rate function $I$.*

As an immediate corollary, we find that if $d$ is a distance on $X$ compatible with the weak topology and $B(x, \delta) = \{y \in X : d(y, x) < \delta\}$,

**Corollary 2.4.** *Assume that for all $x \in X$*

$$-I(x) := \limsup_{\delta\to 0} \limsup_{N\to\infty} \frac{1}{a_N} \log \mu_N(B(x,\delta)) = \liminf_{\delta\to 0} \liminf_{N\to\infty} \frac{1}{a_N} \log \mu_N(B(x,\delta)).$$

*Then, $\mu_N$ satisfies a weak large deviation principle with rate function $I$.*

From a given large deviation principle one can deduce a large deviation principle for other sequences of probability measures by using either the so-called contraction principle or Laplace's method. Namely, let us recall the contraction principle (cf. Theorem 4.2.1 in [41]) :

**Theorem 2.5.** *Assume that $(\mu_N)_{N\in\mathbb{N}}$ satisfies a large deviation principle with good rate function $I$ with speed $a_N$. Then for any function $F : X \to Y$ with values in a Polish space $Y$ which is continuous, the image $(F_\#\mu_N)_{N\in\mathbb{N}} \in \mathcal{P}(Y)^\mathbb{N}$ also satisfies a large deviation principle with the same speed and good rate function given for any $y \in Y$ by*

$$J(y) = \inf\{I(x) : F(x) = y\}.$$

Laplace's method (or Varadhan's Lemma) says the following (cf. Theorem 4.3.1 [41]):

**Theorem 2.6.** *Assume that $(\mu_N)_{N\in\mathbb{N}}$ satisfies a large deviation principle with good rate function $I$. Let $\to \mathbb{R}$ be a bounded continuous function. Then,*

$$\lim_{N\to\infty} \frac{1}{a_N} \log \int e^{a_N F(x)} d\mu_N(x) = \sup_{x\in X} \{F(x) - I(x)\}.$$

*Moreover, the sequence*

$$\nu_N(dx) = \frac{1}{\int e^{a_N F(y)} d\mu_N(y)} e^{a_N F(x)} d\mu_N(x) \in \mathcal{P}(X)$$

*satisfies a large deviation principle with good rate function*

$$J(x) = I(x) - F(x) - \sup_{y \in X}\{F(y) - I(y)\}.$$

Bryc's theorem ([41], Section 4.4) gives an inverse statement to Laplace theorem. Namely, assume that we know that for any bounded continuous function $F : X \to \mathbb{R}$, there exists

$$\Lambda(F) = \lim_{N \to \infty} \frac{1}{a_N} \log \int e^{a_N F(x)} d\mu_N(x) \qquad (2.0.1)$$

Then, Bryc's theorem says that $\mu_N$ satisfies a weak large deviation principle with rate function

$$I(x) = \sup_{F \in \mathcal{C}_b^0(X, \mathbb{R})} \{F(x) - \Lambda(F)\}. \qquad (2.0.2)$$

This actually provides another approach to proving large deviation principles : We see that we need to compute the asymptotics (2.0.1) for as many bounded continuous functions as possible. This in general can easily be done only for some family of functions (for instance, if $\mu_N$ is the law of $N^{-1} \sum_{i=1}^{N} x_i$ for independent equidistributed bounded random variable $x_i$'s, $a_N = N$, such quantities are easy to compute for linear functions $F$). This will always give a weak large deviation upper bound with rate function given as in (2.0.2) but where the supremum is only taken on this family of functions. The point is then to show that in fact this family is sufficient, in particular this restricted supremum is equal to the supremum over all bounded continuous functions.

# Chapter 3

# Large deviations for the spectral measure of large random matrices

## 3.1. Large deviations for the spectral measure of Wigner Gaussian matrices

Let $\mathbf{X}^{N,\beta} = \left(X_{ij}^{N,\beta}\right)$ be $N \times N$ real (resp. complex) Gaussian Wigner matrices when $\beta = 1$ (resp. $\beta = 2$, resp. $\beta = 4$) defined as follows. They are $N \times N$ self-adjoint random matrices with entries

$$X_{kl}^{N,\beta} = \frac{\sum_{i=1}^{\beta} g_{kl}^i e_\beta^i}{\sqrt{\beta N}}, \quad 1 \leq k < l \leq N, \quad X_{kk}^{N,\beta} = \sqrt{\frac{2}{\beta N}} g_{kk} e_\beta^1, \quad 1 \leq k \leq N$$

where $(e_\beta^i)_{1 \leq i \leq \beta}$ is a basis of $\mathbb{R}^\beta$, that is $e_1^1 = 1$, $e_2^1 = 1, e_2^2 = i$. This definition can be extended to the case $\beta = 4$ when $N$ is even by choosing $\mathbf{X}^{N,\beta} = \left(X_{ij}^{N,\beta}\right)_{1 \leq i,j \leq \frac{N}{2}}$ with $X_{kl}^{N,\beta}$ a $2 \times 2$ matrix defined as above but with $(e_\beta^k)_{1 \leq k \leq 4}$ the Pauli matrices

$$e_4^1 = \begin{pmatrix} 1 & 0 \\ 0 & 1 \end{pmatrix}, \quad e_4^2 = \begin{pmatrix} 0 & -1 \\ 1 & 0 \end{pmatrix}, \quad e_4^3 = \begin{pmatrix} 0 & -i \\ -i & 0 \end{pmatrix}, \quad e_4^4 = \begin{pmatrix} i & 0 \\ 0 & -i \end{pmatrix}.$$

$(g_{kl}^i, k \leq l, 1 \leq i \leq \beta)$ are independent equidistributed centered Gaussian variables with variance 1. $(\mathbf{X}^{N,2}, N \in \mathbb{N})$ is commonly referred to as the Gaussian Unitary Ensemble (**GUE**), $(\mathbf{X}^{N,1}, N \in \mathbb{N})$ as the Gaussian Orthogonal Ensemble (**GOE**) and $(\mathbf{X}^{N,4}, N \in \mathbb{N})$ as the Gaussian Symplectic Ensemble (**GSE**) since they can be characterized by the fact that their laws are invariant under the action of the unitary, orthogonal and symplectic group respectively (see [93]).

$\mathbf{X}^{N,\beta}$ has $N$ real eigenvalues $(\lambda_1, \lambda_2, \cdots, \lambda_N)$. Moreover, by invariance of the distribution of $\mathbf{X}^{N,1}$ (resp. $\mathbf{X}^{N,2}$, resp. $\mathbf{X}^{N,4}$) under the action of the orthogonal group $O(N)$ (resp. the unitary group $U(N)$, resp. the symplectic group $S(N)$), it is not hard to check that its eigenvectors will follow the Haar measure $m_N^\beta$ on $O(N)$ (resp. $U(N)$, resp. $S(N)$) in the case $\beta = 1$ (resp. $\beta = 2$, resp. $\beta = 4$). More precisely, a change of variable shows that for any Borel subset $A \subset \mathcal{M}_{N\times N}(\mathbb{R})$ (resp. $\mathcal{M}_{N\times N}(\mathbb{C})$),

$$\mathbb{P}\left(\mathbf{X}^{N,\beta} \in A\right) = \int 1_{UD(\lambda)U^* \in A} dm_\beta^N(U) dQ_\beta^N(\lambda) \qquad (3.1.1)$$

with $D(\lambda) = \text{diag}(\lambda_1, \lambda_2, \cdots, \lambda_N)$ the diagonal matrix with entries ($\lambda_1 \leq \lambda_2 \leq \cdots \leq \lambda_N$) and $Q_\beta^N$ the joint law of the eigenvalues given by

$$Q_\beta^N(d\lambda_1, \cdots, d\lambda_N) = \frac{1}{Z_\beta^N} \Delta(\lambda)^\beta \exp\{-\frac{\beta}{4} N \sum_{i=1}^N \lambda_i^2\} \prod_{i=1}^N d\lambda_i,$$

where $\Delta$ is the Vandermonde determinant $\Delta(\lambda) = \prod_{1 \leq i < j \leq N} |\lambda_i - \lambda_j|$ and $Z_\beta^N$ is the normalizing constant

$$Z_\beta^N = \int \cdot \int \prod_{1 \leq i < j \leq N} |\lambda_i - \lambda_j|^\beta \exp\{-\frac{\beta}{4} N \sum_{i=1}^N \lambda_i^2\} \prod_{i=1}^N d\lambda_i.$$

Such changes of variables are explained in details in the book in preparation of P. Forrester [53].

Using this representation, it was proved in [10] that the law of the spectral measure $\hat\mu^N = \frac{1}{N} \sum_{i=1}^N \delta_{\lambda_i}$, as a probability measure on $\mathbb{R}$, satisfies a large deviation principle.

In the following, we will denote $\mathcal{P}(\mathbb{R})$ the space of probability measure on $\mathbb{R}$ and will endow $\mathcal{P}(\mathbb{R})$ with its usual weak topology. We now can state the main result of [10].

**Theorem 3.1.**
*1) Let $I_\beta(\mu) = \frac{\beta}{4} \int x^2 d\mu(x) - \frac{\beta}{2} \Sigma(\mu) - \frac{3}{8}\beta$ with $\Sigma$ the non-commutative entropy*

$$\Sigma(\mu) = \int \int \log|x - y| d\mu(x) d\mu(y).$$

*Then :*
 *a. $I_\beta$ is well defined on $\mathcal{P}(\mathbb{R})$ and takes its values in $[0, +\infty]$.*
 *b. $I_\beta(\mu)$ is infinite as soon as $\mu$ satisfies one of the following conditions :*
 *b.1 : $\int x^2 d\mu(x) = +\infty$.*
 *b.2 : There exists a subset $A$ of $\mathbb{R}$ of positive $\mu$ mass but null logarithmic capacity, i.e. a set $A$ such that :*

$$\mu(A) > 0 \qquad \gamma(A) := \exp\left\{-\inf_{\nu \in \mathcal{M}_1^+(A)} \int \int \log \frac{1}{|x-y|} d\nu(x) d\nu(y)\right\} = 0.$$

   c. $I_\beta$ is a good rate function, i.e. $\{I_\beta \leq M\}$ is a compact subset of $\mathcal{P}(\mathbb{R})$ for $M \geq 0$.

   d. $I_\beta$ is a convex function on $\mathcal{P}(\mathbb{R})$.

   e. $I_\beta$ achieves its minimum value at a unique probability measure on $\mathbb{R}$ which is described as the Wigner's semicircular law $\sigma_\beta = (2\pi)^{-1}\sqrt{4-x^2}dx$.

2) *The law of the spectral measure* $\hat{\mu}^N = \frac{1}{N}\sum_{i=1}^{N}\delta_{\lambda_i}$ *on* $\mathcal{P}(\mathbb{R})$ *satisfies a full large deviation principle with good rate function* $I_\beta$ *in the scale* $N^2$.

**Proof :** We here skip the proof of 1.b.2. and 1.d and refer to [10] for these two points. The proof of the large deviation principle is rather clear; one writes the density $Q_\beta^N$ of the eigenvalues as

$$Q_\beta^N(d\lambda_1,\cdots,d\lambda_N) = \frac{1}{Z_\beta^N}\exp\{-\frac{\beta}{2}N^2\int_{x\neq y}f(x,y)d\hat{\mu}^N(x)d\hat{\mu}^N(y)\}\prod_{i=1}^{N}e^{-\frac{\beta}{4}\lambda_i^2}d\lambda_i,$$

with

$$f(x,y) = \log|x-y|^{-1} + \frac{1}{4}(x^2+y^2).$$

If the function $x,y \to 1_{x\neq y}f(x,y)$ were bounded continuous, the large deviation principle would result directly from a standard Laplace method (cf. Theorem 2.6), where the entropic term coming from the underlying Lebesgue measure could be neglected since the scale of the large deviation principle we are considering is $N^2 \gg N$. In fact, the main point is to deal with the fact that the logarithmic function blows up at the origin and to control what happens at the diagonal $\Delta := \{x = y\}$. In the sequel, we let $\bar{Q}_\beta^N$ be the non-normalized positive measure $\bar{Q}_\beta^N = Z_\beta^N Q_\beta^N$ and prove upper and lower large deviation estimates with rate function

$$J_\beta(\mu) = -\frac{\beta}{2}\int\int f(x,y)d\mu(x)d\mu(y).$$

This is of course enough to prove the second point of the theorem by taking $F = O = \mathcal{P}(\mathbb{R})$ to obtain

$$\lim_{N\to\infty}\frac{1}{N^2}\log Z_\beta^N = -\frac{\beta}{2}\inf_{\nu\in\mathcal{P}(\mathbb{R})}\int\int f(x,y)d\nu(x)d\nu(y).$$

To obtain that this limit is equal to $\frac{3}{8}\beta$, one needs to show that the infimum is taken at the semi-circle law and then compute its value. Alternatively, Selberg's formula (see [93]) allows to compute $Z_\beta^N$ explicitly from which its asymptotics are easy to get. We refer to [10] for this point. The upper bound is obtained as follows. Noticing that $\hat{\mu}^N \otimes \hat{\mu}^N(\Delta) = N^{-1}$ $Q_\beta^N$-almost surely (since the eigenvalues are almost surely distinct), we see that for any $M \in \mathbb{R}^+$, $Q_\beta^N$ a.s.,

$$\int\int 1_{x\neq y}f(x,y)d\hat{\mu}^N(x)d\hat{\mu}^N(x) \geq \int\int 1_{x\neq y}f(x,y)\wedge M d\hat{\mu}^N(x)d\hat{\mu}^N(x)$$
$$= \int\int f(x,y)\wedge M d\hat{\mu}^N(x)d\hat{\mu}^N(x) - \frac{M}{N}$$

Therefore, for any Borel subset $A \in \mathcal{P}(\mathbb{R})$, any $M \in \mathbb{R}^+$,

$$\bar{Q}_\beta^N \left( \frac{1}{N} \sum_{i=1}^N \delta_{\lambda_i} \in A \right) \leq \frac{1}{\sqrt{\beta\pi}^N} e^{-\frac{\beta N^2}{2} \inf_{\mu \in A}\{\int f(x,y) \wedge M d\mu(x) d\mu(y)\} + NM},$$
(3.1.2)

resulting with

$$\limsup_{N\to\infty} \frac{1}{N^2} \log \left( \bar{Q}_\beta^N(\frac{1}{N}\sum_{i=1}^N \delta_{\lambda_i} \in A) \right) \leq -\frac{\beta}{2} \inf_{\mu \in A} \{\int f(x,y) \wedge M d\mu(x) d\mu(y)\}.$$

We now show that if $A$ is closed, $M$ can be taken equal to infinity in the above right hand side.

We first observe that $I_\beta(\mu) = \frac{\beta}{2} \int f(x,y) d\mu(x) d\mu(y) - \frac{3}{8}\beta$ is a good rate function. Indeed, since $f$ is the supremum of bounded continuous functions, $I_\beta$ is lower semi-continuous, i.e. its level sets are closed. Moreover, because $f$ blows up when $x$ or $y$ go to infinity, its level sets are compact. Indeed, if $m = -\inf f$,

$$[\inf_{x,y \in [-A,A]^c}(f+m)(x,y)]\mu([-A,A]^c)^2 \leq \int (f+m)(x,y) d\mu(x) d\mu(y)$$
$$= \frac{2}{\beta} I_\beta(\mu) + \frac{3}{8} + m$$

resulting with
$$\{I_\beta \leq M\} \subset K_{2M}$$

where $K_M$ is the set

$$K_M = \cap_{A>0} \left\{ \mu \in \mathcal{P}(\mathbb{R}); \mu([-A,A]^c) \leq \sqrt{\frac{2\beta^{-1}M + m + \frac{3}{8}}{\inf_{x,y\in[-A,A]^c}(f+m)(x,y)}} \right\}.$$

$K_M$ is compact since $\inf_{x,y \in [-A,A]^c}(f+m)(x,y)$ goes to infinity with $A$. Hence, $\{I_\beta \leq M\}$ is compact, i.e. $I_\beta$ is a good rate function.

Moreover, (3.1.2) and the above observations show also that

$$\limsup_{M\to\infty} \limsup_{N\to\infty} \bar{Q}_\beta^N(\frac{1}{N}\sum_{i=1}^N \delta_{\lambda_i} \in K_M^c) = -\infty$$

insuring, with the uniform boundedness of $N^{-2} \log Z_\beta^N$, the exponential tightness of $Q_\beta^N$. Hence, me may assume that $A$ is compact, and actually a ball surrounding any given probability measure with arbitrarily small radius (see Chapter 2). Let $B(\mu, \delta)$ be a ball centered at $\mu \in \mathcal{P}(\mathbb{R})$ with radius $\delta$ for a distance compatible with the weak topology,

Since $\mu \to \int \int f(x,y) \wedge M d\mu(x) d\mu(y)$ is continuous for any $\mu \in \mathcal{P}(\mathbb{R})$, (3.1.2) shows that for any probability measure $\mu \in \mathcal{P}(\mathbb{R})$

$$\limsup_{\delta \to 0} \limsup_{N \to \infty} \frac{1}{N^2} \log \bar{Q}_\beta^N (\frac{1}{N} \sum_{i=1}^N \delta_{\lambda_i} \in B(\mu, \delta))$$
$$\leq -\int \int f(x,y) \wedge M d\mu(x) d\mu(y) \qquad (3.1.3)$$

We can finally let $M$ going to infinity and use the monotone convergence theorem which asserts that

$$\lim_{M \to \infty} \int f(x,y) \wedge M d\mu(x) d\mu(y) = \int f(x,y) d\mu(x) d\mu(y)$$

to conclude that

$$\limsup_{\delta \to 0} \limsup_{N \to \infty} \frac{1}{N^2} \log \bar{Q}_\beta^N (\frac{1}{N} \sum_{i=1}^N \delta_{\lambda_i} \in B(\mu, \delta)) \leq -\int \int f(x,y) d\mu(x) d\mu(y)$$

finishing the proof of the upper bound.

To prove the lower bound, we can also proceed quite roughly by constraining the eigenvalues to belong to very small sets. This will not affect the lower bound again because of the fast speed $N^2 \gg N \log N$ of the large deviation principle we are proving. Again, the difficulty lies in the singularity of the logarithm. The proof goes as follows; let $\nu \in \mathcal{P}(\mathbb{R})$. Since $I_\beta(\nu) = +\infty$ if $\nu$ has an atom, we can assume without loss of generality that it does not when proving the lower bound. We construct a discrete approximation to $\nu$ by setting

$$x^{1,N} = \inf\left\{x \mid \nu(]-\infty, x]) \geq \frac{1}{N+1}\right\}$$
$$x^{i+1,N} = \inf\left\{x \geq x^{i,N} \mid \nu(]x^{i,N}, x]) \geq \frac{1}{N+1}\right\} \qquad 1 \leq i \leq N-1.$$

and $\nu^N = \frac{1}{N} \sum_{i=1}^N \delta_{x^{i,N}}$ (note here that the choice of the length $(N+1)^{-1}$ of the intervals rather than $N^{-1}$ is only done to insure that $x^{N,N}$ is finite). Then, $\nu^N$ converges toward $\nu$ as $N$ goes to infinity. Thus, for any $\delta > 0$, for $N$ large enough, if we set $\Delta_N := \{\lambda_1 \leq \lambda_2 \cdots \leq \lambda_N\}$

$$\bar{Q}_\beta^N(\hat{\mu}^N \in B(\nu, \delta)) \geq Z_\beta^N Q_\beta^N(\{\max_{1 \leq i \leq N} |\lambda_i - x^{i,N}| < \frac{\delta}{2}\} \cap \Delta_N)) \qquad (3.1.4)$$
$$= \int_{\{|\lambda_i| < \frac{\delta}{2}\} \cap \Delta_N} \prod_{i<j} |x^{i,N} - x^{j,N} + \lambda_i - \lambda_j|^\beta \exp\{-\frac{N}{2} \sum_{i=1}^N (x^{i,N} + \lambda_i)^2\} \prod_{i=1}^N d\lambda_i$$

But, when $\lambda_1 < \lambda_2 \cdots < \lambda_N$ and since we have constructed the $x^{i,N}$'s such that $x^{1,N} < x^{2,N} < \cdots < x^{N,N}$, we have, for any integer numbers $(i,j)$, the lower bound

$$|x^{i,N} - x^{j,N} + \lambda_i - \lambda_j| > \max\{|x^{j,N} - x^{i,N}|, |\lambda_j - \lambda_i|\}.$$

As a consequence, (3.1.4) gives

$$\overline{Q}_\beta^N \left( \hat{\mu}^N \in B(\nu, \delta) \right) \geq \prod_{i+1<j} |x^{i,N} - x^{j,N}|^\beta$$
$$\times \prod_{i=1}^{N-1} |x^{i+1,N} - x^{i,N}|^{\frac{\beta}{2}} \exp\{-\frac{N}{2} \sum_{i=1}^N (|x^{i,N}| + \delta)^2\}$$
$$\times \int_{[-\frac{\delta}{2},\frac{\delta}{2}]^N \cap \Delta_N} \prod_{i=1}^{N-1} |\lambda_{i+1} - \lambda_i|^{\frac{\beta}{2}} \prod_{i=1}^N d\lambda_i \qquad (3.1.5)$$

Moreover, one can easily bound from below the last term in the right hand side of (3.1.5) and find

$$\int_{[-\frac{\delta}{2},\frac{\delta}{2}]^N \cap \Delta_N} \prod_{i=1}^{N-1} (\lambda_{i+1} - \lambda_i)^{\frac{\beta}{2}} \prod_{i=1}^N d\lambda_i \geq \left(\frac{1}{\beta/2+1}\right)^{(N-1)} \left(\frac{\delta}{2N}\right)^{(\beta/2+1)(N-1)+1} \qquad (3.1.6)$$

Hence, (3.1.5) implies :

$$\overline{Q}_\beta^N \left( \hat{\mu}^N \in B(\nu, \delta) \right) \geq \prod_{i+1<j} |x^{i,N} - x^{j,N}|^\beta \prod_{i=1}^{N-1} |x^{i+1,N} - x^{i,N}|^{\frac{\beta}{2}} e^{-\frac{N}{2} \sum_{i=1}^N (x^{i,N})^2}$$
$$\times \left(\frac{1}{\beta/2+1}\right)^{(N-1)} \left(\frac{\delta}{2N}\right)^{(\beta/2+1)(N-1)+1}$$
$$\times \exp\{-N\delta \sum_{i=1}^N |x^{i,N}| - N^2\delta^2\}. \qquad (3.1.7)$$

Moreover

$$\frac{1}{2(N+1)} \sum_{i=1}^N (x^{i,N})^2 - \frac{\beta}{2(N+1)^2} \sum_{i \neq j} \log|x^{i,N} - x^{j,N}|$$
$$\leq \frac{1}{2} \int x^2 d\nu(x) - \beta \int_{x^{1,N} \leq x < y \leq x^{N,N}} \log(y-x) d\nu(x) d\nu(y) \qquad (3.1.8)$$

Indeed, since $x \to \log(x)$ increases on $\mathbb{R}^+$, we notice that, with $I_i = [x^{i,N}, x^{i+1,N}]$,

$$\int_{x^{1,N} \leq x < y \leq x^{N,N}} \log(y-x) d\nu(x) d\nu(y) \tag{3.1.9}$$

$$\leq \sum_{1 \leq i \leq j \leq N-1} \log(x^{j+1,N} - x^{i,N}) \nu^{\otimes 2}((x,y) \in I_i \times I_j; x < y)$$

$$= \frac{1}{(N+1)^2} \sum_{i<j} \log|x^{i,N} - x^{j+1,N}| + \frac{1}{2(N+1)^2} \sum_{i=1}^{N-1} \log|x^{i+1,N} - x^{i,N}|.$$

The same arguments holds for $\int x^2 d\nu(x)$ and $\sum_{i=1}^{N}(x^{i,N})^2$. We can conclude that

$$\liminf_{N \to \infty} \frac{1}{N^2} \log \overline{Q}_\beta^N \left( \hat{\mu}^N \in B(\nu, \delta) \right)$$

$$\geq -\delta \int |x| d\nu(x) - \delta^2 + \liminf_{N \to \infty} \{\beta \int_{x^{1,N} \leq x < y \leq x^{N,N}} \log(y-x) d\nu(x) d\nu(y)$$

$$- \frac{1}{2} \int x^2 d\nu(x)\}. \tag{3.1.10}$$

$$= -\delta \int |x| d\nu(x) - \delta^2 - \frac{\beta}{2} \int \int f(x,y) d\nu(x) d\nu(y).$$

Letting $\delta$ going to zero gives the result.

Note here that we have used a lot monotonicity arguments to show that our approximation scheme converges toward the right limit. Such arguments can not be used in the setup considered by G. Ben Arous and O. Zeitouni [11] when the eigenvalues are complex; this is why they need to perform first a regularization by convolution.

### 3.2. Discussion and open problems

There are many ways in which one would like to generalize the previous large deviations estimates; the first would be to remove the assumption that the entries are Gaussian. It happens that such a generalization is still an open problem. In fact, when the entries of the matrix are not Gaussian anymore, the law of the eigenvalues is not independent of that of the eigenvectors and becomes rather complicated. With O. Zeitouni [64], I proved however that concentration of measures result hold. In fact, set

$$\mathbf{X}_A(\omega) = ((\mathbf{X}_A(\omega))_{ij})_{1 \leq i,j \leq N}, \quad \mathbf{X}_A(\omega) = \mathbf{X}_A^*(\omega), \quad (\mathbf{X}_A(\omega))_{ij} = \frac{1}{\sqrt{N}} A_{ij} \omega_{ij}$$

with

$$\omega := (\omega^R + i\omega^I) = (\omega_{ij})_{1 \leq i,j \leq N} = (\omega_{ij}^R + \sqrt{-1}\omega_{ij}^I)_{1 \leq i,j \leq N}, \omega_{ij} = \bar{\omega}_{ji},$$

$\{\omega_{ij}, 1 \le i \le j \le N\}$ independent complex random variables with laws $\{P_{ij}, 1 \le i \le j \le N\}$, $P_{ij}$ being a probability measure on $\mathcal{C}$ with

$$P_{ij}(\omega_{ij} \in \cdot) = \int \mathbb{1}_{u+iv \in \cdot} P_{ij}^R(du) P_{ij}^I(dv),$$

and $A$ is a self-adjoint non-random complex matrix with entries $\{A_{ij}, 1 \le i \le j \le N\}$ uniformly bounded by, say, $a$. Our main result is

**Theorem 3.2.** *a) Assume that the $(P_{ij}, i \le j, i, j \in \mathbb{N})$ are uniformly compactly supported, that is that there exists a compact set $K \subset \mathcal{C}$ so that for any $1 \le i \le j \le N$, $P_{ij}(K^c) = 0$. Assume $f$ is convex and Lipschitz. Then, for any $\delta > \delta_0(N) := 8|K|\sqrt{\pi}a|f|_{\mathcal{L}}/N > 0$,*

$$\mathbb{P}^N\left(|tr(f(\mathbf{X}_A(\omega))) - \mathbb{E}^N[tr(f(\mathbf{X}_A))]| \ge \delta\right) \le 4e^{-\frac{1}{16|K|^2 a^2 |f|_{\mathcal{L}}^2} N^2 (\delta - \delta_0(N))^2}.$$

*b) If the $(P_{ij}^R, P_{ij}^I, 1 \le i \le j \le N)$ satisfy the logarithmic Sobolev inequality with uniform constant $c$, then for any Lipschitz function $f$, for any $\delta > 0$,*

$$\mathbb{P}^N\left(|tr(f(\mathbf{X}_A(\omega))) - \mathbb{E}^N[tr(f(\mathbf{X}_A))]| \ge \delta\right) \le 2e^{-\frac{1}{8ca^2 |f|_{\mathcal{L}}^2} N^2 \delta^2}.$$

This result is a direct consequence of standard results about concentration of measure due to Talagrand and Herbst and the observation that if $f : \mathbb{R} \to \mathbb{R}$ is a Lipschitz function, then $\omega \to tr(f(\mathbf{X}_A(\omega)))$ is also Lipschitz and its Lipschitz constant can be evaluated, and that if $f$ is convex, $\omega \to tr(f(\mathbf{X}_A(\omega)))$ is also convex.

Note here that the matrix $A$ can be taken equal to $\{A_{ij} = 1, 1 \le i, j \le N\}$ to recover results for Wigner's matrices. However, the generalization is here costless and allows to include at the same time more general type of matrices such as band matrices or Wishart matrices. See a discussion in [64].

Eventhough this estimate is on the right large deviation scale, it does not precise the rate of deviation toward a given spectral distribution. This problem seems to be very difficult in general. The deviations of a empirical moments of matrices with eventually non-centered entries of order $N^{-1}$ are studied in [39]; in this case, deviations are typically produced by the shift of all the entries and the scaling allows to see the random matrix as a continuous operator. This should not be the case for Wigner matrices.

Another possible generalization is to consider another common model of Gaussian large random matrices, namely Wishart matrices. Sample covariance matrices (or Wishart matrices) are matrices of the form

$$Y_{N,M} = X_{N,M} T_M X_{N,M}^*.$$

Here, $X_{N,M}$ is an $N \times M$ matrix with centered real or complex i.i.d. entries of covariance $N^{-1}$ and $T_M$ is an $M \times M$ Hermitian (or symmetric) matrix.

These matrices are often considered in the limit where $M/N$ goes to a constant $\alpha > 0$. Let us assume that $M \leq N$, and hence $\alpha \in [0, 1]$, to fix the notations. Then, $Y_{N,M}$ has $N - M$ null eigenvalues. Let $(\lambda_1, \cdots, \lambda_M)$ be the $M$ non trivial remaining eigenvalues and denote $\hat{\mu}^M = M^{-1} \sum_{i=1}^{M} \delta_{\lambda_i}$. In the case where $T_M = I$ and the entries of $X_{N,M}$ are Gaussian, Hiai and Petz [71] proved that the law of $\hat{\mu}^M$ satisfies a large deviation principle. In this case, the joint law of the eigenvalues is given by

$$d\sigma_M^\beta(\lambda_1, \cdots, \lambda_M) = \frac{1}{Z_{T_M}^\beta} \prod_{i<j} |\lambda_i - \lambda_j|^\beta \prod_{i=1}^{N} \lambda_i^{\frac{\beta}{2}(N-M+1)-1} e^{-\frac{N}{2}\lambda_i} 1_{\lambda_i \geq 0} d\lambda_i$$

so that the analysis performed in the previous section can be generalized to this model. In the case where $T_M$ is a positive definite matrix, we have the formula

$$\begin{aligned} &d\sigma_M^\beta(\lambda_1, \cdots, \lambda_M) \\ &= \frac{1}{Z_{T_M}^\beta} \prod_{i<j} |\lambda_i - \lambda_j|^\beta I_M^{(\beta)}(\beta^{-1} D(\lambda), T_M^{-1}) \prod 1_{\lambda_i \geq 0} \lambda_i^{\frac{\beta}{2}(N-M+1)-1} d\lambda_i \end{aligned}$$

with $D(\lambda)$ the $M \times M$ diagonal matrix with entries $(\lambda_1, \cdots, \lambda_M)$ and $I_N^{(\beta)}$ the spherical integral

$$I_N^{(\beta)}(D_N, E_N) := \int \exp\{N \frac{\beta}{2} \text{Tr}(U D_N U^* E_N)\} dm_N^\beta(U)$$

with $m_N^\beta$ the Haar measure on the group of orthogonal matrices $O(N)$ if $\beta = 1$ (resp. $U(N)$ if $\beta = 2$, resp. $S(N)$ if $\beta = 4$). $Z_N^\beta(T_M)$ is the normalizing constant such that $\sigma_M^\beta$ has mass one. This formula can be found in [72, (58) and (95)]. Hence, we see that the study of the deviations of the spectral measure $\hat{\mu}^M = \frac{1}{M} \sum_{i=1}^{M} \delta_{\lambda_i}$ when the spectral measure of $T_M$ converges, is equivalent to the study of the asymptotics of the spherical integral $I_N^{(\beta)}(D_N, E_N)$ when the spectral measures of $D_N$ and $E_N$ converge.

The spherical integral also appears when one considers Gaussian Wigner matrices with **non centered** entries. Indeed, if we let

$$\mathbf{Y}^{N,\beta} = \mathbf{M}^N + \mathbf{X}^{N,\beta}$$

with a self adjoint deterministic matrix $\mathbf{M}^N$ (which can be taken diagonal without loss of generality) and a Gaussian Wigner matrix $\mathbf{X}^{N,\beta}$ as considered in the previous section, then the law of the eigenvalues of $\mathbf{Y}^{N,\beta}$ is given by

$$dQ_N(\lambda_1, \cdots, \lambda_N) =$$

$$= \frac{1}{Z^{N,\beta}} \prod_{i<j} |\lambda_i - \lambda_j|^\beta e^{-\frac{N}{2}\text{Tr}((\mathbf{M}^N)^2) - \frac{N}{2}\sum_{i=1}^{N}\lambda_i^2} I_N^{(\beta)}(D(\lambda), \mathbf{M}^N) \prod_{i=1}^{N} d\lambda_i. \tag{3.2.11}$$

We shall in the next section study the asymptotics of the spherical integral by studying the deviations of the spectral measure of $\mathbf{Y}^{N,\beta}$. This in turn provides the large deviation principle for the law of $\hat{\mu}^M$ under $\sigma_M^\beta$ by Laplace's method (see [65], Theorem 1.2).

Let us finish this section by mentioning two other natural generalizations in the context of large random matrices. The first consists in allowing the covariances of the entries of the matrix to depend on their position. More precisely, one considers the matrix $\mathbf{Z}^{N,\beta} = (z_{ij}^{N,\beta})_{1 \leq i,j \leq N}$ with $(z_{ij}^{N,\beta} = \bar{z}_{ji}^{N,\beta}, 1 \leq i \leq j \leq N)$ centered independent Gaussian variables with covariance $N^{-1}c_{ij}^{N,\beta}$. $c_{ij}^{N,\beta}$ is in general chosen to be $c_{ij}^{N,\beta} = \phi(\frac{i}{N}, \frac{j}{N})$. When $\phi(x,y) = 1_{|x+y| \leq c}$, the matrix is called a band matrix, but the case where $\phi$ is smooth can as well be studied as a natural generalization. The large deviation properties of the spectral measure of such a matrix was studied in [60] but only a large deviation upper bound could be obtained so far. Indeed, the non-commutativity of the matrices here plays a much more dramatic role, as we shall discuss later. In fact, it can be seen by the so-called linear trick (see [69] for instance) that studying the deviations of the words of several matrices can be related with the deviations of the spectral measure of matrices with complicated covariances and hence this last problem is intimately related with the study of the so-called microstates entropy discussed in Chapter 7. The second generalization is to consider a matrix where all the entries are independent, and therefore with a priori complex spectrum. Such a generalization was considered by G. Ben Arous and O. Zeitouni [11] in the case of the so-called Ginibre ensemble where the entries are identically distributed standard Gaussian variables, and a large deviation principle for the law of the spectral measure was derived. Let us further note that the method developed in the last section can as well be used if one considers random partitions following the Plancherel measure. In the unitary case, this law can be interpreted as a discretization of the GUE (see S.Kerov [85] or K. Johansson [74]) because the dimension of a representation is given in terms of a Coulomb gas potential. Using this similarity, one can prove large deviations principle for the empirical distribution of these random partitions (see [62]); the rate function is quite similar to that of the GUE except that, because of the discrete nature of the partitions, deviations can only occur toward probability measures which are absolutely continuous with respect to Lebesgue measure and with density bounded by one. More general large deviations techniques have been developed to study random partitions for instance in [40] for uniform distribution.

Let us finally mention that large deviations can as well be obtained for the law of the largest eigenvalue of the Gaussian ensembles (cf. e.g. [9], Theorem 6.2)

# Chapter 4

# Asymptotics of spherical integrals

In this chapter, we shall consider the spherical integral

$$I_N^{(\beta)}(D_N, E_N) := \int \exp\{N\frac{\beta}{2}\text{Tr}(UD_NU^*E_N)\}dm_N^\beta(U)$$

where $m_N^\beta$ denotes the Haar measure on the orthogonal (resp. unitary, resp. symplectic) group when $\beta = 1$ (resp. $\beta = 2$, resp. $\beta = 4$). This object is actually central in many matters; as we have seen, it describes the law of Gaussian Wishart matrices and non centered Gaussian Wigner matrices. It also appears in many matrix models described in physics; we shall describe this point in the next chapter. It is related with the characters of the symmetric group and Schur functions (cf. [107]) because of the determinantal formula below. We shall discuss this point in the next paragraph.

A formula for $I_N^{(2)}$ was obtained by Itzykson and Zuber (and more generally by Harish-Chandra), see [93, Appendix 5]; whenever the eigenvalues of $D_N$ and $E_N$ are distinct then

$$I_N^{(2)}(D_N, E_N) = \frac{\det\{\exp ND_N(i)E_N(j)\}}{\Delta(D_N)\Delta(E_N)},$$

where $\Delta(D_N) = \prod_{i<j}(D_N(j) - D_N(i))$ and $\Delta(E_N) = \prod_{i<j}(E_N(j) - E_N(i))$ are the VanderMonde determinants associated with $D_N, E_N$. Although this formula seems to solve the problem, it is far from doing so, due to the possible cancellations appearing in the determinant so that it is indeed completely unclear how to estimate the logarithmic asymptotics of such a quantity.

To evaluate this integral, we noticed with O. Zeitouni [65] that it was enough to derive a large deviation principle for the law of the spectral measure of non centered Wigner matrices. We shall detail this point in Section 4.1. To prove a large deviation principle for such a law, we improved a strategy initiated with T. Cabanal Duvillard in [29, 30] which consists in considering the matrix $\mathbf{Y}^{N,\beta} = D^N + \mathbf{X}^{N,\beta}$ as the value at time one of a matrix-valued process

$$\mathbf{Y}^{N,\beta}(t) = D^N + \mathbf{H}^{N,\beta}(t)$$

where $\mathbf{H}^{N,2}$ (resp. $\mathbf{H}^{N,1}$, resp. $\mathbf{H}^{N,4}$) is a Hermitian (resp. symmetric, resp. symplectic) Brownian motion, that is a Wigner matrix with Brownian motion entries. More explicitly, $H^{N,\beta}$ is a process with values in the set of $N \times N$ self-adjoint matrices with entries $\{H_{i,j}^{N,\beta}(t), t \geq 0, i \leq j\}$ constructed via independent real valued Brownian motions $(B_{i,j}^k, 1 \leq i \leq j \leq N, 1 \leq k \leq 4)$ by

$$H_{k,l}^{N,\beta} = \begin{cases} \frac{1}{\sqrt{\beta N}} \sum_{i=1}^{N} B_{k,l}^i e_\beta^i, & \text{if } k < l \\ \frac{\sqrt{2}}{\sqrt{\beta N}} B_{l,l}, & \text{if } k = l. \end{cases} \quad (4.0.1)$$

where $(e_\beta^k, 1 \leq k \leq \beta)$ is the basis of $\mathbb{R}^\beta$ described in the previous chapter. The advantage to consider the whole process $\mathbf{Y}^{N,\beta}$ is that we can then use stochastic differential calculus and standard techniques to study deviations of processes by martingales properties as initiated by Kipnis, Olla and Varadhan [88]. The idea to study properties of Gaussian variables by using their characterization as time marginals of Brownian motions is not new. At the level of deviations, it is very natural since we shall construct infinitesimally the paths to follow to create a given deviation. Actually, it seems to be the right way to consider $I_N^{(2)}$ when one realizes that it has a determinantal form according to (1.0.1) and so is by nature related with non-intersecting paths. There is still no more direct study of these asymptotics of the spherical integrals; eventhough B. Collins [34] tried to do it by expanding the exponential into moments and using cumulants calculus, obtaining such asymptotics would still require to be able to control the convergence of infinite signed series. In the physics literature, A. Matytsin [91] derived the same asymptotics for the spherical integrals than these we shall describe. His methods are quite different and only apply a priori in the unitary case. I think they might be written rigorously if one could a priori prove sufficiently strong convergence of the spherical integral as $N$ goes to infinity. As a matter of fact, I do not think there is any other formula for the limiting spherical integral in the physics literature, but mostly saddle point studies of this a priori converging quantity. I do however mention recent works of B. Eynard and als. [50] and M. Bertola [14] who produced a formula of the free energy for the model of matrices coupled in chain by means of residues technology. However, this corresponds to the case where the matrices $E_N, D_N$ of the spherical integral have a random spectrum submitted to a smooth polynomial potential and it is not clear how to apply such technology to the hard constraint case where the spectral measures of $D_N, E_N$ converge to a prescribed limit.

## 4.1. Asymptotics of spherical integrals and deviations of the spectral measure of non centered Gaussian Wigner matrices

Let $\mathbf{Y}^{N,\beta} = \mathbf{D}^N + \mathbf{X}^{N,\beta}$ with a deterministic diagonal matrix $\mathbf{D}^N$ and $\mathbf{X}^{N,\beta}$ a Gaussian Wigner matrix. We now show how the deviations of the spectral measure of $\mathbf{Y}^{N,\beta}$ are related to the asymptotics of the spherical integrals. To this end, we shall make the following hypothesis

**Hypothesis 4.1.**  1. There exists $d_{max} \in \mathbb{R}^+$ such that for any integer number $N$, $\hat{\mu}^N_{D_N}(\{|x| \geq d_{max}\}) = 0$ and that $\hat{\mu}^N_{D_N}$ converges weakly toward $\mu_D \in \mathcal{P}(\mathbb{R})$.
2. $\hat{\mu}^N_{E_N}$ converges toward $\mu_E \in \mathcal{P}(\mathbb{R})$ while $\hat{\mu}^N_{E_N}(x^2)$ stays uniformly bounded.

**Theorem 4.2.** *Under hypothesis 4.1,*
*1) There exists a function $g : [0,1] \times \mathbb{R}^+ \mapsto \mathbb{R}^+$, depending on $\mu_E$ only, such that $g(\delta, L) \to_{\delta \to 0} 0$ for any $L \in \mathbb{R}^+$, and, for $\hat{E}_N, \bar{E}_N$ such that*

$$d(\hat{\mu}^N_{\hat{E}_N}, \mu_E) + d(\hat{\mu}^N_{\bar{E}_N}, \mu_E) \leq \delta/2, \qquad (4.1.2)$$

*and*

$$\int x^2 d\hat{\mu}^N_{\hat{E}_N}(x) + \int x^2 d\hat{\mu}^N_{\bar{E}_N}(x) \leq L, \qquad (4.1.3)$$

*it holds that*

$$\limsup_{N \to \infty} \left| \frac{1}{N^2} \log \frac{I^{(\beta)}_N(D_N, \hat{E}_N)}{I^{(\beta)}_N(D_N, \bar{E}_N)} \right| \leq g(\delta, L).$$

*We define*

$$\bar{I}^{(\beta)}(\mu_D, \mu_E) = \limsup_{N \uparrow \infty} \frac{1}{N^2} \log I^{(\beta)}_N(D_N, E_N),$$

$$\underline{I}^{(\beta)}(\mu_D, \mu_E) = \liminf_{N \uparrow \infty} \frac{1}{N^2} \log I^{(\beta)}_N(D_N, E_N),$$

*By the preceding, $\bar{I}^{(\beta)}(\mu_D, \mu_E)$ and $\underline{I}^{(\beta)}(\mu_D, \mu_E)$ are continuous functions on $\{(\mu_E, \mu_D) \in \mathcal{P}(\mathbb{R})^2 : \int x^2 d\mu_E(x) + \int x^2 d\mu_D(x) \leq L\}$ for any $L < \infty$.*

*2) For any probability measure $\mu \in \mathcal{P}(\mathbb{R})$,*

$$\inf_{\delta \to 0} \liminf_{N \to \infty} \frac{1}{N^2} \log \mathbb{P}\left(d(\hat{\mu}^N_\mathbf{Y}, \mu) < \delta\right)$$

$$= \inf_{\delta \to 0} \limsup_{N \to \infty} \frac{1}{N^2} \log \mathbb{P}\left(d(\hat{\mu}^N_\mathbf{Y}, \mu) < \delta\right) := -J_\beta(\mu_D, \mu).$$

*3) We let, for any $\mu \in \mathcal{P}(\mathbb{R})$,*

$$I_\beta(\mu) = \frac{\beta}{4} \int x^2 d\mu(x) - \frac{\beta}{2} \int \log|x - y| d\mu(x) d\mu(y).$$

If $\hat{\mu}_{E_N}^N$ converges toward $\mu_E \in \mathcal{P}(\mathbb{R})$ with $I_\beta(\mu_E) < \infty$, we have

$$\begin{aligned} I^{(\beta)}(\mu_D, \mu_E) &:= \bar{I}^{(\beta)}(\mu_D, \mu_E) = \underline{I}^{(\beta)}(\mu_D, \mu_E) \\ &= -J_\beta(\mu_D, \mu_E) + I_\beta(\mu_E) - \inf_{\mu \in \mathcal{P}(\mathbb{R})} I_\beta(\mu) + \frac{\beta}{4} \int x^2 d\mu_D(x). \end{aligned}$$

Before going any further, let us point out that these results give interesting asymptotics for Schur functions which are defined as follows.

- a Young shape $\lambda$ is a finite sequence of non-negative integers $(\lambda_1, \lambda_2, \ldots, \lambda_l)$ written in non-increasing order. One should think of it as a diagram whose $i$th line is made of $\lambda_i$ empty boxes: for example,

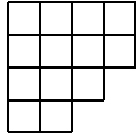 corresponds to $\lambda_1 = 4, \lambda_2 = 4, \lambda_3 = 3, \lambda_4 = 2$.

We denote by $|\lambda| = \sum_i \lambda_i$ the total number of boxes of the shape $\lambda$.
In the sequel, when we have a shape $\lambda = (\lambda_1, \lambda_2, \ldots)$ and an integer $N$ greater than the number of lines of $\lambda$ having a strictly positive length, we will define a sequence $l$ associated to $\lambda$ and $N$, which is an $N$-tuple of integers $l_i = \lambda_i + N - i$. In particular we have that $l_1 > l_2 > \ldots > l_N \geq 0$ and $l_i - l_{i+1} \geq 1$.
- for some fixed $N \in \mathbb{N}$, a Young tableau will be any filling of the Young shape above with integers from 1 to $N$ which is non-decreasing on each line and (strictly) increasing on each column. For each such filling, we define the content of a Young tableau as the $N$-tuple $(\mu_1, \ldots, \mu_N)$ where $\mu_i$ is the number of $i$'s written in the tableau.

For example,
$\begin{array}{|c|c|c|} \hline 1 & 1 & 2 \\ \hline 2 & 3 \\ \cline{1-2} 3 \\ \cline{1-1} \end{array}$
is allowed (and has content $(2, 2, 2)$),

whereas
$\begin{array}{|c|c|c|} \hline 1 & 1 & 2 \\ \hline 1 & 3 \\ \cline{1-2} 3 \\ \cline{1-1} \end{array}$
is not.

Notice that, for $N \in \mathbb{N}$, a Young shape can be filled with integers from 1 to $N$ if and only if $\lambda_i = 0$ for $i > N$.
- for a Young shape $\lambda$ and an integer $N$, the Schur polynomial $s_\lambda$ is an element of $\mathbb{C}[x_1, \ldots, x_N]$ defined by

$$s_\lambda(x_1, \ldots, x_N) = \sum_T x_1^{\mu_1} \ldots x_N^{\mu_N}, \qquad (4.1.4)$$

where the sum is taken over all Young tableaux $T$ of fixed shape $\lambda$ and $(\mu_1, \ldots, \mu_N)$ is the content of $T$. On a statistical point of view, one can think of the filling as the heights of a surface sitting on the tableau $\lambda$, $\mu_i$ being the surface of the height $i$. $s_\lambda$ is then a generating function for these heights when one considers the surfaces uniformly distributed under the constraints prescribed for the filling. Note that $s_\lambda$ is positive whenever the $x_i$'s are and, although it is not obvious from this definition (cf. for example [107] for a proof), $s_\lambda$ is a symmetric function of the $x_i$'s and actually $(s_\lambda, \lambda)$ form a basis of symmetric functions and hence play a key role in representation theory of the symmetric group. If $A$ is a matrix in $\mathcal{M}_N(\mathbb{C})$, then define $s_\lambda(A) \equiv s_\lambda(A_1, \ldots, A_N)$, where the $A_i$'s are the eigenvalues of $A$. Then, by Weyl formula (cf. Theorem 7.5.B of [130]), for any matrices $V, W$,

$$\int s_\lambda(UVU^*W) dm_N(U) = \frac{1}{d_\lambda} s_\lambda(V) s_\lambda(W). \qquad (4.1.5)$$

Then, because $s_\lambda$ also has a determinantal formula, we can see (cf. [107] and [62])

$$s_\lambda(M) = I_N^{(2)}\left(\log M, \frac{l}{N}\right) \Delta\left(\frac{l}{N}\right) \frac{\Delta(\log M)}{\Delta(M)}, \qquad (4.1.6)$$

where $\frac{l}{N}$ denotes the diagonal matrix with entries $N^{-1}(\lambda_i - i + N)$ and $\Delta$ the Vandermonde determinant. Therefore, we have the following immediate corollary to Theorem 4.2 :

**Corollary 4.3.** *Let $\lambda^N$ be a sequence of Young shapes and set $D_N = (N^{-1}(\lambda_i^N - i + N))_{1 \leq i \leq N}$. We pick a sequence of Hermitian matrix $E_N$ and assume that $(D_N, E_N)_{N \in \mathbb{N}}$ satisfy hypothesis 4.1 and that $\Sigma(\mu_D) > -\infty$. Then,*

$$\lim_{N \to \infty} \frac{1}{N^2} \log s_{\lambda^N}(e^{E_N})$$

$$= I^{(2)}(\mu_E, \mu_D) - \frac{1}{2} \int \log\left[\int_0^1 e^{\alpha x + (1-\alpha)y} d\alpha\right] d\mu_E(x) d\mu_E(y) + \frac{1}{2}\Sigma(\mu_D).$$

**Proof of Theorem 4.2 :** To simplify, let us assume that $E_N$ and $\hat{E}_N$ are uniformly bounded by a constant $M$. Let $\delta' > 0$ and $\{A_j\}_{j \in \mathcal{J}}$ be a partition of $[-M, M]$ such that $|A_j| \in [\delta', 2\delta']$ and the endpoints of $A_j$ are continuity points of $\mu_E$. Denote

$$\hat{I}_j = \{i : \hat{E}_N(ii) \in A_j\}, \quad \bar{I}_j = \{i : \bar{E}_N(ii) \in A_j\}.$$

By (4.1.2),

$$|\mu_E(A_j) - |\hat{I}_j|/N| + |\mu_E(A_j) - |\bar{I}_j|/N| \leq \delta.$$

We construct a permutation $\sigma_N$ so that $|\hat{E}(ii) - \bar{E}(\sigma_N(i), \sigma_N(i))| < 2\delta$ except possibly for very few $i$'s as follows. First, if $|\bar{I}_j| \leq |\hat{I}_j|$ then $\tilde{I}_j := \bar{I}_j$, whether if $|\bar{I}_j| > |\hat{I}_j|$ then $|\tilde{I}_j| = |\hat{I}_j|$ while $\tilde{I}_j \subset \bar{I}_j$. Then, choose and fix a permutation $\sigma_N$ such that $\sigma_N(\widetilde{I}_j) \subset \hat{I}_j$. Then, one can check that if $\mathcal{J}_0 = \{i : |\hat{E}(ii) - \bar{E}(\sigma_N(i), \sigma_N(i))| < 2\delta\}$,

$$\begin{aligned}
|\mathcal{J}_0| &\geq |\cup_j \sigma_N(\widetilde{I}_j)| = \sum_j |\sigma_N(\widetilde{I}_j)| \geq N - \sum_j |\bar{I}_j \backslash \widetilde{I}_j| \\
&\geq N - \max_j(|\bar{I}_j| - |\widetilde{I}_j|)|\mathcal{J}| \geq N - 2\delta N \frac{M}{\delta'}
\end{aligned}$$

Next, note the invariance of $I_N^{(\beta)}(D_N, E_N)$ to permutations of the matrix elements of $D_N$. That is,

$$\begin{aligned}
I_N^{(\beta)}(D_N, \bar{E}_N) &= \int \exp\{N\frac{\beta}{2}\text{Tr}(UD_NU^*\bar{E}_N)\}dm_N^\beta(U) \\
&= \int \exp\{N\frac{\beta}{2}\sum_{i,k} u_{ik}^2 D_N(kk)\bar{E}_N(ii)\}dm_N^\beta(U) \\
&= \int \exp\{N\sum_{i,k} u_{ik}^2 D_N(kk)\bar{E}_N(\sigma_N(i)\sigma_N(i))\}dm_N^\beta(U).
\end{aligned}$$

But, with $d_{\max} = \max_k |D_N(kk)|$ bounded uniformly in $N$,

$$\begin{aligned}
&N^{-1}\sum_{i,k} u_{ik}^2 D_N(kk)\bar{E}_N(\sigma_N(i)\sigma_N(i)) \\
&= N^{-1}\sum_{i \in \mathcal{J}_0}\sum_k u_{ik}^2 D_N(kk)\bar{E}_N(\sigma_N(i)\sigma_N(i)) \\
&\quad + N^{-1}\sum_{i \notin \mathcal{J}_0}\sum_k u_{ik}^2 D_N(kk)\bar{E}_N(\sigma_N(i)\sigma_N(i)) \\
&\leq N^{-1}\sum_{i,k} u_{ik}^2 D_N(kk)(\hat{E}_N(ii) + 2\delta) + N^{-1}d_{\max}M|\mathcal{J}_0^c| \\
&\leq N^{-1}\sum_{i,k} u_{ik}^2 D_N(kk)\hat{E}_N(ii) + d_{\max}\frac{M^2\delta}{\delta'}.
\end{aligned}$$

Hence, we obtain, taking $d_{\max}\frac{M^2\delta}{\delta'} = \sqrt{\delta}$,

$$I_N^{(\beta)}(D_N, \bar{E}_N) \leq e^{N\sqrt{\delta}}I_N^{(\beta)}(D_N, \hat{E}_N)$$

and the reverse inequality by symmetry. This proves the first point of the theorem when $(\bar{E}_N, \hat{E}_N)$ are uniformly bounded. The general case (which is not much more complicated) is proved in [65] and follows from first approximating $\bar{E}_N$ and $\hat{E}_N$ by bounded operators using (4.1.3).

The second and the third points are proved simultaneously : in fact, writing

$$\mathbb{P}\left(d(\hat{\mu}_{\mathbf{Y}}^N, \mu) < \delta\right) = \frac{1}{Z_N^\beta} \int_{d(\hat{\mu}_{\mathbf{Y}}^N, \mu) < \delta} e^{-\frac{N\beta}{4}\text{Tr}((\mathbf{Y}^{N,\beta} - D_N)^2)} d\mathbf{Y}^{N,\beta}$$

$$= \frac{e^{-\frac{N\beta}{4}\text{Tr}(D_N^2)}}{Z_N^\beta} \int_{d(\frac{1}{N}\sum_{i=1}^N \delta_{\lambda_i}, \mu) < \delta} I_N^{(\beta)}(D(\lambda), D_N) e^{-\frac{N\beta}{4} \sum_{i=1}^N \lambda_i^2} \Delta(\lambda)^\beta \prod_{i=1}^N d\lambda_i$$

with $Z_N^\beta$ the normalizing constant

$$Z_N^\beta = \int e^{-\frac{N}{2}\text{Tr}((\mathbf{Y}^{N,\beta} - D_N)^2)} d\mathbf{Y}^{N,\beta} = \int e^{-\frac{N}{2}\text{Tr}((\mathbf{Y}^{N,\beta})^2)} d\mathbf{Y}^{N,\beta},$$

we see that the first point gives, since $I_N^{(\beta)}(D(\lambda), D_N)$ is approximately constant on $\{d(N^{-1}\sum \delta_{\lambda_i}, \mu) < \delta\} \cap \{d(\hat{\mu}_{D_N}^N, \mu_D) < \delta\}$,

$$\mathbb{P}\left(d(\hat{\mu}_{\mathbf{Y}}^N, \mu) < \delta\right)$$
$$\approx \frac{e^{N^2(I^{(\beta)}(\mu_D, \mu) - \frac{\beta}{4}\hat{\mu}_{D_N}^N(x^2))}}{Z_N^\beta} \int_{d(\frac{1}{N}\sum_{i=1}^N \delta_{\lambda_i}, \mu) < \delta} e^{-\frac{N\beta}{4} \sum_{i=1}^N \lambda_i^2} \Delta(\lambda)^\beta \prod_{i=1}^N d\lambda_i$$
$$= e^{-\frac{N^2\beta}{4}\hat{\mu}_{D_N}^N(x^2) + N^2 I^{(\beta)}(\mu_D, \mu)} \mathbb{P}\left(d(\hat{\mu}_{\mathbf{X}}^N, \mu) < \delta\right)$$

where $A_{N,\delta} \approx B_{N,\delta}$ means that $N^{-2} \log A_{N,\delta} B_{N,\delta}^{-1}$ goes to zero as $N$ goes to infinity first and then $\delta$ goes to zero.

An equivalent way to obtain this relation is to use (1.0.1) together with the continuity of spherical integrals in order to replace the fixed values at time one of the Brownian motion $B_1 = Y$ by an average over a small ball $\{B_1 : d(\hat{\mu}_{B_1}^N, \mu) \leq \delta\}$.

The large deviation principle proved in the third chapter of these notes shows 2) and 3).

□

Note for 3) that if $I_\beta(\mu_E) = +\infty$, $J(\mu_D, \mu_E) = +\infty$ so that in this case the result is empty since it leads to an indetermination. Still, if $I_\beta(\mu_D) < \infty$, by symmetry of $I^{(\beta)}$, we obtain a formula by exchanging $\mu_D$ and $\mu_E$. If both $I_\beta(\mu_D)$ and $I_\beta(\mu_E)$ are infinite, we can only argue, by continuity of $I^{(\beta)}$, that for any sequence $(\mu_E^\epsilon)_{\epsilon>0}$ of probability measures with uniformly bounded variance and finite entropy $I_\beta$ converging toward $\mu_E$,

$$I^{(\beta)}(\mu_D, \mu_E) = \lim_{\epsilon \to \infty} \{-J_\beta(\mu_D, \mu_E^\epsilon) + I_\beta(\mu_E^\epsilon)\} - \inf I_\beta + \frac{\beta}{4} \int x^2 d\mu_D(x).$$

A more explicit formula is not yet available.

Note here that the convergence of the spherical integral is in fact not obvious and here given by the fact that we have convergence of the probability of deviation toward a given probability measure for the law of the spectral measure of non-centered Wigner matrices.

### 4.2. Large deviation principle for the law of the spectral measure of non-centered Wigner matrices

The goal of this section is to prove the following theorem.

**Theorem 4.4.** *Assume that $D_N$ is uniformly bounded with spectral measure converging toward $\mu_D$. Then the law of the spectral measure $\hat{\mu}_Y^N$ of the Wigner matrix $\mathbf{Y}^{N,\beta} = D_N + \mathbf{X}^{\beta,N}$ satisfies a large deviation principle with good rate function $J_\beta(\mu_D, .)$ in the scale $N^2$.*

By Bryc's theorem (2.0.2), it is clear that the above large deviation principle statement is equivalent to the fact that for any bounded continuous function $f$ on $\mathcal{P}(\mathbb{R})$,

$$\Lambda(f) = \lim_{N \to \infty} \frac{1}{N^2} \log \int e^{N^2 f(\hat{\mu}_Y^N)} d\mathbb{P}$$

exists and is given by $-\inf\{J_\beta(\mu_D, \nu) - f(\nu)\}$. It is not clear how one could a priori study such limits, except for very trivial functions $f$. However, if we consider the matrix valued process $\mathbf{Y}^{N,\beta}(t) = D^N + \mathbf{H}^{N,\beta}(t)$ with Brownian motion $\mathbf{H}^{N,\beta}$ described in (4.0.1) and its spectral measure process

$$\hat{\mu}_t^N = \hat{\mu}_{\mathbf{Y}^{N,\beta}(t)}^N = \frac{1}{N} \sum_{i=1}^N \delta_{\lambda_i(\mathbf{Y}^{N,\beta}(t))} \in \mathcal{P}(\mathbb{R}),$$

we may construct martingales by use of Itô's calculus. Continuous martingales lead to exponential martingales, which have constant expectation, and therefore allows one to compute the exponential moments of a whole family of functionals of $\hat{\mu}_.^N$. This idea will give easily a large deviation upper bound for the law of $(\hat{\mu}_t^N, t \in [0,1])$, and therefore for the law of $\hat{\mu}_Y^N$, which is the law of $\hat{\mu}_1^N$. The difficult point here is to check that it is enough to compute the exponential moments of this family of functionals in order to obtain the large deviation lower bound.

Let us now state more precisely our result. We shall consider $\{\hat{\mu}^N(t), t \in [0,1]\}$ as an element of the set $\mathcal{C}([0,1], \mathcal{P}(\mathbb{R}))$ of continuous processes with values in $\mathcal{P}(\mathbb{R})$. The rate function for these deviations shall be given as follows. For any $f, g \in \mathcal{C}_b^{2,1}(\mathbb{R} \times [0,1])$, any $s \leq t \in [0,1]$, and any $\nu_. \in \mathcal{C}([0,1], \mathcal{P}(\mathbb{R}))$, we let

$$\begin{aligned}
S^{s,t}(\nu, f) &= \int f(x,t) d\nu_t(x) - \int f(x,s) d\nu_s(x) \\
&\quad - \int_s^t \int \partial_u f(x,u) d\nu_u(x) du \\
&\quad - \frac{1}{2} \int_s^t \int\int \frac{\partial_x f(x,u) - \partial_x f(y,u)}{x-y} d\nu_u(x) d\nu_u(y) du, \quad (4.2.7)
\end{aligned}$$

$$<f, g>_\nu^{s,t} = \int_s^t \int \partial_x f(x,u) \partial_x g(x,u) d\nu_u(x) du, \quad (4.2.8)$$

and
$$\bar{S}^{s,t}(\nu, f) = S^{s,t}(\nu, f) - \frac{1}{2} < f, f >^{\nu}_{s,t} . \qquad (4.2.9)$$

Set, for any probability measure $\mu \in \mathcal{P}(\mathbb{R})$,

$$S_{\mu}(\nu) := \begin{cases} +\infty, & \text{if } \nu_0 \neq \mu, \\ S^{0,1}(\nu) := \sup_{f \in \mathcal{C}_b^{2,1}(\mathbb{R} \times [0,1])} \sup_{0 \leq s \leq t \leq 1} \bar{S}^{s,t}(\nu, f), & \text{otherwise.} \end{cases}$$

Then, the main theorem of this section is the following

**Theorem 4.5.** *1) For any $\mu \in \mathcal{P}(\mathbb{R})$, $S_\mu$ is a good rate function on $\mathcal{C}([0,1], \mathcal{P}(\mathbb{R}))$, i.e. $\{\nu \in \mathcal{C}([0,1], \mathcal{P}(\mathbb{R})); S_\mu(\nu) \leq M\}$ is compact for any $M \in \mathbb{R}^+$.*
*2) Assume that*

$$\text{there exists } \epsilon > 0 \text{ such that } \sup_N \hat{\mu}^N_{D_N}(|x|^{5+\epsilon}) < \infty, \quad \hat{\mu}^N_{D_N} \text{ converges toward } \mu_D,$$
(4.2.10)

*then the law of $(\hat{\mu}^N_t, t \in [0,1])$ satisfies a large deviation principle in the scale $N^2$ with good rate function $S_{\mu_D}$.*

<u>Remark</u> 4.6: In [65], the large deviation principle was only obtained for marginals; it was proved at the level of processes in [66].

Note that the application $(\mu_t, t \in [0,1]) \to \mu_1$ is continuous from $\mathcal{C}([0,1], \mathcal{P}(\mathbb{R}))$ into $\mathcal{P}(\mathbb{R})$, so that Theorem 4.5 and the contraction principle Theorem 2.5, imply that

**Theorem 4.7.** *Under assumption 4.2.10, Theorem 4.4 is true with*

$$J_\beta(\mu_D, \mu_E) = \frac{\beta}{2} \inf\{S_{\mu_D}(\nu_.); \nu_1 = \mu\}.$$

The main point to prove Theorem 4.5 is to observe that the evolution of $\hat{\mu}^N$ is described, thanks to Itô's calculus, by an autonomous differential equation. This is easily seen from the fact observed by Dyson [45] (see also [93], Theorem 8.2.1) that the eigenvalues $(\lambda^i_t, 1 \leq i \leq N, 0 \leq t \leq 1)$ of $(Y^{N,\beta}(t), 0 \leq t \leq 1)$ are described as the strong solution of the interacting particle system

$$d\lambda^i_t = \frac{\sqrt{2}}{\sqrt{\beta N}} dB^i_t + \frac{1}{N} \sum_{j \neq i} \frac{1}{\lambda^i_t - \lambda^j_t} dt \qquad (4.2.11)$$

with $\text{diag}(\lambda^1_0, \cdots, \lambda^N_0) = D^N$ and $\beta = 1, 2$ or $4$. This is the starting point to use Kipnis-Olla-Varadhan's papers ideas [88, 87]. These papers concerns the case where the diffusive term is not vanishing ($\beta N$ is of order one). The large deviations for the law of the empirical measure of the particles following (4.2.11) in such a scaling have been recently studied by Fontbona [52] in the context of Mc Kean-Vlasov diffusion with singular interaction. We shall first recall for the reader these techniques when one considers the empirical measures of independent Brownian motions as presented in [87]. We will then describe the necessary changes to adapt this strategy to our setting.

### 4.2.1. Large deviations from the hydrodynamical limit for a system of independent Brownian particles

Note that the deviations of the law of the empirical measure of independent Brownian motions on path space

$$L_N = \frac{1}{N} \sum_{i=1}^{N} \delta_{B^i_{[0,1]}} \in \mathcal{P}(\mathcal{C}([0,1], \mathbb{R}))$$

are well known by Sanov's theorem which yields (cf. [41], Section 6.2)

**Theorem 4.8.** *Let $\mathcal{W}$ be the Wiener law. Then, the law $(L_N)_\# \mathcal{W}^{\otimes N}$ of $L_N$ under $\mathcal{W}^{\otimes N}$ satisfies a large deviation principle in the scale $N$ with rate function given, for $\mu \in \mathcal{P}(\mathcal{C}([0,1], \mathbb{R}))$, by $I(\mu|\mathcal{W})$ which is infinite if $\mu$ is not absolutely continuous with respect to Lebesgue measure and otherwise given by*

$$I(\mu|\mathcal{W}) = \int \log \frac{d\mu}{d\mathcal{W}} \log \frac{d\mu}{d\mathcal{W}} d\mathcal{W}.$$

Thus, if we consider

$$\hat{\mu}_t^N = \frac{1}{N} \sum_{i=1}^{N} \delta_{B^i_t}, \quad t \in [0,1],$$

since $L_N \to (\hat{\mu}_t^N, t \in [0,1])$ is continuous from $\mathcal{P}(\mathcal{C}([0,1], \mathbb{R}))$ into $\mathcal{C}([0,1], \mathcal{P}(\mathbb{R}))$, the contraction principle shows immediately that the law of $(\hat{\mu}_t^N, t \in [0,1])$ under $\mathcal{W}^{\otimes N}$ satisfies a large deviation principle with rate function given, for $p \in \mathcal{C}([0,1], \mathcal{P}(\mathbb{R}))$, by

$$S(p) = \inf\{I(\mu|\mathcal{W}) \quad : \quad (x_t)_\# \mu = p_t \quad \forall t \in [0,1]\}.$$

Here, $(x_t)_\# \mu$ denotes the law of $x_t$ under $\mu$. It was shown by Föllmer [51] that in fact $S(p)$ is infinite unless there exists $k \in L^2(p_t(dx)dt)$ such that

$$\inf_{f \in \mathcal{C}^{1,1}(\mathbb{R} \times [0,1])} \int_0^1 \int (\partial_x f(x,t) - k(x,t))^2 p_t(dx) dt = 0, \quad (4.2.12)$$

and for all $f \in \mathcal{C}^{2,1}(\mathbb{R} \times [0,1])$,

$$\partial_t p_t(f_t) = p_t(\partial_t f_t) + \frac{1}{2} p_t(\partial_x^2 f_t) + p_t(\partial_x f_t k_t).$$

Moreover, we then have

$$S(p) = \frac{1}{2} \int_0^1 p_t(k_t^2) dt. \quad (4.2.13)$$

Kipnis and Olla [87] proposed a direct approach to obtain this result based on exponential martingales. Its advantage is to be much more robust and to adapt to many complicated settings encountered in hydrodynamics (cf. [86]). Let us now summarize it. It follows the following scheme

- *Exponential tightness and study of the rate function S* Since the rate function $S$ is the contraction of the relative entropy $I(.|\mathcal{W})$, it is clearly a good rate function. This can be proved directly from formula (4.2.13) as we shall detail it in the context of the eigenvalues of large random matrices. Similarly, we shall not detail here the proof that $\hat{\mu}_\#^N \mathcal{W}^{\otimes N}$ is exponentially tight which reduces the proof of the large deviation principle to the proof of a weak large deviation principle and thus to estimate the probability of deviations into small open balls (cf. Chapter 2). We will now concentrate on this last point.
- *Itô's calculus:* Itô's calculus (cf. [82], Theorem 3.3 p. 149) implies that for any function $F$ in $\mathcal{C}_b^{2,1}(\mathbb{R}^N \times [0,1])$, any $t \in [0,1]$

$$F(B_t^1, \cdots, B_t^N, t) = F(0, \cdots, 0) + \int_0^t \partial_s F(B_s^1, \cdots, B_s^N, s) ds$$

$$+ \sum_{i=1}^N \int_0^t \partial_{x_i} F(B_s^1, \cdots, B_s^N, s) dB_s^i + \frac{1}{2} \sum_{1 \le i,j \le N} \int_0^t \partial_{x_i} \partial_{x_j} F(B_s^1, \cdots, B_s^N, s) ds.$$

Moreover, $M_t^F = \sum_{i=1}^N \int_0^t \partial_{x_i} F(B_s^1, \cdots, B_s^N, s) dB_s^i$ is a martingale with respect to the filtration of the Brownian motion, with bracket

$$< M^F >_t = \sum_{i=1}^N \int_0^t [\partial_{x_i} F(B_s^1, \cdots, B_s^N, s)]^2 ds.$$

Taking $F(x^1, \cdots, x^N, t) = N^{-1} \sum_{i=1}^N f(B_t^i, t) = \int f(x,t) d\hat{\mu}_t^N(x) = \hat{\mu}_t^N(f_t)$, we deduce that for any $f \in \mathcal{C}_b^{2,1}(\mathbb{R} \times [0,1])$,

$$M_f^N(t) = \hat{\mu}_t^N(f_t) - \hat{\mu}_0^N(f_0) - \int_0^t \hat{\mu}_s^N(\partial_s f_s) ds - \int_0^t \hat{\mu}_s^N(\frac{1}{2}\partial_x^2 f_s) ds$$

is a martingale with bracket

$$< M_f^N >_t = \frac{1}{N} \int_0^t \hat{\mu}_s^N((\partial_x f_s)^2) ds.$$

The last ingredient of stochastic calculus we want to use is that (cf. [82], Problem 2.28, p. 147) for any bounded continuous martingale $m_t$ with bracket $< m >_t$, any $\lambda \in \mathbb{R}$,

$$\left\{ \exp(\lambda m_t - \frac{\lambda^2}{2} < m >_t), t \in [0,1] \right\}$$

is a martingale. In particular, it has constant expectation. Thus, we deduce that for all $f \in \mathcal{C}_b^{2,1}(\mathbb{R} \times [0,1])$, all $t \in [0,1]$,

$$\mathbb{E}[\exp\{N(M_f^N(t) - \frac{1}{2} < M_f^N >_t)\}] = 1. \qquad (4.2.14)$$

- *Weak large deviation upper bound*
  We equip $\mathcal{C}([0,1], \mathcal{P}(\mathbb{R}))$ with the weak topology on $\mathcal{P}(\mathbb{R})$ and the uniform topology on the time variable. It is then a Polish space. A distance compatible with such a topology is for instance given, for any $\mu, \nu \in \mathcal{C}([0,1], \mathcal{P}(\mathbb{R}))$, by

$$D(\mu, \nu) = \sup_{t \in [0,1]} d(\mu_t, \nu_t)$$

with a distance $d$ on $\mathcal{P}(\mathbb{R})$ compatible with the weak topology such as

$$d(\mu_t, \nu_t) = \sup_{|f|_\mathcal{L} \leq 1} |\int f(x) d\mu_t(x) - \int f(x) d\nu_t(x)|$$

where $|f|_\mathcal{L}$ is the Lipschitz constant of $f$:

$$|f|_\mathcal{L} = \sup_{x \in \mathbb{R}} |f(x)| + \sup_{x \neq y \in \mathbb{R}} |\frac{f(x) - f(y)}{x - y}|.$$

We prove here that

**Lemma 4.9.** *For any $p \in \mathcal{C}([0,1], \mathcal{P}(\mathbb{R}))$,*

$$\limsup_{\delta \to 0} \limsup_{N \to \infty} \frac{1}{N} \log \mathcal{W}^{\otimes N} \left( D(\hat{\mu}^N, p) \leq \delta \right) \leq -S(p).$$

**Proof :** Let $p \in \mathcal{C}([0,1], \mathcal{P}(\mathbb{R}))$. Observe first that if $p_0 \neq \delta_0$, since $\hat{\mu}_0^N = \delta_0$ almost surely,

$$\limsup_{\delta \to 0} \limsup_{N \to \infty} \frac{1}{N} \log \mathcal{W}^{\otimes N} \left( \sup_{t \in [0,1]} d(\hat{\mu}_t^N, p_t) \leq \delta \right) = -\infty.$$

Therefore, let us assume that $p_0 = \delta_0$. We set

$$B(p, \delta) = \{\mu \in \mathcal{C}([0,1], \mathcal{P}(\mathbb{R})) : D(\mu, p) \leq \delta\}.$$

Let us denote, for $f, g \in \mathcal{C}_b^{2,1}(\mathbb{R} \times [0,1])$, $\mu \in \mathcal{C}([0,1], \mathcal{P}(\mathbb{R}))$, $0 \leq t \leq 1$,

$$T^{0,t}(f, \mu) = \mu_t(f_t) - \mu_0(f_0) - \int_0^t \mu_s(\partial_s f_s) ds - \int_0^t \mu_s(\frac{1}{2} \partial_x^2 f_s) ds$$

and

$$<f, g>_\mu^{0,t} := \int_0^t \mu_s(\partial_x f_s \partial_x g_s) ds.$$

Then, by (4.2.14), for any $t \leq 1$,

$$\mathbb{E}[\exp\{N(T^{0,t}(f, \hat{\mu}^N) - \frac{1}{2} <f, f>_{\hat{\mu}^N}^{0,t})\}] = 1.$$

Therefore, if we denote in short $T(f, \mu) = T^{0,1}(f, \mu) - \frac{1}{2} < f, f >^{0,1}_\mu$,

$$\mathcal{W}^{\otimes N}\left(D(\hat{\mu}^N, p) \leq \delta\right) = \mathcal{W}^{\otimes N}\left(1_{D(\hat{\mu}^N, p) \leq \delta} \frac{e^{NT(f,\hat{\mu}^N)}}{e^{NT(f,\hat{\mu}^N)}}\right)$$

$$\leq \exp\{-N \inf_{B(p,\delta)} T(f, .)\} \mathcal{W}^{\otimes N}\left(1_{D(\hat{\mu}^N, p) \leq \delta} e^{NT(f,\hat{\mu}^N)}\right)$$

$$\leq \exp\{-N \inf_{B(p,\delta)} T(f, .)\} \mathcal{W}^{\otimes N}\left(e^{NT(f,\hat{\mu}^N)}\right) \quad (4.2.15)$$

$$= \exp\{-N \inf_{\mu \in B(p,\delta)} T(f, \mu)\}$$

Since $\mu \to T(f, \mu)$ is continuous when $f \in \mathcal{C}^{2,1}_b(\mathbb{R} \times [0,1])$, we arrive at

$$\limsup_{\delta \to 0} \limsup_{N \to \infty} \frac{1}{N} \log \mathcal{W}^{\otimes N}\left(\sup_{t \in [0,1]} d(\hat{\mu}^N_t, p_t) \leq \delta\right) \leq -T(f, p).$$

We now optimize over $f$ to obtain a weak large deviation upper bound with rate function

$$S(p) = \sup_{f \in \mathcal{C}^{2,1}_b(\mathbb{R} \times [0,1])} (T^{0,1}(f, p) - \frac{1}{2} < f, f >^{0,1}_p)$$

$$= \sup_{f \in \mathcal{C}^{2,1}_b(\mathbb{R} \times [0,1])} \sup_{\lambda \in \mathbb{R}} (\lambda T^{0,1}(f, p) - \frac{\lambda^2}{2} < f, f >^{0,1}_p)$$

$$= \frac{1}{2} \sup_{f \in \mathcal{C}^{2,1}_b(\mathbb{R} \times [0,1])} \frac{T^{0,1}(f, p)^2}{< f, f >^{0,1}_p} \quad (4.2.16)$$

From the last formula, one sees that any $p$ such that $S(p) < \infty$ is such that $f \to T_f(p)$ is a linear map which is continuous with respect to the norm $||f||^{0,1}_p = (< f, f >^{0,1}_p)^{\frac{1}{2}}$. Hence, Riesz's theorem asserts that there exists a function $k$ verifying (4.2.12, 4.2.13).
• *Large deviation lower bound* The derivation of the large deviation upper bound was thus fairly easy. The lower bound is a bit more sophisticated and relies on the proof of the following points
(a) The solutions to the heat equations with a smooth drift are unique.
(b) The set described by these solutions is dense in $\mathcal{C}([0, 1], \mathcal{P}(\mathbb{R}))$.
(c) The entropy behaves continuously with respect to the approximations by elements of this dense set.
We now describe more precisely these ideas. In the previous section (see (4.2.15)), we have merely obtained the large deviation upper bound from the observation that for all $\nu \in \mathcal{C}([0, 1], \mathcal{P}(\mathbb{R}))$, all $\delta > 0$ and any $f \in \mathcal{C}^{2,1}_b([0, 1], \mathbb{R})$,

$$\mathbb{E}[1_{\hat{\mu}^N \in B(\nu, \delta)} \exp\left(N(T^{0,1}(\hat{\mu}^N, f) - \frac{1}{2} < f, f >^{0,1}_{\hat{\mu}^N})\right)]$$

$$\leq \mathbb{E}[\exp\left(N(T^{0,1}(\hat{\mu}^N, f) - \frac{1}{2} <f,f>_{\hat{\mu}^N}^{0,1})\right)] = 1.$$

To make sure that this upper bound is sharp, we need to check that for any $\nu \in \mathcal{C}([0,1], \mathcal{P}(\mathbb{R}))$ and $\delta > 0$, this inequality is almost an equality for some $k$, i.e there exists $k \in \mathcal{C}_b^{2,1}([0,1], \mathbb{R})$,

$$\liminf_{N \to \infty} \frac{1}{N} \log \frac{\mathbb{E}[1_{\hat{\mu}^N \in B(\nu, \delta)} \exp\left(N(T^{0,1}(\hat{\mu}^N, k) - \frac{1}{2} <k,k>_{\hat{\mu}^N}^{0,1})\right)]}{\mathbb{E}[\exp\left(N(T^{0,1}(\hat{\mu}^N, k) - \frac{1}{2} <k,k>_{\hat{\mu}^N}^{0,1})\right)]} \geq 0.$$

In other words that we can find a $k$ such that the probability that $\hat{\mu}_{\cdot}^N$ belongs to a small neighborhood of $\nu$ under the shifted probability measure

$$\mathbb{P}^{N,k} = \frac{\exp\left(N(T^{0,1}(\hat{\mu}^N, k) - \frac{1}{2} <k,k>_{\hat{\mu}^N}^{0,1})\right)}{\mathbb{E}[\exp\left(N(T^{0,1}(\hat{\mu}^N, k) - \frac{1}{2} <k,k>_{\hat{\mu}^N}^{0,1})\right)]}$$

is not too small. In fact, we shall prove that for good processes $\nu$, we can find $k$ such that this probability goes to one by the following argument. Take $k \in \mathcal{C}_b^{2,1}(\mathbb{R} \times [0,1])$. Under the shifted probability measure $\mathbb{P}^{N,k}$, it is not hard to see that $\hat{\mu}_{\cdot}^N$ is exponentially tight (indeed, for $k \in \mathcal{C}_b^{2,1}(\mathbb{R} \times [0,1])$, the density of $\mathbb{P}^{N,k}$ with respect to $\mathbb{P}$ is uniformly bounded by $e^{C(k)N}$ with a finite constant $C(k)$ so that $\mathbb{P}^{N,k} \circ (\hat{\mu}_{\cdot}^N)^{-1}$ is exponentially tight since $\mathbb{P} \circ (\hat{\mu}_{\cdot}^N)^{-1}$ is). As a consequence, $\hat{\mu}_{\cdot}^N$ is almost surely tight. We let $\mu_{\cdot}$ be a limit point. Now, by Itô's calculus, for any $f \in \mathcal{C}_b^{2,1}(\mathbb{R} \times [0,1])$, any $0 \leq t \leq 1$,

$$T^{0,t}(\hat{\mu}^N, f) = \int_0^t \int \partial_x f_u(x) \partial_x k_u(x) d\hat{\mu}_u^N(x) du + M_t^N(f)$$

with a martingale $(M_t^N(f), t \in [0,1])$ with bracket $(N^{-2} \int_0^t \int (\partial_x f(x))^2 d\hat{\mu}_s^N(x) ds, t \in [0,1])$. Since the bracket of $M_t^N(f)$ goes to zero, the martingale $(M_t^N(f), t \in [0,1])$ goes to zero uniformly almost surely. Hence, any limit point $\mu_{\cdot}$ must satisfy

$$T^{0,1}(\mu, f) = \int_0^1 \int \partial_x f_u(x) \partial_x k_u(x) d\mu_u(x) du \qquad (4.2.17)$$

for any $f \in \mathcal{C}_b^{2,1}(\mathbb{R} \times [0,1])$.
When $(\mu, k)$ satisfies (4.2.17) for all $f \in \mathcal{C}_b^{2,1}(\mathbb{R} \times [0,1])$, we say that $k$ is the field associated with $\mu$.
Therefore, if we can prove that there exists a unique solution $\nu_{\cdot}$ to (4.2.17), we see that $\hat{\mu}_{\cdot}^N$ converges almost surely under $\mathbb{P}^{N,k}$ to this solution. This proves the lower bound at any measure-valued path $\nu_{\cdot}$ which is the unique solution of (4.2.17), namely for any $k \in \mathcal{C}_b^{2,1}(\mathbb{R} \times [0,1])$ such that there exists a unique solution $\nu_k$ to (4.2.17),

$$\liminf_{\delta \to 0} \liminf_{N \to \infty} \frac{1}{N} \log \mathcal{W}^{\otimes N} \left( \sup_{t \in [0,1]} d(\hat{\mu}_t^N, \nu_k) < \delta \right)$$
$$= \liminf_{\delta \to 0} \liminf_{N \to \infty} \frac{1}{N} \log \mathbb{P}^{N,k} \left( 1_{\sup_{t \in [0,1]} d(\hat{\mu}_t^N, \nu_k) < \delta} e^{-NT(k, \hat{\mu}^N)} \right)$$
$$\geq -T(k, \nu_k) + \liminf_{\delta \to 0} \liminf_{N \to \infty} \frac{1}{N} \log \mathbb{P}^{N,k} \left( 1_{\sup_{t \in [0,1]} d(\hat{\mu}_t^N, \nu_k) < \delta} \right)$$
$$\geq -S(\nu_k). \tag{4.2.18}$$

where we used in the second line the continuity of $\mu \to T(\mu, k)$ due to our assumption that $k \in \mathcal{C}_b^{2,1}(\mathbb{R} \times [0,1])$ and the fact that $\mathbb{P}^{N,k} \left( 1_{\sup_{t \in [0,1]} d(\hat{\mu}_t^N, \nu_k) < \delta} \right)$ goes to one in the third line. Hence, the question boils down to uniqueness of the weak solutions of the heat equation with a drift. This problem is not too difficult to solve here and one can see that for instance for fields $k$ which are analytic within a neighborhood of the real line, there is at most one solution to this equation. To generalize (4.2.18) to any $\nu \in \{S < \infty\}$, it is not hard to see that it is enough to find, for any such $\nu$, a sequence $\nu_{k_n}$ for which (4.2.18) holds and such that

$$\lim_{n \to \infty} \nu_{k_n} = \nu, \quad \lim_{n \to \infty} S(\nu_{k_n}) = S(\nu). \tag{4.2.19}$$

Now, observe that $S$ is a convex function so that for any probability measure $p_\epsilon$,

$$S(\mu * p_\epsilon) \leq \int S((. - x)_\# \mu) p_\epsilon(dx) = S(\mu) \tag{4.2.20}$$

where in the last inequality we neglected the condition at the initial time to say that $S((. - x)_\# \mu) = S(\mu)$ for all $x$. Hence, since $S$ is also lower semicontinuous, one sees that $S(\mu * p_\epsilon)$ will converge toward $S(\mu)$ for any $\mu$ with finite entropy $S$. Performing also a regularization with respect to time and taking care of the initial conditions allows to construct a sequence $\nu_n$ with analytic fields satisfying (4.2.19). This point is quite technical but still manageable in this context. Since it will be done quite explicitly in the case we are interested in, we shall not detail it here.

### 4.2.2. Large deviations for the law of the spectral measure of a non-centered large dimensional matrix-valued Brownian motion

To prove a large deviation principle for the law of the spectral measure of Hermitian Brownian motions, the first natural idea would be, following (4.2.11), to prove a large deviation principle for the law of the spectral measure of $\widetilde{L}^N : t \to N^{-1} \sum_{i=1}^N \delta_{\sqrt{N}^{-1} B^i(t)}$, to use Girsanov theorem to show that the law we are considering is absolutely continuous with respect to the law of the independent Brownian motions with a density which only depend on $\widetilde{L}^N$ and conclude by Laplace's method (cf. Chapter 2). However, this approach presents difficulties due to the singularity of the interacting potential, and thus of the density. Here, the techniques developed in [87] will however be very efficient because they only rely on smooth functions of the empirical measure since the empirical measure are taken as distributions so that the interacting potential is smoothed by the test functions (Note however that this strategy would not have worked with more singular potentials). According to (4.2.11), we can in fact follow the very same approach. We here mainly develop the points which are different.

*Itô's calculus*

With the notations of (4.2.7) and (4.2.8), we have

**Theorem 4.10 ([29, Lemma 1.1]).** *1) When $\beta = 2$, for any $N \in \mathbb{N}$, any $f \in \mathcal{C}_b^{2,1}(\mathbb{R} \times [0,1])$ and any $s \in [0,1)$, $\left(S^{s,t}(\hat{\mu}^N, f), s \leq t \leq 1\right)$ is a bounded martingale with quadratic variation*

$$< S^{s,\cdot}(\hat{\mu}^N, f) >_t = \frac{1}{N^2} < f, f >_{\hat{\mu}^N}^{s,t}.$$

*2) When $\beta = 1$ or $4$, for any $N \in \mathbb{N}$, any $f \in \mathcal{C}_b^{2,1}(\mathbb{R} \times [0,1])$ and any $s \in [0,1)$, $\left(S^{s,t}(\hat{\mu}^N, f) + \frac{(-1)^\beta}{2N} \int_s^t \int \partial_x^2 f(y,s) d\hat{\mu}_s^N(x) ds, s \leq t \leq 1\right)$ is a bounded martingale with*
*quadratic variation*

$$< S^{s,\cdot}(\hat{\mu}^N, f) >_t = \frac{2}{\beta N^2} < f, f >_{\hat{\mu}^N}^{s,t}.$$

**Proof :** It is easy to derive this result from Itô's calculus and (4.2.11). Let us however point out how to derive it directly in the case where $f(x) = x^k$ with an integer number $k \in \mathbb{N}$ and $\beta = 2$. Then, for any $(i,j) \in \{1,..,N\}$, Itô 's calculus gives

$$\begin{aligned}
d(\mathbf{H}^{N,2}(t)^k)_{ij} &= \sum_{l=0}^{k-1} \sum_{p,n=1}^{N} (\mathbf{H}^{N,2}(t)^l)_{ip} d(\mathbf{H}^{N,2}(t))_{pn} (\mathbf{H}^{N,2}(t)^{k-l-1})_{nj} \\
&+ \frac{1}{N} \sum_{l+m=0}^{k-2} \sum_{p,n=1}^{N} (\mathbf{H}^{N,2}(t)^l)_{ip} (\mathbf{H}^{N,2}(t)^m)_{nn} (\mathbf{H}^{N,2}(t)^{k-l-m-2})_{pj} dt \\
&= \sum_{l=0}^{k-1} \left(\mathbf{H}^{N,2}(t)^l d\mathbf{H}^{N,2}(t) \mathbf{H}^{N,2}(t)^{k-l-1}\right)_{ij} \\
&+ \sum_{l=0}^{k-2} (k-l-1) \mathrm{tr}(\mathbf{H}^{N,2}(t)^l)(\mathbf{H}^{N,2}(t)^{k-l-2})_{ij} dt
\end{aligned}$$

Let us finally compute the martingale bracket of the normalized trace of the above martingale. We have

$$\langle \int_0^\cdot \sum_{l=0}^{k-1} \text{tr}\left(\mathbf{H}^{N,2}(s)^l d\mathbf{H}^{N,2}(s)\mathbf{H}^{N,2}(s)^{k-l-1}\right)\rangle_t$$
$$= \frac{1}{N^2}k^2 \sum_{ij,mn}\langle \int_0^\cdot (\mathbf{H}^{N,2}(s)^{k-1})_{ij} d\mathbf{H}^{N,2}(s)_{ij}, \int_0^\cdot (\mathbf{H}^{N,2}(s)^{k-1})_{mn} d\mathbf{H}^{N,2}(s)_{mn}\rangle_t$$
$$= \frac{1}{N^2}k^2 \text{tr}\left(\int_0^t (\mathbf{H}^{N,2}(s)^{2(k-1)}) ds\right)$$

Similar computations give the bracket of more general polynomial functions. □

*Remark* 4.11: Observe that if the entries were not Brownian motions but diffusions described for instance as solution of a SDE

$$dx_t = dB_t + U(x_t)dt,$$

then the evolution of the spectral measure of the matrix would not be autonomous anymore. In fact, our strategy is strongly based on the fact that the variations of the spectral measure under small variations of time only depends on the spectral measure, allowing us to construct exponential martingales which are functions of the process of the spectral measure only. It is easy to see that if the entries of the matrix are not Gaussian, the variations of the spectral measures will depend on much more general functions of the entries than those of the spectral measure.

However, this strategy can also be used to study the spectral measure of other Gaussian matrices as emphasized in [29, 60].

From now on, we shall consider the case where $\beta = 2$ and drop the subscript 2 in $\mathbf{H}^{N,2}$, which is slightly easier to write down since there are no error terms in Itô's formula, but everything extends readily to the cases $\beta = 1$ or $4$. The only point to notice is that

$$S_\beta(\mu) = \sup_{\substack{f \in \mathcal{C}^{2,1}_\beta(\mathbb{R}\times[0,1]) \\ 0\leq s\leq t\leq 1}} \{S^{s,t}(\mu, f) - \frac{1}{\beta} <f,f>^{s,t}_\mu\} = \frac{\beta}{2}S_2(\mu)$$

where the last equality is obtained by changing $f$ into $2^{-1}\beta f$.

*Large deviation upper bound*

From the previous Itô's formula, one can deduce by following the ideas of [88] (see Section 4.2.1) a large deviation upper bound for the measure valued process $\hat{\mu}^N_. \in \mathcal{C}([0,1], \mathcal{P}(\mathbb{R}))$. To this end, we shall make the following assumption on the initial condition $D_N$;

(H)
$$C_D := \sup_{N\in\mathbb{N}} \hat{\mu}^N_{D_N}(\log(1+|x|^2)) < \infty,$$

implying that $(\hat{\mu}^N_{D_N}, N \in \mathbb{N})$ is tight. Moreover, $\hat{\mu}^N_{D_N}$ converges weakly, as $N$ goes to infinity, toward a probability measure $\mu_D$.

Then, we shall prove, with the notations of (4.2.7)-(4.2.9), the following

**Theorem 4.12.** *Assume (H). Then*
*(1) $S_{\mu_D}$ is a good rate function on $\mathcal{C}([0,1], \mathcal{P}(\mathbb{R}))$.*
*(2) For any closed set $F$ of $\mathcal{C}([0,1], \mathcal{P}(\mathbb{R}))$,*

$$\limsup_{N\to\infty} \frac{1}{N^2} \log \mathbb{P}\left(\hat{\mu}^N_{\cdot} \in F\right) \leq -\inf_{\nu \in F} S_{\mu_D}(\nu).$$

**Proof :** We first prove that $S_{\mu_D}$ is a good rate function. Then, we show that exponential tightness holds and then obtain a weak large deviation upper bound, these two arguments yielding (2) (cf. Chapter 2).

(a) Let us first observe that $S_{\mu_D}(\nu)$ is also given, when $\nu_0 = \mu_D$, by

$$S_{\mu_D}(\nu) = \frac{1}{2} \sup_{f \in \mathcal{C}^{2,1}_b(\mathbb{R}\times[0,1])} \sup_{0\leq s\leq t\leq 1} \frac{S^{s,t}(\nu, f)^2}{<f,f>_\nu^{s,t}}. \qquad (4.2.21)$$

Consequently, $S_{\mu_D}$ is non negative. Moreover, $S_{\mu_D}$ is obviously lower semi-continuous as a supremum of continuous functions.

Hence, we merely need to check that its level sets are contained in relatively compact sets. For $K$ and $C$ compact subsets of $\mathcal{P}(\mathbb{R})$ and $\mathcal{C}([0,1], \mathbb{R})$, respectively, set

$$\mathcal{K}(K) = \{\nu \in \mathcal{C}([0,1], \mathcal{P}(\mathbb{R})), \nu_t \in K \ \forall t \in [0,1]\}$$

and

$$\mathcal{C}(C, f) = \{\nu \in \mathcal{C}([0,1], \mathcal{P}(\mathbb{R})), (t \to \nu_t(f)) \in C\}.$$

With $(f_n)_{n\in\mathbb{N}}$ a family of bounded continuous functions dense in the set $\mathcal{C}_c(\mathbb{R})$ of compactly supported continuous functions, and $K_M$ and $C_n$ compact subsets of $\mathcal{P}(\mathbb{R})$ and $\mathcal{C}([0,1], \mathbb{R})$, respectively, recall (see [29, Section 2.2]) that the sets

$$\mathcal{K} = \mathcal{K}(K_M) \bigcap \left( \bigcap_{n \in \mathbb{N}} \mathcal{C}(C_n, f_n) \right)$$

are relatively compact subsets of $\mathcal{C}([0,1], \mathcal{P}(\mathbb{R}))$.

Recall now that $K_L = \cap_n \{\nu \in \mathcal{P}(\mathbb{R}) : \nu([-L_n, L_n]^c) \leq n^{-1}\}$ (resp. $C_{\delta,M} = \cap_n \{f \in \mathcal{C}([0,1], \mathbb{R}) : \sup_{|s-t|\leq \delta_n} |f(s) - f(t)| \leq n^{-1}\} \cap \{||f||_\infty \leq M\}$) are compact subsets of $\mathcal{P}(\mathbb{R})$ (resp. $\mathcal{C}([0,1], \mathbb{R})$) for any choice of sequences $(L_n)_{n\in\mathbb{N}} \in (\mathbb{R}^+)^{\mathbb{N}}$, $(\delta_n)_{n\in\mathbb{N}} \in (\mathbb{R}^+)^{\mathbb{N}}$ and positive constant $M$. Thus, following the above description of relatively compact subsets of $\mathcal{C}([0,1], \mathcal{P}(\mathbb{R}))$, to achieve our proof, it is enough to show that, for any $M > 0$,

•1) For any integer $m$, there is a positive real number $L^M_m$ so that for any $\nu \in \{S_{\mu_D} \leq M\}$,

$$\sup_{0\leq s\leq 1} \nu_s(|x| \geq L^M_m) \leq \frac{1}{m} \qquad (4.2.22)$$

proving that $\nu_s \in K_{L^M}$ for all $s \in [0,1]$.

• 2) For any integer $m$ and $f \in \mathcal{C}_b^2(\mathbb{R})$, there exists a positive real number $\delta_m^M$ so that for any $\nu \in \{S_{\mu_D} \leq M\}$,

$$\sup_{|t-s| \leq \delta_m^M} |\nu_t(f) - \nu_s(f)| \leq \frac{1}{m}. \qquad (4.2.23)$$

showing that $s \to \nu_s(f) \in \mathbb{C}_{\delta^M, ||f||_\infty}$.

To prove (4.2.22), we consider, for $\delta > 0$, $f_\delta(x) = \log\left(x^2(1+\delta x^2)^{-1} + 1\right) \in \mathcal{C}_b^{2,1}(\mathbb{R} \times [0,1])$. We observe that

$$C := \sup_{0<\delta\leq 1} ||\partial_x f_\delta||_\infty + \sup_{0<\delta\leq 1} ||\partial_x^2 f_\delta||_\infty$$

is finite and, for $\delta \in (0,1]$,

$$\left|\frac{\partial_x f_\delta(x) - \partial_x f_\delta(y)}{x - y}\right| \leq C.$$

Hence, (4.2.21) implies, by taking $f = f_\delta$ in the supremum, that for any $\delta \in (0,1]$, any $t \in [0,1]$, any $\mu. \in \{S_{\mu_D} \leq M\}$,

$$\mu_t(f_\delta) \leq \mu_0(f_\delta) + 2Ct + 2C\sqrt{Mt}.$$

Consequently, we deduce by the monotone convergence theorem and letting $\delta$ decrease to zero that for any $\mu. \in \{S_{\mu_D} \leq M\}$,

$$\sup_{t \in [0,1]} \mu_t(\log(x^2 + 1)) \leq \mu_D(\log(x^2 + 1)) + 2C(1 + \sqrt{M}).$$

Chebycheff's inequality and hypothesis (H) thus imply that for any $\mu. \in \{S_{\mu_D} \leq M\}$ and any $K \in \mathbb{R}^+$,

$$\sup_{t \in [0,1]} \mu_t(|x| \geq K) \leq \frac{C_D + 2C(1+\sqrt{M})}{\log(K^2 + 1)}$$

which finishes the proof of (4.2.22).

The proof of (4.2.23) again relies on (4.2.21) which implies that for any $f \in \mathcal{C}_b^2(\mathbb{R})$, any $\mu. \in \{S_{\mu_D} \leq M\}$ and any $0 \leq s \leq t \leq 1$,

$$|\mu_t(f) - \mu_s(f)| \leq ||\partial_x^2 f||_\infty |t-s| + 2||\partial_x f||_\infty \sqrt{M}\sqrt{|t-s|}. \qquad (4.2.24)$$

(b) Exponential tightness

**Lemma 4.13.** *For any integer number $L$, there exists a finite integer number $N_0 \in \mathbb{N}$ and a compact set $\mathcal{K}_L$ in $\mathcal{C}([0,1], \mathcal{P}(\mathbb{R}))$ such that $\forall N \geq N_0$,*

$$\mathbb{P}(\hat{\mu}^N \in \mathcal{K}_L^c) \leq \exp\{-LN^2\}.$$

**Proof :** In view of the previous description of the relatively compact subsets of $\mathcal{C}([0,1], \mathcal{P}(\mathbb{R}))$, we need to show that

- a) For every positive real numbers $L$ and $m$, there is an $N_0 \in \mathbb{N}$ and a positive real number $M_{L,m}$ so that $\forall N \geq N_0$

$$\mathbb{P}\Big(\sup_{0 \leq t \leq 1} \hat{\mu}_t^N(|x| \geq M_{L,m}) \geq \frac{1}{m}\Big) \leq \exp(-LN^2)$$

- b) For any $f \in \mathcal{C}_b^2(\mathbb{R})$, for any positive real numbers $L$ and $m$, there exists an $N_0 \in \mathbb{N}$ and a positive real number $\delta_{L,m,f}$ such that $\forall N \geq N_0$

$$\mathbb{P}\Big(\sup_{|t-s| \leq \delta_{L,m,f}} |\hat{\mu}_t^N(f) - \hat{\mu}_s^N(f)| \geq \frac{1}{m}\Big) \leq \exp(-LN^2)$$

The proof is rather classical (it uses Doob's inequality but otherwise is closely related to the proof that $S_{\mu_D}$ is a good rate function); we shall omit it here (see the first section of [29] for details).

(c) Weak large deviation upper bound : We here summarize the main arguments giving the weak large deviation upper bound.

**Lemma 4.14.** *For every process $\nu$ in $\mathcal{C}([0,1], \mathcal{P}(\mathbb{R}))$, if $B_\delta(\nu)$ denotes the open ball with center $\nu$ and radius $\delta$ for the distance $D$, then*

$$\lim_{\delta \to 0} \limsup_{N \to \infty} \frac{1}{N^2} \log \mathbb{P}(\hat{\mu}^N \in B_\delta(\nu)) \leq -S_{\mu_D}(\nu)$$

The arguments are exactly the same as in Section 4.2.1.

*Large deviation lower bound*

We shall prove at the end of this section that

**Lemma 4.15.** *Let*

$$\mathcal{MF}_\infty = \{h \in \mathcal{C}_b^{\infty,1}(\mathbb{R} \times [0,1]) \cap \mathcal{C}([0,1], L^2(\mathbb{R}));$$

$$\exists (C, \epsilon) \in (0, \infty); \sup_{t \in [0,1]} |\hat{h}_t(\lambda)| \leq Ce^{-\epsilon|\lambda|}\}$$

*where $\hat{h}_t$ stands for the Fourier transform of $h_t$. Then, for any field $k$ in $\mathcal{MF}_\infty$, there exists a unique solution $\nu_k$ to*

$$S^{s,t}(f, \nu) = < f, k >_\nu^{s,t} \tag{4.2.25}$$

*for any $f \in \mathcal{C}_b^{2,1}(\mathbb{R} \times [0,1])$. We set $\mathcal{MC}([0,1], \mathcal{P}(\mathbb{R}))$ to be the subset of $\mathcal{C}([0,1], \mathcal{P}(\mathbb{R}))$ consisting of such solutions.*

Note that $h$ belongs to $\mathcal{MF}_\infty$ iff it can be extended analytically to $\{z : |\Im(z)| < \epsilon\}$.

As a consequence of Lemma 4.15, we find that for any open subset $O \in \mathcal{C}([0,1], \mathcal{P}(\mathbb{R}))$, any $\nu \in O \cap \mathcal{MC}([0,1], \mathcal{P}(\mathbb{R}))$, there exists $\delta > 0$ small enough so that

$$\begin{aligned}
\mathbb{P}\left(\hat{\mu}_{\cdot}^{N} \in O\right) &\geq \mathbb{P}\left(d(\hat{\mu}_{\cdot}^{N}, \nu) < \delta\right) \\
&= \mathbb{P}^{N,k}\left(1_{d(\hat{\mu}_{\cdot}^{N}, \nu) < \delta} e^{-N^2(S^{0,1}(\hat{\mu}_{\cdot}^{N}, k) - \frac{1}{2}<k, k>_{\hat{\mu}_{\cdot}^{N}}^{0,1})}\right) \\
&\geq e^{-N^2(S^{0,1}(\nu, k) - \frac{1}{2}<k, k>_{\nu}^{0,1}) - g(\delta)N^2} \mathbb{P}^{N,k}\left(1_{d(\hat{\mu}_{\cdot}^{N}, \nu) < \delta}\right)
\end{aligned}$$

with a function $g$ going to zero at zero. Hence, for any $\nu \in O \cap \mathcal{MC}([0,1], \mathcal{P}(\mathbb{R}))$

$$\liminf_{N \to \infty} \frac{1}{N^2} \log \mathbb{P}\left(\hat{\mu}_{\cdot}^{N} \in O\right) \geq -(S^{0,1}(\nu, k) - \frac{1}{2} <k, k>_{\nu}^{0,1}) = -S_{\mu_D}(\nu)$$

and therefore

$$\liminf_{N \to \infty} \frac{1}{N^2} \log \mathbb{P}\left(\hat{\mu}_{\cdot}^{N} \in O\right) \geq - \inf_{O \cap \mathcal{MC}([0,1], \mathcal{P}(\mathbb{R}))} S_{\mu_D}. \qquad (4.2.26)$$

To complete the lower bound, it is therefore sufficient to prove that for any $\nu \in \mathcal{C}([0,1], \mathcal{P}(\mathbb{R}))$, there exists a sequence $\nu^n \in \mathcal{MC}([0,1], \mathcal{P}(\mathbb{R}))$ such that

$$\lim_{n \to \infty} \nu^n = \nu \text{ and } \lim_{n \to \infty} S_{\mu_D}(\nu^n) = S_{\mu_D}(\nu). \qquad (4.2.27)$$

The rate function $S_{\mu_D}$ is not convex a priori since it is the supremum of **quadratic** functions of the measure-valued path $\nu$ so that there is no reason why it should be reduced by standard convolution as in the classical setting (cf. Section 4.2.1). Thus, it is now unclear how we can construct the sequence (4.2.27). Further, we begin with a degenerate rate function which is infinite unless $\nu_0 = \mu_D$.

To overcome the lack of convexity, we shall remember the origin of the problem; in fact, we have been considering the spectral measure of matrices and should not forget the special features of operators due to the matrices structure. By definition, the differential equation satisfied by a Hermitian Brownian motion should be invariant if we translate the entries, that is translate the Hermitian Brownian motion by a self-adjoint matrix. The natural limiting framework of large random matrices is free probability, and the limiting spectral measure of the sum of a Hermitian Brownian motion and a deterministic self-adjoint matrix converges toward the free convolution of their respective limiting spectral measure. Intuitively, we shall therefore expect (and in fact we will show) that the rate function $S^{0,1}$ decreases by free convolution, generalizing the fact that standard convolution was decreasing the Brownian motion rate function (cf. (4.2.20)). However, because free convolution by a Cauchy law is equal to the standard convolution by a Cauchy law, we shall regularize our laws by convolution by Cauchy laws. Free probability shall be developed in Chapter 6.

Let us here outline the main steps of the proof of (4.2.27):

1. We find that convolution by Cauchy laws $(P^\epsilon)_{\epsilon>0}$ decreases the entropy and prove (this is very technical and proved in [65]) that for any $\nu \in \{S_{\mu_D} < \infty\}$, any given partition $0 = t_1 < t_2 < \ldots < t_n = 1$ with $t_i = (i-1)\Delta$, the measure valued path given, for $t \in [t_k, t_{k+1}[$, by

$$\nu_t^{\epsilon,\Delta} = P^\epsilon * \nu_{t_k} + \frac{(t-t_k)}{\Delta}[P^\epsilon * \nu_{t_{k+1}} - P^\epsilon * \nu_{t_k}],$$

   satisfies

$$\lim_{\epsilon \to 0} \lim_{\Delta \to 0} S^{0,1}(\nu^{\epsilon,\Delta}) = S^{0,1}(\nu), \quad \lim_{\epsilon \to 0} \lim_{\Delta \to 0} \nu^{\epsilon,\Delta} = \nu$$

2. We prove that for $\nu \in \mathcal{A}$,

$$\mathcal{A} = \{\mu \in \mathcal{C}([0,1], \mathcal{P}(\mathbb{R})) : \exists \epsilon > 0; \sup_{t \in [0,1]} \nu_t(|x|^{5+\epsilon}) < \infty\}, \qquad (4.2.28)$$

   $\nu^{\epsilon,\Delta} \in \mathcal{MC}([0,1], \mathcal{P}(\mathbb{R}))$ (actually a slightly weaker form since the field $h^{\epsilon,\Delta}$ may have time discontinuities at the points $\{t_1, ..., t_n\}$, which does not affect the uniqueness statement of Lemma 4.15).
   Here, the choice of Cauchy laws is not innocent; it is a good choice because free convolution by Cauchy laws is, exceptionally, the standard convolution. Hence, it is easier to work with. Moreover, to prove the above result, it is convenient that $\nu_t^{\epsilon,\Delta}$ is absolutely continuous with respect to Lebesgue measure with non vanishing density (see (4.2.31)). But it is not hard to see that for any measure $\mu$ with compact support, $\mu \boxplus p$, if it has a density with respect to Lebesgue measure for many choices of laws $p$, is likely to have holes in its density unless $\int x dp(x) = +\infty$. For these two reasons, the Cauchy law is a wise choice. However, we have to pay for it as we shall see below.

3. Everything looks nice except that we modified the initial condition from $\mu_D$ into $\mu_D * P_\epsilon$, so that in fact $S_{\mu_D}(\nu^{\epsilon,\Delta}) = +\infty$! and moreover, the empirical measure-valued process can not deviate toward processes of the form $\nu^{\epsilon,\Delta}$ even after some time because these processes do not have finite second moment. To overcome this problem, we first note that this result will still give us a large deviation lower bound if we change the initial data of our matrices. Namely, let, for $\epsilon > 0$, $C_\epsilon^N$ be a $N \times N$ diagonal matrix with spectral measure converging toward the Cauchy law $P_\epsilon$ and consider the matrix-valued process

$$\mathbf{X}_t^{N,\epsilon} = \mathbf{U}_N C_\epsilon^N \mathbf{U}_N^* + D_N + \mathbf{H}^N(t)$$

   with $\mathbf{U}_N$ a $N \times N$ unitary measure following the Haar measure $m_2^N$ on $U(N)$. Then, it is well known (see Voiculescu [122]) that the spectral distribution of $\mathbf{U}_N C_\epsilon^N \mathbf{U}_N^* + D_N$ converges toward the free convolution $P_\epsilon \boxplus \mu_D = P_\epsilon * \mu_D$.
   Hence, we can proceed as before to obtain the following large deviation estimates on the law of the spectral measure $\hat{\mu}_t^{N,\epsilon} = \hat{\mu}_{X_t^{N,\epsilon}}^N$

**Corollary 4.16.** *For any $\epsilon > 0$, for any closed subset $F$ of $\mathcal{C}([0,1], \mathcal{P}(\mathbb{R}))$,*

$$\limsup_{N \to \infty} \frac{1}{N^2} \log \mathbb{P}\left(\hat{\mu}_{.}^{N,\epsilon} \in F\right) \leq -\inf\{S_{P_\epsilon * \mu_D}(\nu), \nu \in F\}.$$

*Further, for any open set $O$ of $\mathcal{C}([0,1], \mathcal{P}(\mathbb{R}))$,*

$$\liminf_{N \to \infty} \frac{1}{N^2} \log \mathbb{P}\left(\hat{\mu}_{.}^{N,\epsilon} \in O\right) \geq$$

$$\geq -\inf\{S_{P_\epsilon * \mu_D}(\nu), \nu \in O, \nu = P_\epsilon * \mu, \mu \in \mathcal{A} \cap \{S_{\mu_D} < \infty\}\}.$$

4. To deduce our result for the case $\epsilon = 0$, we proceed by exponential approximation. In fact, we have the following lemma, whose proof is fairly classical and omitted here (cf. [65], proof of Lemma 2.11).

**Lemma 4.17.** *Consider, for $L \in \mathbb{R}^+$, the compact set $K_L$ of $\mathcal{P}(\mathbb{R})$ given by*

$$K_L = \{\mu \in \mathcal{P}(\mathbb{R}); \mu(\log(x^2 + 1)) \leq L\}.$$

*Then, on $\mathcal{K}_\epsilon^N(K_L) := \bigcap_{t \in [0,1]} \left\{\{\hat{\mu}_t^{N,\epsilon} \in K_L\} \cap \{\hat{\mu}_t^N \in K_L\}\right\}$,*

$$D(\hat{\mu}_{.}^{N,\epsilon}, \hat{\mu}_{.}^N) \leq f(N, \epsilon)$$

*where*

$$\limsup_{\epsilon \to 0} \limsup_{N \to \infty} f(N, \epsilon) = 0.$$

We then can prove that

**Theorem 4.18.** *Assume that $\hat{\mu}_{D_N}^N$ converges toward $\mu_D$ while*

$$\sup_{N \in \mathbb{N}} \hat{\mu}_{D_N}^N(x^2) < \infty.$$

*Then, for any $\mu_{.} \in \mathcal{A}$*

$$\lim_{\delta \to 0} \liminf_{N \to \infty} \frac{1}{N^2} \log \mathbb{P}\left(D(\hat{\mu}_{.}^N, \mu_{.}) \leq \delta\right) \geq -S_{\mu_D}(\mu_{.})$$

*so that for any open subset $O \in \mathcal{C}([0,1], \mathcal{P}(\mathcal{P}(I\!R)))$,*

$$\liminf_{N \to \infty} \frac{1}{N^2} \log \mathbb{P}\left(\hat{\mu}_{.}^N \in O\right) \geq -\inf_{O \cap \mathcal{A}} S_{\mu_D}$$

**Proof of Theorem 4.18 :** Following Lemma 4.13, we deduce that for any $M \in \mathbb{R}^+$, we can find $L_M \in \mathbb{R}^+$ such that for any $L \geq L_M$,

$$\sup_{0 \leq \epsilon \leq 1} \mathbb{P}(\mathcal{K}_\epsilon^N(K_L)^c) \leq e^{-MN^2}. \tag{4.2.29}$$

Fix $M > S_{\mu_D}(\mu) + 1$ and $L \geq L_M$. Let $\delta > 0$ be given. Next, observe that $P_\epsilon * \mu.$ converges weakly toward $\mu.$ as $\epsilon$ goes to zero and choose consequently $\epsilon$ small enough so that $D(P_\epsilon * \mu., \mu.) < \frac{\delta}{3}$. Then, write

$$\begin{aligned}
\mathbb{P}\left(\hat{\mu}_.^N \in B(\mu., \delta)\right) &\geq \mathbb{P}\left(D(\hat{\mu}_.^N, \mu.) < \frac{\delta}{3}, \hat{\mu}_.^{N,\epsilon} \in B(P_\epsilon * \mu., \frac{\delta}{3}), \mathcal{K}_\epsilon^N(K_L)\right) \\
&\geq \mathbb{P}\left(\hat{\mu}_.^{N,\epsilon} \in B(P_\epsilon * \mu., \frac{\delta}{3})\right) - \mathbb{P}(\mathcal{K}_\epsilon^N(K_L)^c) \\
&\quad - \mathbb{P}\left(D(\hat{\mu}_.^{N,\epsilon}, \hat{\mu}_.^N) \geq \frac{\delta}{3}, \mathcal{K}_\epsilon^N(K_L)\right) = I - II - III.
\end{aligned}$$

(4.2.29) implies, up to terms of smaller order, that

$$II \leq e^{-N^2(S_{\mu_D}(\mu)+1)}.$$

Lemma 4.17 shows that $III = 0$ for $\epsilon$ small enough and $N$ large, while Corollary 4.16 imply that for any $\eta > 0$, $N$ large and $\epsilon > 0$

$$I \geq e^{-N^2 S_{P_\epsilon * \mu_D}(P_\epsilon * \mu) - N^2 \eta} \geq e^{-N^2 S_{\mu_D}(\mu) - N^2 \eta}.$$

Theorem 4.18 is proved. □

5. To complete the lower bound, we need to prove that for any $\nu \in \{S_{\mu_D} < \infty\}$, there exists a sequence of $\nu_n \in \mathcal{A}$ such that

$$\lim_{n \to \infty} S_{\mu_D}(\nu_n) = S_{\mu_D}(\nu), \quad \lim_{n \to \infty} \nu_n = \nu. \quad (4.2.30)$$

This is done in [66] by the following approximation : we let,

$$\begin{aligned}
\nu_t^{\eta,\epsilon}(dx) &= \mu_D \boxplus \sigma_t, \quad \text{if } t \leq \eta \\
&= [\mu_D \boxplus \sigma_\eta]^{(t-\eta)}, \quad \text{if } \eta \leq t \leq \eta + \epsilon \\
&= [\sigma_\eta \boxplus \nu_{t-\eta-\epsilon}]^\epsilon, \quad \text{if } t \geq \eta + \epsilon
\end{aligned}$$

where, for $\mu \in \mathcal{P}(\mathbb{R})$, $[\mu]^\eta$ is the probability measure given for $f \in \mathcal{C}_b(\mathbb{R})$ by

$$[\mu]^\eta(f) := \int f((1 + \eta x^2)^{-a} x) d\mu(x)$$

for some $a < \frac{1}{2}$. Then, we show that we can choose $(\eta_n, \epsilon_n)_{n \in \mathbb{N}}$ so that $\nu_n = \nu_t^{\eta_n, \epsilon_n}$ satisfies (4.2.30). Moreover $\nu_n \in \mathcal{A}$ for $5(1 - 2a) < 2$ because $\sup_{t \in [0,1]} \nu_t(x^2) < \infty$ as $S_{\mu_D}(\nu) < \infty$, $\sigma_t$ is compactly supported and we assumed $\mu_D(|x|^{5+\epsilon}) < \infty$.

6. To finish the proof, we need to complete lemmas 4.15, 4.17, and the points 1) and 2) (see (4.2.22),(4.2.23)) of our program. We prove below Lemma 4.15. We provide in Chapter 6 (see sections 6.6 and 6.7) part of the proofs of 1), 2).

In fact we show in Chapter 6 that $S^{0,1}(\nu^{\epsilon,\Delta})$ converges toward $S^{0,1}(\nu)$. The fact that the field $h^{\epsilon,\Delta}$ associated with $\nu^{\epsilon,\Delta}$ satisfies the necessary conditions so that $\nu^{\epsilon,\Delta} \in \mathcal{MC}([0,1], \mathcal{P}(\mathbb{R}))$ is proved in [65]. We shall not detail it here but let us just point out the basic idea which is to observe that (4.2.25) is equivalent to write that, if $\mu_t(dx) = \rho_t(x)dx$ with a smooth density $\rho_.$,

$$\partial_t \rho_t(x) = -\partial_x(\rho_t(x)H\rho_t(x) + \rho_t(x)\partial_x k_t(x)),$$

where $H\nu$ is the Hilbert transform

$$H\nu(x) = PV \int \frac{1}{x-y} d\nu(y) = \lim_{\epsilon \to 0} \int \frac{(x-y)}{(x-y)^2 + \epsilon^2} d\nu(y).$$

In other words,

$$\partial_x k_t(x) = \frac{\int_x^\infty \partial_t \rho_t(y) dy}{\rho_t(x)} - H\rho_t(x). \qquad (4.2.31)$$

Hence, we see that $\partial_x k_.$ is smooth as soon as $\rho$ is, that its Hilbert transform behaves well and that $\rho$ does not vanish. To study the Fourier transform of $\partial_x k_t$, we need it to belong to $L^1(dx)$, which we can only show when the original process $\nu$ possess at least finite fifth moment. More details are given in [65].

**Proof of Lemma 4.15 :** Following [29], we take $f(x,t) := e^{i\lambda x}$ in (4.2.25) and denote by $\mathcal{L}_t(\lambda) = \int e^{i\lambda x} d\nu_t(x)$ the Fourier transform of $\nu_t$. $\nu \in \mathcal{MC}([0,1], \mathcal{P}(\mathbb{R}))$ implies that if $k$ is the field associated with $\nu$, $|\hat{k}_t(\lambda)| \leq Ce^{-\epsilon|\lambda|}$ with a given $\epsilon > 0$. Then, we find that for $t \in [0 = t_1, t_2]$,

$$\mathcal{L}_t(\lambda) = \mathcal{L}_0(\lambda) - \frac{\lambda^2}{2} \int_0^t \int_0^1 \mathcal{L}_s(\alpha\lambda)\mathcal{L}_s((1-\alpha)\lambda) d\alpha ds + i\lambda \int_0^t \int \mathcal{L}_s(\lambda+\lambda')\hat{k}(\lambda',s) d\lambda' ds. \qquad (4.2.32)$$

Multiplying both sides of this equality by $e^{-\frac{\epsilon}{4}|\lambda|}$ gives, with $\mathcal{L}_t^\epsilon(\lambda) = e^{-\frac{\epsilon}{4}|\lambda|}\mathcal{L}_t(\lambda)$,

$$\begin{aligned}\mathcal{L}_t^\epsilon(\lambda) &= \mathcal{L}_0^\epsilon(\lambda) - \frac{\lambda^2}{2} \int_0^t \int_0^1 \mathcal{L}_s^\epsilon(\alpha\lambda)\mathcal{L}_s^\epsilon((1-\alpha)\lambda) d\alpha ds \\ &+ i\lambda \int_0^t \int \mathcal{L}_s^\epsilon(\lambda+\lambda')e^{\frac{\epsilon}{4}|\lambda+\lambda'|-\frac{\epsilon}{4}|\lambda|}\hat{k}(\lambda',s) d\lambda' ds. \qquad (4.2.33)\end{aligned}$$

Therefore, if $\nu, \widetilde{\nu}'$ are two solutions with Fourier transforms $\mathcal{L}$ and $\widetilde{\mathcal{L}}$ respectively and if we set $\Delta_t^\epsilon(\lambda) = |\mathcal{L}_t^\epsilon(\lambda) - \widetilde{\mathcal{L}}_t^\epsilon(\lambda)|$, we deduce from (4.2.33) (see [29], proof

of Lemma 2.6, for details) that

$$\begin{aligned}
\Delta_t^\epsilon(\lambda) &\leq \lambda^2 \int_0^t \int \Delta_s^\epsilon(\alpha\lambda) e^{-\frac{1}{4}(1-\alpha)\epsilon\lambda} d\alpha ds \\
&\quad + C|\lambda| \int_0^t \int \Delta_s^\epsilon(\lambda+\lambda') e^{\frac{\epsilon}{4}|\lambda+\lambda'|-\frac{\epsilon}{4}|\lambda|-\epsilon|\lambda'|} d\lambda' ds \\
&\leq (C+\frac{4}{\epsilon})|\lambda| \int_0^t \sup_{|\lambda'|\leq|\lambda|} \Delta_s^\epsilon(\lambda')ds + 3t|\lambda|e^{-\frac{\epsilon}{4}|\lambda|}
\end{aligned}$$

where we used that $\Delta_t^\epsilon(\lambda) \leq 2e^{-\frac{\epsilon}{4}|\lambda|}$. Considering $\bar{\Delta}_t^\epsilon(R) = \sup_{|\lambda'|\leq R} \Delta_s^\epsilon(\lambda')$, we therefore obtain

$$\bar{\Delta}_t^\epsilon(R) \leq (C+\frac{4}{\epsilon})R \int_0^t \bar{\Delta}_s^\epsilon(R)ds + 3tRe^{-\frac{\epsilon}{4}R}.$$

By Gronwall's lemma, we deduce that

$$\bar{\Delta}_t^\epsilon(R) \leq 3Re^{-\frac{\epsilon}{4}R}e^{(C+\frac{4}{\epsilon})Rt}$$

and thus that $\bar{\Delta}_t^\epsilon(\infty) = 0$ for $t < \tau \equiv \frac{\epsilon^2}{4(4+\epsilon C)}$. By induction over the time, we conclude that $\bar{\Delta}_t^\epsilon(\infty) = 0$ for any time $t \leq t_2$, and then any time $t \leq 1$, and therefore that $\nu = \tilde{\nu}$. □

### 4.3. Discussion and open problems

We have seen in this section how the non-intersecting paths description (or equivalently the Dyson equations (4.2.11)) of spherical integrals can be used to study their first order asymptotic behaviour and understand where the integral concentrates.

A related question is to study the law of the spectral measure of the random matrix $M$ with law

$$d\mu_N(M) = \frac{1}{Z_N} e^{-N\text{Tr}(V(M))+N\text{Tr}(AM)} dM.$$

In the case where $V(M) = \frac{1}{2}M^2$, $M = W + A$ with $W$ a Wigner matrix, the asymptotic distribution of the eigenvalues are given by the free convolution of the semi-circular distribution and the limiting spectral distribution of $A$. A more detailed study based on Riemann-Hilbert techniques gives the limiting eigenvalue distribution correlations when $\hat{\mu}_A^N = \alpha_N \delta_a + (1-\alpha_N)\delta_{-a}$ (cf. [19, 3]). It is natural to wonder whether such a result could be derived from the Brownian paths description of the matrix. Our result allows (as a mild generalization of the next chapter) to describe the limiting spectral measure for general potentials $V$. However, the study of correlations requires different sets of techniques.

It would be very interesting to understand the relation between this limit and the expansion obtained by B. Collins [34] who proved that

$$\lim_{N\to\infty} \partial_\lambda^p \frac{1}{N^2} \log I_N^{(2)}(\lambda E_N, D_N)|_{\lambda=0} = a_p(\mu_D, \mu_E)$$

with some complicated functions $a_p(\mu_D, \mu_E)$ of interest. The physicists answer to such questions is that

$$a_p(\mu_D, \mu_E) = \partial_\lambda^p \lim_{N\to\infty} \frac{1}{N^2} \log I_N^{(2)}(\lambda E_N, D_N)|_{\lambda=0} := \widetilde{a}_p(\mu_D, \mu_E)$$

which would validate the whole approach to use the asymptotics of $I_N^{(\beta)}$ to compute and study the $a_p(\mu_D, \mu_E)$'s. However, it is not yet known whether such an interchange of limit and derivation is rigorous. A related topic is to understand whether the asymptotics we obtained extends to non Hermitian matrices. This is not true in general following a counterexample of S. Zelditch [133] but still could hold when the spectral norm of the matrices is small enough. In the case where one matrix has rank one, I have shown with M. Maida [63] that such an analytic extension was true.

Other models such as domino tilings can be represented by non-intersecting paths (cf. [77] for instance). It is rather tempting to hope that similar techniques could be used in this setting and give a second approach to [84]. This however seems slightly more involved because time and space are discrete and have the same scaling (making the approximation by Brownian motions irrelevant), so that the whole machinery borrowed from hydrodynamics technology does not seem to be well adapted.

It would be as well very interesting to obtain second order corrections terms for the spherical integrals, problem related with the enumeration of maps with genus $g \geq 1$ as we shall see in the next section.

Finally, it would also be nice to get a better understanding of the limiting value of the spherical integrals, namely of $J_\beta(\mu_D, \mu_E)$. We shall give some elements in this direction in the next section but wish to emphasize already that this quantity remains rather mysterious. We have not yet been able for instance to obtain a simple formula in the case of Bernouilli measures (which are somewhat degenerate cases in this context).

## Chapter 5

# Matrix models and enumeration of maps

It appears since the work of 't Hooft that matrix integrals can be seen, via Feynman diagrams expansion, as generating functions for enumerating maps (or triangulated surfaces). We refer here to the very nice survey of A. Zvonkin's [136]. One matrix integrals are used to enumerate maps with a given genus and given vertices degrees distribution whereas several matrices integrals can be used to consider the case where the vertices can additionally be colored (i.e. can take different states).

Matrix integrals are usually of the form

$$Z_N(P) = \int e^{-N\text{Tr}(P(A_1^N, \cdots, A_d^N))} dA_1^N \cdots dA_d^N$$

with some polynomial function $P$ of $d$-non-commutative variables and the Lebesgue measure $dA$ on some well chosen ensemble of $N \times N$ matrices such as the set $\mathcal{H}_N$ (resp. $\mathcal{S}_N$, resp. $\mathcal{S}ymp_N$) of $N \times N$ Hermitian (resp. symmetric, resp. symplectic) matrices.

We shall describe in the next section how such integrals are related with the enumeration of maps. Then, we shall apply the results of the previous section to compute the first order asymptotics of such integrals in some cases where the polynomial function $P$ has a quadratic interaction.

**5.1. Relation with the enumeration of maps**

Following A. Zvonkin's [136], let us consider the case where $d = 1$, $P(x) = tx^4 + \frac{1}{4}x^2$ and integration holds over $\mathcal{H}_N$. Then, it is argued that
<u>Conjecture</u> 5.1:[see A. Zvonkin's [136]] For any $N \in \mathbb{N}$, any $t \in \mathbb{R}$

$$\frac{1}{N^2} \log \int e^{-N\text{Tr}(tA^4 + \frac{1}{4}A^2)} dA = \sum_{n \geq 0} \sum_{g \geq 0} \frac{(-t)^n}{n! N^{2g}} C(n, g)$$

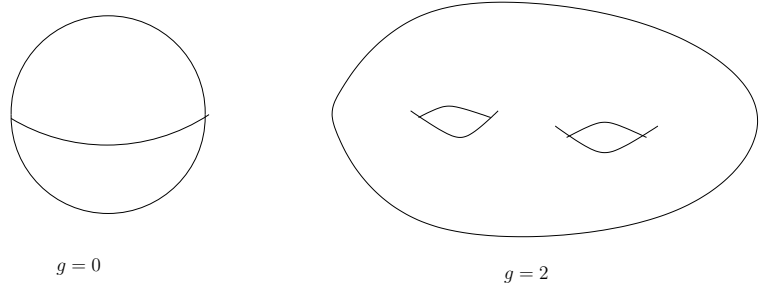

Fig 5.1. *Surfaces of genus $g = 0, 2$*

with
$$C(n, g) = \sharp\{\text{maps of genus } g \text{ with } n \text{ vertices of degree } 4\}.$$

Here, a map of genus $g$ is an oriented connected graph drawn on a surface of genus $g$ modulo equivalent classes.

To be more precise, a surface is a compact oriented two-dimensional manifold without boundary. The surfaces are classified according to their genera, which is characterized by a genus (the number of 'handles', see figure 5.1). There is only one surface with a given genus up to homeomorphism. Following [136], Definition 4.1, a map is then a graph which is 'drawn' (or embedded into) a surface in such a way that the edges do not intersect and if we cut the surface along the edges, we get a disjoint union of sets which are homeomorphic to an open disk (these sets are the faces of the map). Note that in Conjecture 5.1, the maps are counted up to homeomorphisms (that is modulo equivalent classes).

The **formal** proof of Conjecture 5.1 goes as follows. One begins by expanding all the non quadratic terms in the exponential

$$Z_N(tA^4) = \int e^{-N\mathrm{Tr}(tA^4 + \frac{1}{4}A^2)} dA = Z_N(0) \sum_{n \geq 0} \frac{(-tN)^n}{n!} \mu_N\left((tr(A^4))^n\right)$$

with $\mu_N$ the Gaussian law

$$\mu_N(dA) = Z_N(0)^{-1} e^{-\frac{N}{4}\mathrm{Tr}(A^2)} dA$$

that is the law of Wigner's Hermitian matrices given by

$$A_{ji} = \bar{A}_{ij}, \quad A_{ij} = \mathcal{N}(0, N^{-1}), 1 \leq i \leq j \leq N.$$

Writing

$$\mu_N\left((tr(A^4))^n\right) = \sum_{\substack{i_1^k, i_2^k, i_3^k, i_4^k = 1 \\ 1 \leq k \leq n}} \mu_N\left(\prod_{k=1}^n A_{i_1^k i_2^k} A_{i_2^k i_3^k} A_{i_3^k i_4^k} A_{i_4^k i_1^k}\right)$$

and using Wick's formula to compute the expectation over the Gaussian variables, gives the result.

One can alternatively use the graphical representation introduced by Feynman to compute such expectations. It goes as shown on figure 5.2.

The last equivalence in figure 5.2 results from the observation that $\mu_N(A_{ij}) = 0$ for all $i, j$ and that $\mu_N(A_{ij}A_{kl}) = \delta_{ij=lk}N^{-1}$, so that each end of the crosses has to be connected with another one. As a consequence from this construction, one can see that only oriented fat graphs will contribute to the sum. Moreover, since on each face of the graph, the indices are constant, we see that each given graph will appear $N^{\sharp\text{faces}}$ times. Hence, if we denote

$$G(n, F) = \{\text{oriented graphs with n vertices with degree 4 and F faces}\}$$

we obtain

$$Z_N(tA^4) = \sum_{n\geq 0}\sum_{F\geq 0} \frac{(-t)^n}{N^{n-F}} \sharp G(n, F).$$

Taking logarithm, it is well known that we get only contributions from connected graphs:

$$\log Z_N(tA^4) = \sum_{n\geq 0}\sum_{F\geq 0} \frac{(-t)^n}{N^{n-F}} \sharp\{G(n, F) \cap \text{ connected graphs }\}.$$

Finally, since the genus of a map is related with its number of faces by $n - F = 2(g - 1)$, the result is proved. When considering several matrices, we see that we have additionally to decide at each vertex which matrix contributed, which corresponds to give different states to the vertices.

Of course, this derivation is formal and it is not clear at all that such an expansion of the free energy exists. Indeed, the series here might well not be summable (this is clear when $t < 0$ since the integral diverges). In the one matrix case, this result was proved very recently by J. Mc Laughlin et N. Ercolani [46] who have shown that such an expansion is valid by mean of Riemann-Hilbert problems techniques, under natural assumptions over the potential. The first step in order to use Riemann-Hilbert problems techniques, is to understand the limiting behavior of the spectral measures of the matrices following the corresponding Gibbs measure

$$\mu_N^P(d\mathbf{A}_1, \cdots, d\mathbf{A}_d) = \frac{1}{Z_N(P)} e^{-N\text{Tr}(P(\mathbf{A}_1,\cdots,\mathbf{A}_d))} d\mathbf{A}_1 \cdots d\mathbf{A}_d.$$

This is clearly needed to localize the integral. The understanding of such asymptotics is the subject of the next section.

Let us remark before embarking in this line of attack that a more direct combinatorial strategy can be developed. For instance, G. Scheaffer and M. Bousquet Melou [22] studied the Ising model on planar random graph as a generating function for the enumeration of colored planar maps, generalizing Tutte's approach. The results are then more explicitly described by an algebraic

We can represent $A_{i_1i_2}A_{i_2i_3}A_{i_3i_4}A_{i_4i_1}$ by an oriented cross

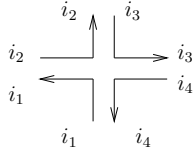

Therefore, $\Lambda = \prod_{i=1}^n A_{i_1^k i_2^k} A_{i_2^k i_3^k} A_{i_3^k i_4^k} A_{i_4^k i_1^k}$ is represented by $n$ such crosses

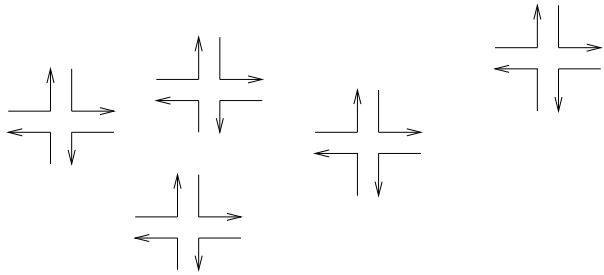

Since $\mu_N(A_{ij}) = 0$, $\mu_N(A_{ij}A_{kl}) = N^{-1}1_{ij=lk}$, the expectation of $\Lambda$ will be null except if the indices form an oriented diagram equation for this generating function. However, matrix models approach is more general a priori since it allows to consider many different models and arbitrary genus, eventhough the mathematical understanding of these points is still far from being achieved.

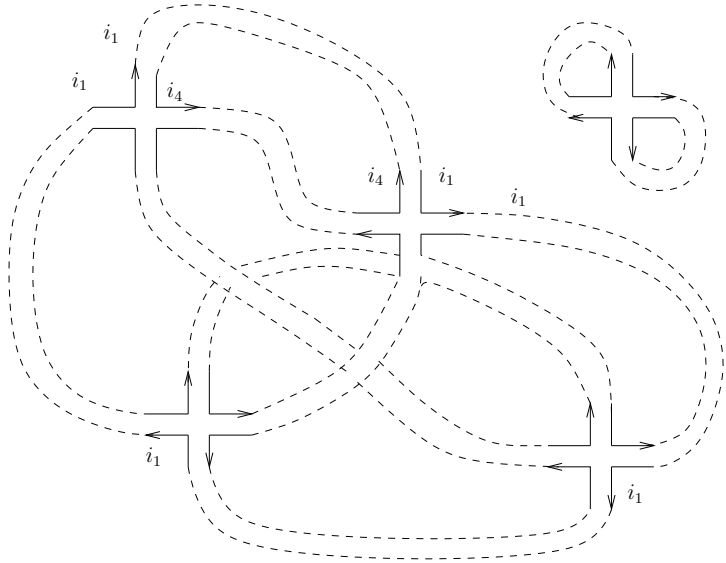

Fig 5.2. *Computing Gaussian integrals*

Finally, let us notice that the relation between the enumeration of maps and the matrix models should a priori be true only when the weights of non quadratic Gaussian terms are small (the expansion being in fact a formal analytic expansion around the origin) but that matrix integrals are of interest in other regimes for instance in free probability.

### 5.2. Asymptotics of some matrix integrals

We would like to consider integrals of more than one matrix. The simplest interaction that one can think of is the quadratic one. Such an interaction describes already several classical models in random matrix theory; We refer here to the works of M. Mehta, A. Matytsin, A. Migdal, V. Kazakov, P. Zinn Justin and B. Eynard for instance.

- The random Ising model on random graphs is described by the Gibbs measure

$$\mu_{Ising}^N(dA, dB) = \frac{1}{Z_{Ising}^N} e^{N\text{Tr}(AB) - N\text{Tr}(P_1(A)) - N\text{Tr}(P_2(B))} dA dB$$

with $Z_{Ising}^N$ the partition function

$$Z_{Ising}^N = \int e^{N\text{Tr}(AB) - N\text{Tr}(P_1(A)) - N\text{Tr}(P_2(B))} dA dB$$

and two polynomial functions $P_1, P_2$. The limiting free energy for this model was calculated by M. Mehta [93] in the case $P_1(x) = P_2(x) = x^2 + gx^4$ and integration holds over $\mathcal{H}_N$. However, the limiting spectral measures of $A$ and $B$ under $\mu_{Ising}^N$ were not considered in that paper. A discussion about this problem can be found in P. Zinn Justin [135].

- One can also define the $q-1$ Potts model on random graphs described by the Gibbs measure

$$\mu_{Potts}^N(d\mathbf{A}_1, ..., d\mathbf{A}_q) =$$

$$= \frac{1}{Z_{Potts}^N} \prod_{i=2}^{q} e^{N\text{Tr}(\mathbf{A}_1\mathbf{A}_i) - N\text{Tr}(P_i(\mathbf{A}_i))} d\mathbf{A}_i e^{-N\text{Tr}(P_1(\mathbf{A}_1))} d\mathbf{A}_1.$$

The limiting spectral measures of $(\mathbf{A}_1, \cdots, \mathbf{A}_q)$ are discussed in [135] when $P_i = gx^3 - x^2$ (!).

- As a straightforward generalization, one can consider matrices coupled by a chain following S. Chadha, G. Mahoux and M. Mehta [94] given by

$$\mu_{chain}^N(d\mathbf{A}_1, ..., d\mathbf{A}_q) =$$

$$= \frac{1}{Z_{chain}^N} \prod_{i=2}^{q} e^{N\text{Tr}(\mathbf{A}_{i-1}\mathbf{A}_i) - N\text{Tr}(P_i(\mathbf{A}_i))} d\mathbf{A}_i e^{-N\text{Tr}(P_1(\mathbf{A}_1))} d\mathbf{A}_1.$$

$q$ can eventually go to infinity as in [92].

The first order asymptotics of these models can be studied thanks to the control of spherical integrals obtained in the last chapter.

**Theorem 5.2.** *Assume that $P_i(x) \geq c_i x^4 + d_i$ with $c_i > 0$ and some finite constants $d_i$. Hereafter, $\beta = 1$ (resp. $\beta = 2$, resp. $\beta = 4$) when $dA$ denotes the Lebesgue measure on $\mathcal{S}_N$ (resp. $\mathcal{H}_N$, resp. $\mathcal{H}_N$ with $N$ even). Then, with $c = \inf_{\nu \in \mathcal{P}(\mathbb{R})} I_\beta(\nu)$,*

$$\begin{aligned}
F_{Ising} &= \lim_{N \to \infty} \frac{1}{N^2} \log Z_{Ising}^N \\
&= -\inf\{\mu(P) + \nu(Q) - I^{(\beta)}(\mu, \nu) - \frac{\beta}{2}\Sigma(\mu) - \frac{\beta}{2}\Sigma(\nu)\} - 2c \quad (5.2.1)
\end{aligned}$$

$$\begin{aligned}
F_{Potts} &= \lim_{N \to \infty} \frac{1}{N^2} \log Z_{Potts}^N \\
&= -\inf\{\sum_{i=1}^{q} \mu_i(P_i) - \sum_{i=2}^{q} I^{(\beta)}(\mu_1, \mu_i) - \frac{\beta}{2}\sum_{i=1}^{q}\Sigma(\mu_i)\} - qc \quad (5.2.2)
\end{aligned}$$

$$\begin{aligned}
F_{chain} &= \lim_{N \to \infty} \frac{1}{N^2} \log Z_{chain}^N \\
&= -\inf\{\sum_{i=1}^{q} \mu_i(P_i) - \sum_{i=2}^{q} I^{(\beta)}(\mu_{i-1}, \mu_i) - \frac{\beta}{2}\sum_{i=1}^{q}\Sigma(\mu_i)\} - qc \quad (5.2.3)
\end{aligned}$$

<u>Remark</u> 5.3: The above theorem actually extends to polynomial functions going to infinity like $x^2$. However, the case of quadratic polynomials is trivial since it boils down to the Gaussian case and therefore the next interesting case is quartic polynomial as above. Moreover, Theorem 5.4 fails in the case where $P, Q$ go to infinity only like $x^2$. However, all our proofs would extends easily for any continuous functions $P_i's$ such that $P_i(x) \geq a|x|^{2+\epsilon} + b$ with some $a > 0$ and $\epsilon > 0$. In particular, we do not need any analyticity assumptions which can be required for instance to obtain the so-called Master loop equations (see Eynard and als. [50, 49]).

**Proof of Theorem 5.2 :** It is enough to notice that, when diagonalizing the matrices $\mathbf{A}_i$'s, the interaction is expressed in terms of spherical integrals by (3.1.1). Laplace's (or saddle point) method then gives the result (up to the boundedness of the matrices $A_i$'s in the spherical integrals, which can be obtained by approximation). We shall not detail it here and refer the reader to [61] . □

We shall then study the variational problems for the above energies; indeed, by standard large deviation considerations, it is clear that the spectral measures of the matrices $(\mathbf{A}_i)_{1\leq i\leq d}$ will concentrate on the set of the minimizers defining the free energies, and in particular converge to these minimizers when they are unique. We prove the following for the Ising model.

**Theorem 5.4.** *Assume $P_1(x) \geq ax^4 + b, P_2(x) \geq ax^4 + b$ for some positive constant $a$. Then*

*0) The infimum in $F_{Ising}$ is achieved at a unique couple $(\mu_A, \mu_B)$ of probability measures.*

*1) $(\hat{\mu}_A^N, \hat{\mu}_B^N)$ converges almost surely toward $(\mu_A, \mu_B)$.*

*2) $(\mu_A, \mu_B)$ are compactly supported with finite non-commutative entropy*

$$\Sigma(\mu) = \int\int \log|x-y|d\mu(x)d\mu(y).$$

*3) There exists a couple $(\rho^{A\to B}, u^{A\to B})$ of measurable functions on $\mathbb{R}\times(0,1)$ such that $\rho_t^{A\to B}(x)dx$ is a probability measure on $\mathbb{R}$ for all $t \in (0,1)$ and $(\mu_A, \mu_B, \rho^{A\to B}, u^{A\to B})$ are characterized uniquely as the minimizer of a strictly convex function under a linear constraint.*

*In particular, $(\rho^{A\to B}, u^{A\to B})$ are solution of the Euler equation for isentropic flow with negative pressure $p(\rho) = -\frac{\pi^2}{3}\rho^3$ such that, for all $(x,t)$ in the interior of $\Omega = \{(x,t) \in \mathbb{R}\times[0,1]; \rho_t^{A\to B}(x) \neq 0\}$,*

$$\begin{cases} \partial_t \rho_t^{A\to B} + \partial_x(\rho_t^{A\to B} u_t^{A\to B}) = 0 \\ \partial_t(\rho_t^{A\to B} u_t^{A\to B}) + \partial_x(\rho_t^{A\to B}(u_t^{A\to B})^2 - \frac{\pi^2}{3}(\rho_t^{A\to B})^3) = 0 \end{cases} \quad (5.2.4)$$

*with the probability measure $\rho_t^{A\to B}(x)dx$ weakly converging toward $\mu_A(dx)$ (resp. $\mu_B(dx)$) as $t$ goes to zero (resp. one). Moreover, we have*

$$P'(x) - x - \frac{\beta}{2}u_0^{A\to B}(x) - \frac{\beta}{2}H\mu_A(x) = 0 \quad \mu_A - a.s$$

*and* $\quad Q'(x) - x + \frac{\beta}{2}u_1^{A\to B}(x) - \frac{\beta}{2}H\mu_B(x) = 0 \quad \mu_B - a.s.$

For the other models, uniqueness of the minimizers is not always clear. For instance, we obtain uniqueness of the minimizers for the $q$-Potts models only for $q \leq 2$ whereas it is also expected for $q = 3$. For the description of these minimizers, I refer the reader to [61].

To prove Theorem 5.4, we shall study more carefully the entropy

$$J_\beta(\mu, \mu_D) = \frac{\beta}{2}\inf\{S_{\mu_D}(\nu), \nu_1 = \mu\}$$

and in particular understand where the infimum is taken. The main observation is that by (4.2.21), for any $f \in \mathcal{C}_b^{2,1}(\mathbb{R}\times[0,1])$,

$$S^{0,1}(\nu, f)^2 \leq 2S^{0,1}(\nu) <f,f>_{0,1}^\nu$$

so that the linear form $f \to S^{0,1}(\nu, f)$ is bounded and Riesz's theorem asserts that there exists $k \in H^1_\nu = \overline{\mathcal{C}^{2,1}_b(\mathbb{R} \times [0,1])}^{<\cdot,\cdot>^\nu_{0,1}}$ such that for any $f \in \mathcal{C}^{2,1}_b(\mathbb{R} \times [0,1])$

$$S^{0,1}(\nu, f) = <f, k>^\nu_{0,1}. \tag{5.2.5}$$

We then say that $(\nu, k)$ satisfies (5.2.5). Then

$$S^{0,1}(\nu) = \frac{1}{2} <k, k>^{0,1}_\nu = \frac{1}{2} \int_0^1 \int (\partial_x k_t)^2 d\nu_t(x) dt.$$

**Property 5.5.** *Let $\mu_0 \in \{\mu \in \mathcal{P}(\mathbb{R}) : \Sigma(\mu) > -\infty\}$ and $\nu_. \in \{S_{\mu_0} < \infty\}$. If $k$ is the field such that $(\nu, k)$ satisfies (5.2.5), we set $u_t := \partial_x k_t(x) + H\nu_t(x)$. Then, $\nu_t(dx) \ll dx$ for almost all $t$ and*

$$\begin{aligned}
S_{\mu_0}(\nu) &= \frac{1}{2} \int_0^1 \int [(u_t(x))^2 + (H\nu_t(x))^2] d\nu_t(x) dt - \frac{1}{2}(\Sigma(\nu_1) - \Sigma(\mu_0)) \quad (5.2.6)\\
&= \frac{1}{2} \int_0^1 \int (u_t(x))^2 d\nu_t(x) dt + \frac{\pi^2}{6} \int_0^1 \int \left(\frac{d\nu_t(x)}{dx}\right)^3 dx dt \\
&\quad - \frac{1}{2}(\Sigma(\nu_1) + \Sigma(\mu_0)).
\end{aligned}$$

**Proof:** We shall here only prove formula (5.2.6), assuming the first part of the property which is actually proved by showing that formula (5.2.6) yields for smooth approximations $\nu^\epsilon$ of the measure-valued path $\nu$ such that $S_{\mu_D}(\nu^\epsilon)$ approximate $S_{\mu_D}(\nu)$, yielding a uniform control on the $L_3$ norm of their densities in terms of $S_{\mu_D}(\nu)$. This uniform controls allow us to show that any path $\nu \in \{S_{\mu_D} < \infty\}$ is absolutely continuous with respect to Lebesgue measure and with density in $L_3(dxdt)$.

Let us denote by $k$ the field associated with $\nu$, i.e such that for any $f \in \mathcal{C}^{2,1}_b(\mathbb{R} \times [0,1])$,

$$\begin{aligned}
\int f(x,t) d\nu_t(x) - \int f(x,s) d\nu_s(x) &= \int_s^t \int \partial_v f(x,v) d\nu_v(x) ds \\
&\quad + \frac{1}{2} \int_s^t \int \int \frac{\partial_x f(x,v) - \partial_x f(y,v)}{x-y} d\nu_v(x) d\nu_v(y) dv \\
&\quad + \int_s^t \int \partial_x f(x,v) \partial_x k(x,v) d\nu_v(x) dv \quad (5.2.7)
\end{aligned}$$

with $\partial_x k \in L^2(d\nu_t(x) \times dt)$. Observe that by [119], p. 170, for any $s \in [0,1]$ such that $\nu_s$ is absolutely continuous with respect to Lebesgue measure with density $\rho_s \in L^3(dx)$, for any compactly supported measurable function $\partial_x f(.,s)$,

$$\int \int \frac{\partial_x f(x,s) - \partial_x f(y,s)}{x-y} d\nu_s(x) d\nu_s(y) = 2 \int \partial_x f(x,s) H\nu_s(x) dx ds.$$

*A. Guionnet/Large deviations for random matrices* 133Hence, since we assumed that $\nu_s(dx) = \rho_s(x)dx$ for almost all $s$ with $\rho \in L^3(dxdt)$, (5.2.7) shows that for any $f \in \mathcal{C}_b^{2,1}(\mathbb{R} \times [0,1])$

$$\int f(x,t)\rho_t(x)dx - \int f(x,s)\rho_s(x)ds = \int_s^t \int \partial_v f(x,v)\rho_v(x)dxds \quad (5.2.8)$$
$$+ \int_s^t \int \partial_x f(x,v)u(x,v)\rho_v(x)dxdv,$$

i.e. that in the sense of distributions on $\mathbb{R} \times [0,1]$,

$$\partial_s \rho_s + \partial_x(u_s \rho_s) = 0, \quad (5.2.9)$$

Moreover, since $H\nu_\cdot$ belongs to $L^2(d\nu_s \times ds)$ when $\rho \in L^3(dxdt)$, we can write

$$2S_{\mu_0}(\nu_\cdot) = <k,k>_{0,1}^\nu = \int_0^1 \int_\mathbb{R} (u_s(x))^2 d\nu_s(x)ds + \int_0^1 \int_\mathbb{R} (H\nu_s(x))^2 d\nu_s(x)ds$$
$$-2\int_0^1 \int_\mathbb{R} H\nu_s(x)u_s(x)d\nu_s(x)ds \quad (5.2.10)$$

We shall now see that the last term in the above right hand side only depends on $(\mu_0, \mu_1)$. To simplify, let us assume that $\nu$ is a smooth path such that $f_t(x) = \int \log|x-y|d\nu_t(y)$ is in $\mathcal{C}_b^{2,1}(\mathbb{R} \times [0,1])$. Then, (5.2.7) yields

$$\Sigma(\nu_1) - \Sigma(\nu_0) = 2 \int_0^1 \int_\mathbb{R} H(\nu_s)(x)u_s(x)d\nu_s(x)ds \quad (5.2.11)$$

which gives the result. The general case is obtained by smoothing by free convolution, as can be seen in [61], p. 537-538. □

To prove Theorem 5.4, one should therefore write

$$F_{Ising} = -\inf\{\mathcal{L}(\mu,\nu,\rho^*,m^*);$$
$$\mu_t^*(dx) = \rho_t^*(x)dx \in \mathcal{C}([0,1],\mathcal{P}(\mathbb{R})), \mu_0^* = \mu, \mu_1^* = \nu, \partial_t \rho_t^* + \partial_x m_t^* = 0\}$$

with

$$\mathcal{L}(\mu,\nu,\rho^*,m^*) := \mu(P_1 - \frac{1}{2}x^2) + \nu(P_2 - \frac{1}{2}x^2) - \frac{\beta}{4}(\Sigma(\mu) + \Sigma(\nu))$$
$$+\frac{\beta}{4}\left(\int_0^1 \int \frac{(m_t^*(x))^2}{\rho_t^*(x)}dxdt + \frac{\pi^2}{3}\int_0^1 \int \rho_t^*(x)^3 dxdt\right)$$

It is easy to see that

$$(\mu,\nu,\rho^*,m^*) \in \mathcal{P}(\mathbb{R})^2 \times \mathcal{C}([0,1],\mathcal{P}(\mathbb{R})) \times L^2((\rho_t^*(x))^{-1}dxdt) \to \mathcal{L}(\mu,\nu,\rho^*,m^*)$$

is strictly convex. Therefore, since $F_{Ising}$ is its infimum under a **linear** constraint, this infimum is taken at a unique point $(\mu_A, \mu_B, \rho^{A\to B}, u^{A\to B})$. To describe such a minimizer, one classically performs variations.

The only point to take care of here is that perturbations should be constructed so that the constraint remains satisfied. This type of problems has been dealt with before. In the case above, I met three possible strategies.

The first is to use a target type perturbation, which is a standard perturbation on the space of probability measure, viewed as a subspace of the vector space of measures. The second is to make a perturbation with respect to the source. This strategy was followed by D. Serre in [108]. The idea is basically to set $a_t(x) = a(t,x) = (\rho_t^*(x), \rho_t^*(x)u_t^*(x))$ so that the constraint reads $\text{div}(a_t(x)) = 0$ and perturb $a$ by considering a family

$$a^g = J_g(a.\nabla_{x,t}h) \circ g = J_g(\rho^*(\partial_t h + u^*\partial_x h)) \circ g$$

with a $\mathcal{C}^\infty$ diffeomorphism $g$ of $[0,1] \times \mathbb{R}$ with inverse $h = g^{-1}$ and Jacobian $J_g$. Note here that $div(a_t^g(x)) = 0$. Such an approach yields the Euler's equation (5.2.4).

The last way is to use convex analysis, following for instance Y. Brenier (see [26], Section 2). These two last strategies can only be applied when we know a priori that $(\mu_A, \mu_B)$ are compactly supported (this indeed guarantees some a priori bounds on $(\rho^*, u^*\rho^*)$ for instance). When the minimizers are smooth, all these strategies should give the same result. It is therefore important to study these regularity properties.

It is quite hard to obtain good smoothness properties of the minimizers directly from the fact that they minimize $\mathcal{L}$. An alternative possibility is to come back to the matricial integrals and find directly there some of these properties (cf. [62]). What I did in [61] was to obtain directly the following informations

**Property 5.6.** *1) If $P(x) \geq a|x|^4 + b$, $Q(x) \geq a|x|^4 + b$ for some $a > 0$ there exists a finite constant $C$ such that for any $N \in \mathbb{N}$, there exists $k(N)$, $k(N)$ going to infinity with $N$, such that*

$$\frac{1}{N} Tr(\mathbf{A}_N^{2p}) \leq C^p, \quad \frac{1}{N} Tr(\mathbf{B}_N^{2p}) \leq C^p, \qquad p \leq k(N)$$

*for $\mu_{Ising}^N$-almost all $(\mathbf{A}_N, \mathbf{B}_N)$*

*2) There exists a sequence $(\widetilde{\mathbf{A}}_N, \widetilde{\mathbf{B}}_N)$ of matrices such that $\hat{\mu}_{\widetilde{\mathbf{A}}_N}^N$ (resp.$\hat{\mu}_{\widetilde{\mathbf{B}}_N}^N$) converges toward $\mu_A$ (resp. $\mu_B$) and a Wigner matrix $\mathbf{X}_N$, independent from $(\widetilde{\mathbf{A}}_N, \widetilde{\mathbf{B}}_N)$ such that $\hat{\mu}_{t\widetilde{\mathbf{B}}_N + (1-t)\widetilde{\mathbf{A}}_N + \sqrt{t(1-t)}X_N}^N$ converges toward $\mu_t^*(dx) = \rho_t^*(x)dx$.*

The first result tells us that the polynomial potentials $P, Q$ force the limiting laws $(\mu_A, \mu_B)$ to be compactly supported. The second point is quite natural also when we think that we are looking at matrices with Gaussian entries with given values at time zero and one; the brownian bridge is the path which takes the lowest energy to do that and therefore it is natural that the minimizer of $S_{\mu_D}$ should be the limiting distribution of matrix-valued Brownian bridge. As we shall see in Chapter 6, such a distribution has a limit in free probability as soon as $\hat{\mu}^N_{t\mathbf{B}_N+(1-t)\mathbf{A}_N}$ converges for all $t \in [0,1]$, it is nothing but the distribution of a free Brownian bridge $tB + (1-t)A + \sqrt{t(1-t)}S$ between $A$ and $B$, with a semi-circular variable $S$, free with $(A, B)$. Note however that unless the joint law of $(A, B)$ is given, this result does not describe entirely the law $\mu_t^*$. It however shows, because the limiting law is a free convolution by a semicircular law, that we have the following

**Corollary 5.7.** *1)$\mu_A$ and $\mu_B$ are compactly supported.*

*2) a) There exists a compact set $K \subset \mathbb{R}$ so that for all $t \in [0,1]$, $\mu_t^*(K^c) = 0$. For all $t \in (0,1)$, the support of $\mu_t^*$ is the closure of its interior (in particular it does not put mass on points)*

*b) $\mu_t^*(dx) \ll dx$ for all $t \in (0,1)$. Denote $\rho_t^*(x) = \frac{d\mu_t^*(x)}{dx}$.*

*c) There exists a finite constant $C$ (independent of $t$) so that, $\mu_t^*$ almost surely,*

$$\rho_t^*(x)^2 + (H\mu_t^*(x))^2 \leq (t(1-t))^{-1}$$

*and*

$$|u_t^*(x)| \leq C(t(1-t))^{-\frac{1}{2}}.$$

*d) $(\rho^*, u^*)$ are analytic in the interior of $\Omega = \{x, t \in \mathbb{R} \times [0,1] : \rho_t^*(x) > 0\}$.*

*e) At the boundary of $\Omega_t = \{x \in \mathbb{R} : \rho_t^*(x) > 0\}$, for $x \in \Omega_t$,*

$$|\rho_t^*(x)^2 \partial_x \rho_t^*(x)| \leq \frac{1}{4\pi^3 t^2(1-t)^2} \quad \Rightarrow \quad \rho_t^*(x) \leq \left(\frac{3}{4\pi^3 t^2(1-t)^2}\right)^{\frac{1}{3}} (x-x_0)^{\frac{1}{3}}$$

*if $x_0$ is the nearest point of $x$ in $\Omega_t^c$.*

All these properties are due to the free convolution by the semi-circular law, and are direct consequences of [15]. Once Corollary 5.7 is given, the variational study of $F_{Ising}$ is fairly standard and gives Theorem 5.4. We do not detail this proof here. Hence, free probability arises naturally when we deal with the study of the rate function $S_{\mu_D}$ and more generally when we consider traces of matrices with size going to infinity. We describe the basis of this rapidly developing field in the next chapter.

### 5.3. Discussion and open problems

Matrix models have been much more studied in physics than in mathematics and raise numerous open questions.

As in the last chapter, it would be extremely interesting to understand the relation between the asymptotics of the free energy and the asymptotics of its derivatives at the origin, which are the objects of primary interest. Once this question would be settled, one should understand how to retrieve informations on this serie of numbers from the limiting free energy. In fact, it would be tempting for instance to describe the radius of convergence of this serie via the phase transition of the model, since both should coincide with a default of analyticity of the free energy/the generating function of this serie. However, for instance for the Ising model, one should realize that the criticality is reached at negative temperature where the integral actually diverges. Eventhough the free energy can still be defined in this domain, its description as a generating function becomes even more unclear. There is however some challenging works on this subject in physics (cf. [21] for instance).

There are many other matrix models to be understood with nice combinatorial interpretations. A few were solved in physics literature by means of character expansions for instance in the work of V. Kazakov, M. Staudacher, I. Kostov or P. Zinn Justin. However, such expansions are signed in general and the saddle points methods used rarely justified. In one case, M. Maida and myself [62] could give a mathematical understanding to this method. In general, even the question of the existence of the free energies for most matrix models (that is the convergence of $N^{-2}\log Z_N(P)$) is open and would actually be of great interest in free probability (see the discussion at the end of Chapter 7).

Related with string theory is the question of the understanding of the full expansion of the partition function in terms of the dimension $N$ of the matrices, and the definition of critical exponents. This is still far from being understood on a rigorous ground, eventhought the present section showed that at least for $AB$ interaction models, a good strategy could be to understand better non-intersecting Brownian paths. However, it is yet not clear how to concatenate such an approach with the technology, currently used for one matrix model given in [46], which is based on orthogonal polynomials and Riemann-Hilbert techniques. It would be very tempting to try to generalize precise Laplace's methods which are commonly used to understand the second order corrections of the free energy of mean field interacting particle systems [20]. However, such an approach until now failed even in the one matrix case due to the singularity of the logarithmic interacting potential (cf. [33]). Another approach to this problem has recently been proposed by B. Eynard et all[50, 49] and M. Bertola [14].

On a more analytic point of view, it would be interesting to understand better the properties of complex Burgers equations; we have here deduced most of the smoothness properties of the solution by recalling its realization in free probability terms. A direct analysis should be doable. Moreover, it would be nice to understand how holes in the initial density propagates along time; this might as well be related with the phase transition phenomena according to A. Matytsin and P. Zaugg [92].

## Chapter 6

# Large random matrices and free probability

Free probability is a probability theory for non-commutative variables. In this field, random variables are operators which we shall assume hereafter selfadjoint. For the sake of completeness, but actually not needed for our purpose, we shall recall some notions of operator algebra. We shall then describe free probability as a probability theory on non-commutative functionals, a point of view which forgets the space of realizations of the laws. Then, we will see that free probability is the right framework to consider large random matrices. Finally, we will sketch the proofs of some results we needed in the previous chapter.

**6.1. A few notions about von Neumann algebras**

**Definition 6.1 (Definition 37, [16]).** *A $C^*$-algebra $(\mathcal{A}, *)$ is an algebra equipped with an involution $*$ and a norm $||.||_\mathcal{A}$ which furnishes it with a Banach space structure and such that for any $\mathbf{X}, \mathbf{Y} \in \mathcal{A}$,*

$$\|\mathbf{XY}\|_\mathcal{A} \leq \|\mathbf{X}\|_\mathcal{A} \|\mathbf{Y}\|_\mathcal{A}, \quad \|\mathbf{X}^*\|_\mathcal{A} = \|\mathbf{X}\|_\mathcal{A}, \quad \|\mathbf{XX}^*\|_\mathcal{A} = \|\mathbf{X}\|_\mathcal{A}^2.$$

$\mathbf{X} \in \mathcal{A}$ is self-adjoint iff $\mathbf{X}^* = \mathbf{X}$. $\mathcal{A}_{sa}$ denote the set of self-adjoint elements of $\mathcal{A}$. A $C^*$-algebra $(\mathcal{A}, *)$ is said unital if it contains a neutral element denoted $I$.

$\mathcal{A}$ can always be realized as a sub-$C^*$-algebra of the space $B(H)$ of bounded linear operators on a Hilbert space $H$. For instance, if $\mathcal{A}$ is a unital $C^*$-algebra furnished with a positive linear form $\tau$, one can always construct such a Hilbert space $H$ by completing and separating $L^2(\tau)$ (this is the Gelfand-Naimark-Segal (GNS) construction, see [115], Theorem 2.2.1).

We shall restrict ourselves to this case in the sequel and denote by $H$ a Hilbert space equipped with a scalar product $<.,.>_H$ such that $\mathcal{A} \subset B(H)$.

**Definition 6.2.** *If $\mathcal{A}$ is a sub-$C^*$-algebra of $B(H)$, $\mathcal{A}$ is a von Neumann algebra iff it is closed for the weak topology, generated by the semi-norms family $\{p_{\xi,\eta}(\mathbf{X}) = <\mathbf{X}\xi, \eta>_H, \xi, \eta \in H\}$.*

Let us notice that by definition, a von Neumann algebra contains only bounded operators. The theory nevertheless allows us to consider unbounded operators thanks to the notion of affiliated operators. An operator **X** on $H$ is said to be affiliated to $\mathcal{A}$ iff for any Borel function $f$ on the spectrum of **X**, $f(\mathbf{X}) \in \mathcal{A}$ (see [104], p. 164). Here, $f(\mathbf{X})$ is well defined for any operator **X** as the operator with the same eigenvectors than **X** and eigenvalues given by the image of those of **X** by the map $f$. Note also that if **X** and **Y** are affiliated with $\mathcal{A}$, $a\mathbf{X} + b\mathbf{Y}$ is also affiliated with $\mathcal{A}$ for any $a, b \in \mathbb{R}$.

A state $\tau$ on a unital von Neumann algebra $(\mathcal{A}, *)$ is a linear form on $\mathcal{A}$ such that $\tau(\mathcal{A}_{sa}) \subset \mathbb{R}$ and

1. **Positivity** $\tau(\mathbf{AA}^*) \geq 0$, for any $\mathbf{A} \in \mathcal{A}$.
2. **Total mass** $\tau(I) = 1$

A tracial state satisfies the additional hypothesis

3. **Traciality** $\tau(\mathbf{AB}) = \tau(\mathbf{BA})$ for any $\mathbf{A}, \mathbf{B} \in \mathcal{A}$.

The couple $(\mathcal{A}, \tau)$ of a von Neumann algebra equipped with a state $\tau$ is called a $W^*$- probability space.

**Example 6.3.**   1. Let $n \in \mathbb{N}$, and consider $\mathcal{A} = M_n(\mathbb{C})$ as the set of bounded linear operators on $\mathbb{C}^n$. For any $v \in \mathbb{C}^n$, $\|v\|_{\mathbb{C}^n} = 1$,

$$\tau_v(M) = <v, Mv>_{\mathbb{C}^n}$$

is a state. There is a unique tracial state on $M_n(\mathbb{C})$ which is the normalized trace

$$tr(M) = \frac{1}{n} \sum_{i=1}^{n} M_{ii}.$$

2. Let $(X, \Sigma, d\mu)$ be a classical probability space. Then $\mathcal{A} = L^{\infty}(X, \Sigma, d\mu)$ equipped with the expectation $\tau(f) = \int f d\mu$ is a (non-)commutative probability space. Here, $L^{\infty}(X, \Sigma, d\mu)$ is identified with the set of bounded linear operators on the Hilbert space $H$ obtained by separating $L^2(X, \Sigma, d\mu)$ (by the equivalence relation $f \simeq g$ iff $\mu((f-g)^2) = 0$). The identification follows from the multiplication operator $M(f)g = fg$. Observe that it is weakly closed for the semi-norms $(<f, .g>_H, f, g \in L^2(\mu))$ as $L^{\infty}(X, \Sigma, d\mu)$ is the dual of $L^1(X, \Sigma, d\mu)$.

3. Let $G$ be a discrete group, and $(e_h)_{h \in G}$ be a basis of $\ell^2(G)$. Let $\lambda(h)e_g = e_{hg}$. Then, we take $\mathcal{A}$ to be the von Neumann algebra generated by the linear span of $\lambda(G)$. The (tracial) state is the linear form such that $\tau(\lambda(g)) = 1_{g=e}$ ($e$ = neutral element).

We refer to [128] for further examples and details.

The notion of law $\tau_{\mathbf{X}_1, \ldots, \mathbf{X}_m}$ of $m$ operators $(\mathbf{X}_1, \ldots, \mathbf{X}_m)$ in a $W^*$-probability space $(\mathcal{A}, \tau)$ is simply given by the restriction of the trace $\tau$ to the algebra generated by $(\mathbf{X}_1, \ldots, \mathbf{X}_m)$, that is by the values

$$\tau_{\mathbf{X}_1,\ldots,\mathbf{X}_m}(P) = \tau(P(\mathbf{X}_1,\ldots,\mathbf{X}_m)), \quad \forall P \in \mathbb{C}\langle X_1, \ldots X_m\rangle$$

where $\mathbb{C}\langle X_1, \ldots X_m\rangle$ is the set of polynomial functions of $m$ non-commutative variables.

### 6.2. Space of laws of $m$ non-commutative self-adjoint variables

Following the above description, laws of $m$ non-commutative self-adjoint variables can be seen as elements of the set $\mathcal{M}^{(m)}$ of linear forms on the set of polynomial functions of $m$ non-commutative variables $\mathbb{C}\langle X_1, \ldots X_m\rangle$ furnished with the involution

$$(X_{i_1} X_{i_2} \cdots X_{i_n})^* = X_{i_n} X_{i_{n-1}} \cdots X_{i_1}$$

and such that

1. **Positivity** $\tau(PP^*) \geq 0$, for any $P \in \mathbb{C}\langle X_1, \ldots X_m\rangle$.
2. **Traciality** $\tau(PQ) = \tau(QP)$ for any $P, Q \in \mathbb{C}\langle X_1, \ldots X_m\rangle$.
3. **Total mass** $\tau(I) = 1$

This point of view is identical to the previous one. Indeed, by the Gelfand-Naimark-Segal construction, being given $\mu \in \mathcal{M}^{(m)}$, we can construct a $W^*$-probability space $(\mathcal{A}, \tau)$ and operators $(\mathbf{X}_1, \cdots, \mathbf{X}_m)$ such that

$$\mu = \tau_{\mathbf{X}_1,\ldots,\mathbf{X}_m}. \tag{6.2.1}$$

This construction can be summarized as follows. Consider the bilinear form on $\mathbb{C}\langle X_1, \ldots X_m\rangle^2$ given by

$$<P, Q>_\tau = \tau(PQ^*).$$

We then let $H$ be the Hilbert space obtained as follows. We set $L^2(\tau) = \overline{\mathbb{C}\langle X_1, \ldots X_m\rangle}^{||.||_\tau}$ to be the set of $\mathbb{C}\langle X_1, \ldots X_m\rangle$ for the norm $||.||_\tau = <.,.>_\tau^{\frac{1}{2}}$. We then separate $L^2(\tau)$ by taking the quotient by the left ideal

$$L_\mu = \{F \in L^2(\tau) : ||F||_\tau = 0\}.$$

Then $H = L^2(\tau)/L_\mu$ is a Hilbert space with scalar product $<.,.>_\tau$. The non-commutative polynomials $\mathbb{C}\langle X_1, \ldots X_m\rangle$ act by left multiplication on $L^2(\tau)$ and we can consider the completion of these multiplication operators for the semi-norms $\{<P,.Q>_H, P, Q \in L^2(\tau)\}$, which form a von Neumann algebra $\mathcal{A}$ equipped with a tracial state $\tau$ satisfying (6.2.1). In this sense, we can think about $\mathcal{A}$ as the set of bounded measurable functions $L^\infty(\tau)$.

The topology under consideration is usually in free probability the $\mathbb{C}\langle X_1, \ldots X_m\rangle^*$-topology that is $(\tau_{\mathbf{X}_1^n,\ldots,\mathbf{X}_m^n}, m \in \mathbb{N})$ converges toward $\tau_{\mathbf{X}_1,\ldots,\mathbf{X}_m}$ iff for every $P \in \mathbb{C}\langle X_1, \ldots X_m\rangle$,

$$\lim_{n\to\infty} \tau_{\mathbf{X}_1^n,\ldots,\mathbf{X}_m^n}(P) = \tau_{\mathbf{X}_1,\ldots,\mathbf{X}_m}(P).$$

If $(\mathbf{X}_1^n, \ldots, \mathbf{X}_m^n)_{n \in \mathbb{N}}$ are non-commutative variables whose law $\tau_{\mathbf{X}_1^n, \ldots, \mathbf{X}_m^n}$ converges toward $\tau_{\mathbf{X}_1, \ldots, \mathbf{X}_m}$, then we shall also say that $(\mathbf{X}_1^n, \ldots, \mathbf{X}_m^n)_{n \in \mathbb{N}}$ converges in law or in distribution toward $(\mathbf{X}_1, \ldots, \mathbf{X}_m)$.

Such a topology is reasonable when one deals with uniformly bounded non-commutative variables. In fact, if we consider for $R \in \mathbb{R}^+$,

$$\mathcal{M}_R^{(m)} := \{\mu \in \mathcal{M}^{(m)} : \mu(X_i^{2p}) \leq R^p, \ \forall p \in \mathbb{N}, \ 1 \leq i \leq m\}$$

then it is not hard to see that $\mathcal{M}_R^{(m)}$, equipped with this weak-* topology, is a Polish space (i.e. a complete metric space). A distance can for instance be given by

$$d(\mu, \nu) = \sum_{n \geq 0} \frac{1}{2^n} |\mu(P_n) - \nu(P_n)|$$

where $P_n$ is a dense sequence of polynomials with operator norm bounded by one when evaluated at any set of self-adjoint operators with operator norms bounded by $R$.

This notion is the generalization of laws of $m$ real-valued variables bounded say by a given finite constant $R$, in which case the weak-* topology driven by polynomial functions is the same as the standard weak topology. Actually, it is not hard to check that $\mathcal{M}_R^{(1)} = \mathcal{P}([-R, R])$. However, it can be usefull to consider more general topologies compatible with the existence of unbounded operators, as might be encountered for instance when considering the deviations of large random matrices. Then, the only point is to change the set of test functions. In [29], we considered for instance the complex vector space $\mathcal{CC}_{st}^m(\mathbb{C})$ generated by the Stieljes functionals

$$ST^m(\mathbb{C}) = \{ \prod_{1 \leq i \leq n}^{\rightarrow} (z_i - \sum_{k=1}^m \alpha_i^k \mathbf{X}_k)^{-1}; \quad z_i \in \mathbb{C}\backslash\mathbb{R}, \alpha_i^k \in \mathbb{Q}, n \in \mathbb{N}\} \quad (6.2.2)$$

It can be checked easily that, with such type of test functions, $\mathcal{M}^{(m)}$ is again a Polish space.

**Example 6.4.** *Let $N \in \mathbb{N}$ and consider $m$ Hermitian matrices $A_1^N, \cdots, A_m^N \in \mathcal{H}_N^m$ with spectral radius $\|A_i^N\|_\infty \leq R$, $1 \leq i \leq m$. Then, set*

$$\hat{\mu}_{A_1^N, \cdots, A_m^N}^N(P) = tr\left(P(A_1^N, \cdots, A_m^N)\right), \quad \forall P \in \mathbb{C}\langle X_1, \cdots X_m\rangle.$$

*Clearly, $\hat{\mu}_{A_1^N, \cdots, A_m^N}^N \in \mathcal{M}_R^{(m)}$. Moreover, if $(A_1^N, \cdots, A_m^N)_{N \in \mathbb{N}}$ is a sequence such that*

$$\lim_{N \to \infty} \hat{\mu}_{A_1^N, \cdots, A_m^N}^N(P) = \tau(P), \quad \forall P \in \mathbb{C}\langle X_1, \cdots X_m\rangle$$

*then $\tau \in \mathcal{M}_R^{(m)}$ since $\mathcal{M}_R^{(m)}$ is complete.*

*It is actually a long standing question posed by A. Connes to know whether all $\tau \in \mathcal{M}^{(m)}$ can be approximated in such a way.*

*In the case $m = 1$, the question amounts to ask if for all $\mu \in \mathcal{P}([-R, R])$, there exists a sequence $(\lambda_1^N, \cdots, \lambda_N^N)_{N \in \mathbb{N}}$ such that*

$$\lim_{N \to \infty} \frac{1}{N} \sum_{i=1}^{N} \delta_{\lambda_i^N} = \mu.$$

*This is well known to be true by Birkhoff's theorem (which is based on Krein-Milman's theorem), but still an open question when $m \geq 2$.*

### 6.3. Freeness

Free probability is not only a theory of probability for non-commutative variables; it contains also the central notion of freeness, which is the analogue of independence in standard probability.

**Definition 6.5.** *The variables $(\mathbf{X}_1, \ldots, \mathbf{X}_m)$ and $(\mathbf{Y}_1, \ldots, \mathbf{Y}_n)$ are said to be free iff for any $(P_i, Q_i)_{1 \leq i \leq n} \in (\mathbb{C}\langle X_1, \cdots, X_m \rangle \times \mathbb{C}\langle X_1, \ldots, X_n \rangle)^n$,*

$$\tau\left(\overrightarrow{\prod_{1 \leq i \leq n}} P_i(\mathbf{X}_1, \ldots, \mathbf{X}_m) Q_i(\mathbf{Y}_1, \ldots, \mathbf{Y}_n)\right) = 0$$

*as soon as*

$$\tau(P_i(\mathbf{X}_1, \ldots, \mathbf{X}_m)) = 0, \qquad \tau(Q_i(\mathbf{Y}_1, \ldots, \mathbf{Y}_n)) = 0, \qquad \forall i \in \{1, \ldots, n\}.$$

<u>*Remark*</u> 6.6:
1) The notion of freeness defines uniquely the law of $\{\mathbf{X}_1, \ldots, \mathbf{X}_m, \mathbf{Y}_1, \ldots, \mathbf{Y}_n\}$ once the laws of $(\mathbf{X}_1, \ldots, \mathbf{X}_m)$ and $(\mathbf{Y}_1, \ldots, \mathbf{Y}_n)$ are given (in fact, check that every expectation of any polynomial is given uniquely by induction over the degree of this polynomial).

2) If $\mathbf{X}$ and $\mathbf{Y}$ are free variables with joint law $\tau$, and $P, Q \in \mathbb{C}\langle X \rangle$ such that $\tau(P(\mathbf{X})) = 0$ and $\tau(Q(\mathbf{Y})) = 0$, it is clear that $\tau(P(\mathbf{X})Q(\mathbf{Y})) = 0$ as it should for independent variables, but also $\tau(P(\mathbf{X})Q(\mathbf{Y})P(\mathbf{X})Q(\mathbf{Y})) = 0$ which is very different from what happens with usual independent commutative variables where $\mu(P(\mathbf{X})Q(\mathbf{Y})P(\mathbf{X})Q(\mathbf{Y})) = \mu(P(\mathbf{X})^2 Q(\mathbf{Y})^2) > 0$.

3) The above notion of freeness is related with the usual notion of freeness in groups as follows. Let $(x_1, ..x_m, y_1, \cdots, y_n)$ be elements of a group. Then, $(x_1, \cdots, x_m)$ is said to be free from $(y_1, \cdots, y_n)$ if any non trivial words in these elements is not the neutral element of the group, that is that for every monomials $P_1, \cdots, P_k \in \mathbb{C}\langle X_1, \cdots, X_m \rangle$ and $Q_1, \cdots, Q_k \in \mathbb{C}\langle X_1, \cdots, X_n \rangle$, $P_1(x)Q_1(y)P_2(x)\cdots Q_k(y)$ is not the neutral element as soon as the $Q_k(y)$ and the $P_i(x)$ are not the neutral element. If we consider, following example 6.3.3), the map which is one on trivial words and zero otherwise and extend it by linearity to polynomials, we see that this define a tracial state on the operators of left multiplication by the elements of the group and that the two notions of freeness coincide.

4) We shall see below that examples of free variables naturally show up when considering random matrices with size going to infinity.

### 6.4. Large random matrices and free probability

We have already seen in example 6.3 that if we consider $(M_1^N, \ldots, M_m^N) \in \mathcal{H}_N^m$, their empirical distribution $\hat{\mu}_{M_1^N, \ldots, M_m^N}^N$ given for any $P \in \mathbb{C}\langle X_1, \ldots, X_m \rangle$ by

$$\hat{\mu}_{M_1^N, \cdots, M_m^N}^N(P) := \operatorname{tr}\left(P(M_1^N, \ldots, M_m^N)\right).$$

belongs to $\mathcal{M}^{(m)}$.

Moreover, if we take a sequence $\{(M_1^N, \ldots, M_m^N) \in \mathcal{H}_N^m\}_{N \in \mathbb{N}}$ such that for any $P \in \mathbb{C}\langle X_1, \ldots, X_m \rangle$ the limit

$$\tau(P) := \lim_{N \to \infty} \operatorname{tr}\left(P(M_1^N, \ldots, M_m^N)\right)$$

exists, then $\tau \in \mathcal{M}^{(m)}$. Hence, $\mathcal{M}^{(m)}$ is the natural space in which one should consider large matrices, as far as their trace is concerned. As we already stressed in example 6.4, it is still unknown whether any element of $\mathcal{M}^{(m)}$ can be approximated by a sequence of empirical distribution of self-adjoint matrices in the case $m \geq 2$. Reciprocally, large random matrices became an important source of examples of operators algebras.

The fact that free probability is particularly well suited to study large **random** matrices is due to an observation of Voiculescu [121] who proved that if $(\mathbf{A}_N, \mathbf{B}_N)_{N \in \mathbb{N}}$ is a sequence of uniformly bounded diagonal matrices with converging spectral distribution, and $\mathbf{U}_N$ a unitary matrix following Haar measure $m_2^N$, then the empirical distribution of $(\mathbf{A}_N, \mathbf{U}_N \mathbf{B}_N \mathbf{U}_N^*)$ converges toward the law of $(A, B)$, $A$ and $B$ being free and each of their law being given by their limiting spectral distribution. This convergence holds in expectation with respect to the unitary matrices $U_N$ [121] and then almost surely (as can be checked by Borel-Cantelli's lemma and by controlling
$\int (\hat{\mu}_{\mathbf{A}_N, \mathbf{U}_N \mathbf{B}_N \mathbf{U}_N^*}^N(P) - \int \hat{\mu}_{\mathbf{A}_N, \widetilde{\mathbf{U}}_N \mathbf{B}_N \widetilde{\mathbf{U}}_N^*}^N(P) dm_2^N(\widetilde{\mathbf{U}}_N))^2 dm_2^N(\mathbf{U}_N)$ as $N$ goes to infinity).

As a consequence, if one considers the joint distribution of $m$ independent Wigner matrices $(\mathbf{X}_1^N, \cdots, \mathbf{X}_m^N)$, their empirical distribution converges almost surely toward $(\mathbf{S}_1, \cdots, \mathbf{S}_m)$, $m$ free variables distributed according to the semi-circular law $\sigma(dx) = C\sqrt{4 - x^2} dx$.

Hence, freeness appears very naturally in the context of large random matrices. The semi-circular law $\sigma$, which we have seen to be the asymptotic spectral distribution of Gaussian Wigner matrices, is in fact also deeply related with the notion of freeness; it plays the role that Gaussian law has with the notion of independence in the sense that it gives the limit law of the analogue of central limit theorem. Indeed, let $(\mathcal{A}, \tau)$ be a $W^*$-probability space and $\{\mathbf{X}_i, i \in \mathbb{N}\} \in \mathcal{A}$ be free random variables which are centered ($\tau(X_i) = 0$) and with covariance 1 ($\tau(X_i^2) = 1$). Then, the sum $\sqrt{n}^{-1}(\mathbf{X}_1 + \cdots + \mathbf{X}_n)$ converges in distribution toward a semi-circular variable. We shall be interested in the next section in the free Itô's calculus which appeared naturally in our problems.

Before that, let us introduce the notation for free convolution : if $\mathbf{X}$ (resp. $\mathbf{Y}$) is a random variable with law $\mu$ (resp $\nu$) and $\mathbf{X}$ and $\mathbf{Y}$ are free, we denote $\mu \boxplus \nu$ the law of $\mathbf{X} + \mathbf{Y}$. There is a general formula to describe the law $\mu \boxplus \nu$; in fact, analytic functions $R_\mu$ were introduced by Voiculescu as an analogue of the logarithm of Fourier transform in the sense that

$$R_{\mu \boxplus \nu}(z) = R_\mu(z) + R_\nu(z)$$

and that $R_\mu$ defines $\mu$ uniquely. I refer the reader to [128] for more details. Convolution by a semicircular variable was precisely studied by Biane [15].

### 6.5. Free processes

The notion of freeness allows us to construct a free Brownian motion such as

**Definition 6.7.** *A free Brownian motion $\{S_t, t \geq 0\}$ is a process such that*
  1) $S_0 = 0$.
  2) *For any $t \geq s \geq 0$, $S_t - S_s$ is free with the algebra $\sigma(S_u, u \leq s)$ generated by $(S_u, u \leq s)$.*
  3) *For any $t \geq s \geq 0$, the law of $S_t - S_s$ is a semi-circular distribution with covariance $s - t$ ; $\sigma_{s-t}(dx) = (\pi(s-t))^{-1}\sqrt{4(s-t) - x^2}dx$.*

The Hermitian Brownian motion in particular converges toward the free Brownian motion according to the previous section (using convergence on cylinder functions).

As for the Brownian motion, one can make sense of the free differential equation

$$dX_t = dS_t + b_t(X_t)dt$$

and show existence and uniqueness of solution to this equation when $b_t$ is Lipschitz operator in the sense that if $||.||$ denotes the operator norm on the algebra $\mathcal{A}$ ($||A|| = \lim \tau(A^{2n})^{\frac{1}{2n}}$) on which the free Brownian motion $\{S_t, t \geq 0\}$ lives,

$$||b_t(X) - b_t(Y)|| \leq C||X - Y||$$

with a finite constant $C$(one just uses a Picard argument).

With the intuition given by stochastic calculus, we shall give some outline of the proof of some results needed in Chapter 4.

### 6.6. Continuity of the rate function under free convolution

In the classical case, the entropy $S$ of the deviations of the law of the empirical measure of independent Brownian motion decreases by convolution (see (4.2.20)). We want here to generalize this result to our eigenvalues setting. The intuition coming from the classical case, adapted to the free probability setting, will help to show the following result :

**Lemma 6.8.** *For any $p \in \mathcal{P}(\mathbb{R})$, any $\nu \in \mathcal{C}([0,1], \mathcal{P}(\mathbb{R}))$*

$$S^{0,1}(\nu \boxplus p) \leq S^{0,1}(\nu).$$

**Proof :** We shall give the philosophy of the proof via the formula (see (5.2.5))

$$S^{0,1}(\nu) = \frac{1}{2} <k,k>_{0,1}^{\nu} = \frac{1}{2}\int_0^1 \int (\partial_x k_t)^2 d\nu_t(x) dt.$$

Namely, let us assume that $\nu_t$ can be represented as the law at time $t$ of the free stochastic differential equation (FSDE)

$$dX_t = dS_t + k_t(X_t)dt$$

It can be checked by free Itô's calculus that this FSDE satisfies the same free Fokker-Planck equation (5.2.5) (get the intuition by replacing $S$ by the Hermitian Brownian motion and $X_{\cdot}$ by a matrix-valued process). However, until uniqueness of the solutions of (5.2.5) is proved, it is not clear that $\nu$ is indeed the law of this FSDE.

Now, let $C$ be a random variable with law $p$, free with $S$ and $X_0$. $Y = X + C$ satisfies the same FSDE

$$dY_t = dS_t + \partial_x k_t(X_t)dt$$

and therefore its law $\mu_t = \nu_t \boxplus p$ satisfies for any $\mathcal{C}_b^{2,1}(\mathbb{R} \times [0,1])$,

$$\begin{aligned} S^{0,1}(\mu, f) &= \int_0^1 \tau(\partial_x f(X_t + C)\partial_x k_t(X_t))dt \\ &= \int_0^1 \tau(\partial_x f(X_t + C)\tau(\partial_x k_t(X_t)|X_t + C))dt \end{aligned}$$

where $\tau(\ |X_t+C)$ is the orthogonal projection in $L^2(\tau)$ on the algebra generated by $X_t + C$ (recall the definition of $L^2(\tau)$ given in Section 6.2). From this, we see

that

$$\begin{aligned}
S^{0,1}(\mu) &= \sup_{f \in \mathcal{C}_b^{2,1}(\mathbb{R}\times[0,1])} \{S^{0,1}(\mu, f) - \frac{1}{2} <f,f>_{0,1}^{\mu}\} \\
&= \sup_{f \in \mathcal{C}_b^{2,1}(\mathbb{R}\times[0,1])} \{\int_0^1 \tau(\partial_x f_t(X_t + C)\tau(\partial_x k_t(X_t)|X_t + C))dt \\
&\quad -\frac{1}{2}\int_0^1 \tau[(\partial_x f_t(X_t+C))^2]dt\} \\
&\leq \frac{1}{2}\int_0^1 \tau(\tau(\partial_x k_t(X_t)|X_t+C)^2)dt \\
&\leq \frac{1}{2}\int_0^1 \tau(\partial_x k_t(X_t)^2)dt = S^{0,1}(\nu)
\end{aligned}$$

Of course, such an inequality has nothing to do with the existence and uniqueness of a strong solution of our free Fokker-Planck equation and we can indeed prove this result by means of $R$-transform theory for any $\nu \in \{S^{0,1} < \infty\}$ (see [30]).

□

### 6.7. The infimum of $S_{\mu_D}$ is achieved at a free Brownian bridge

Let us state more precisely the theorem obtained in this section. A free Brownian bridge between $\mu_0$ and $\mu_1$ is the law of

$$X_t = (1-t)X_0 + tX_1 + \sqrt{t(1-t)}S \qquad (6.7.3)$$

with a semicircular variable $S$, free with $X_0$ and $X_1$, with law $\mu_0$ and $\mu_1$ respectively. We let $FBB(\mu_0, \mu_1) \subset \mathcal{C}([0,1], \mathcal{P}(\mathbb{R}))$ denote the set of such laws (which depend of course not only on $\mu_0, \mu_1$ but on the joint distribution of $(X_0, X_1)$ too). Then, we shall prove (cf. Appendix 4 in [61] that

**Theorem 6.9.** *Assume $\mu_0, \mu_1$ compactly supported, with support included into $[-R, R]$.*
*1) Then,*

$$\begin{aligned}
J_\beta(\mu_0, \mu_1) &= \frac{\beta}{2}\inf\{S(\nu); \nu_0 = \mu_0, \nu_1 = \mu_1\} \\
&= \frac{\beta}{2}\inf\{S(\nu); \nu \in FBB(\mu_0, \mu_1)\}.
\end{aligned}$$

*2) Since $FBB(\mu_0, \mu_1)$ is a closed subset of $\mathcal{C}([0,1], \mathcal{P}(\mathbb{R}))$, the unique minimizer in $J_\beta(\mu_0, \mu_1)$ belongs to $FBB(\mu_0, \mu_1)$.*

The proof of Theorem 6.9 is rather technical and goes back through the large random matrices origin of $J_\beta$. Let us consider the case $\beta = 2$. By definition, if

$$X_t^N = X_0^N + H_t^N$$

with a real diagonal matrix $X_0^N$ with spectral measure $\hat\mu_0^N$ and a Hermitian Brownian motion $H^N$, if we denote $\hat\mu_t^N$ the spectral measure of $X_t^N$, then, if $\hat\mu_0^N$ converges toward a compactly supported probability measure $\mu_0$, for any $\mu_1 \in \mathcal{P}(\mathbb{R})$,

$$\limsup_{\delta \to 0} \limsup_{N \to \infty} \frac{1}{N^2} \log \mathbb{P}(d(\hat\mu_1^N, \mu_1) < \delta) \le -\inf\{S_{\mu_D}(\nu_\cdot); \nu_1 = \mu_1\}.$$

Let us now reconsider the above limit and show that the limit must be taken at a free Brownian bridge. More precisely, we shall see that, if $\tau$ denotes the joint law of $(X_0, X_1)$ and $\mu^\tau$ the law of the free Brownian bridge (6.7.3) associated with the distribution $\tau$ of $(X_0, X_1)$,

$$\limsup_{\delta \to 0} \limsup_{N \to \infty} \frac{1}{N^2} \log \mathbb{P}(d(\hat\mu_1^N, \mu_1) < \delta)$$
$$\le \sup_{\substack{\tau \circ X_0^{-1} = \mu_0 \\ \tau \circ X_1^{-1} = \mu_1}} \limsup_{\delta \to 0} \limsup_{N \to \infty} \frac{1}{N^2} \log \mathbb{P}(\max_{1 \le k \le n} d(\hat\mu_{t_k}^N, \mu_{t_k}^\tau) \le \delta)$$

for any family $\{t_1, \cdots, t_n\}$ of times in $[0, 1]$. Therefore, Theorem 4.5.2).a) implies that

$$\limsup_{\delta \to 0} \limsup_{N \to \infty} \frac{1}{N^2} \log \mathbb{P}(d(\hat\mu_1^N, \mu_1) < \delta)$$
$$\le \inf\{S(\mu^\tau), \tau \circ X_0^{-1} = \mu_0, \tau \circ X_1^{-1} = \mu_1\}.$$

The lower bound estimate obtained in Theorem 4.5.2).b) therefore guarantees that

$$\inf\{S(\nu), \nu_0 = \mu_0, \nu_1 = \mu_1\} \ge \inf\{S(\mu^\tau), \tau \circ X_0^{-1} = \mu_0, \tau \circ X_1^{-1} = \mu_1\}.$$

and therefore the equality since the other bound is trivial.

Let us now be more precise. We consider the empirical distribution $\hat\mu_{0,1}^N = \hat\mu_{X_0^N, X_1^N}^N$ of the couple of the initial and final matrices of our process as an element of $\mathcal{M}_1^{(2)}$, equipped with the topology of the Stieljes functionals. It is not hard to see that $\mathcal{M}_1^{(2)}$ is a compact metric space by the Banach Alaoglu theorem. Let $D$ be a distance on $\mathcal{M}_1^{(2)}$.

Then, for any $\epsilon > 0$, we can find $M \in \mathbb{N}$, $(\tau_k)_{1 \le k \le M}$ so that $\mathcal{M}_1^{(2)} \subset \cup_{1 \le k \le M}\{\tau : D(\tau, \tau_k) < \epsilon\}$ and therefore

$$\limsup_{\delta \to 0} \limsup_{N \to \infty} \frac{1}{N^2} \log \mathbb{P}(d(\hat\mu_1^N, \mu_1) < \delta)$$
$$\le \max_{1 \le k \le M} \limsup_{N \to \infty} \frac{1}{N^2} \log \mathbb{P}(d(\hat\mu_1^N, \mu_1) < \delta; D(\hat\mu_{0,1}^N, \tau_k) < \epsilon)$$

Now, conditionally on $X_1^N$,

$$dX_t^N = dH_t^N - \frac{X_t^N - X_1^N}{1-t}dt$$

or equivalently

$$X_t^N = tX_1^N + (1-t)X_0^N + (1-t)\int_0^t (1-s)^{-1}dH_s^N.$$

It is not hard to see that when $\hat{\mu}_{0,1}^N$ converges toward $\tau$, $\hat{\mu}_{X_t^N}^N$ converges toward $\mu_t^\tau$ for all $t \in [0,1]$.

Therefore, for any $\kappa > 0$, any $t_1, \cdots, t_n \in [0,1]$, there exists $\epsilon > 0$ such that for any $(X_0^N, X_1^N) \in \{D(\hat{\mu}_{0,1}^N, \tau) < \epsilon\}$,

$$\mathbb{P}(\max_{1 \leq k \leq n} d(\hat{\mu}_{X_{t_k}^N}^N, \mu_{t_k}^\tau) > \eta | X_1^N) \leq \kappa$$

Hence for any $\eta$, when $\epsilon$ is small enough and $N$ large enough, taking $\kappa = \frac{1}{2}$,

$$\mathbb{P}(d(\hat{\mu}_1^N, \mu_1) < \delta; D(\hat{\mu}_{0,1}^N, \tau) < \epsilon)$$

$$\leq 2\mathbb{P}(d(\hat{\mu}_1^N, \mu_1) < \delta; D(\hat{\mu}_{0,1}^N, \tau) < \epsilon, \max_{1 \leq k \leq n} d(\hat{\mu}_{X_{t_k}^N}^N, \mu_{t_k}^\tau) < \eta).$$

We arrive at, for $\epsilon$ small enough and any $\tau \in \mathcal{M}_{0,1}$,

$$\limsup_{N \to \infty} \frac{1}{N^2} \log \mathbb{P}(d(\hat{\mu}_1^N, \mu_1) < \delta; D(\hat{\mu}_{0,1}^N, \tau) < \epsilon)$$

$$\leq \limsup_{N \to \infty} \frac{1}{N^2} \log \mathbb{P}(\max_{1 \leq k \leq n} d(\hat{\mu}_{t_k}^N, \mu_{t_k}^\tau) < \delta).$$

Using the large deviation upper bound for the law of $(\hat{\mu}_t^N, t \in [0,1])$ of Theorem 4.5.2), we deduce

$$\limsup_{N \to \infty} \frac{1}{N^2} \log \mathbb{P}(d(\hat{\mu}_1^N, \mu_1) < \delta) \quad \leq \quad -\frac{\beta}{2} \min_{1 \leq p \leq M} \inf_{\max_{1 \leq k \leq n} d(\nu_{t_k}, \mu_{t_k}^{\tau_p}) \leq \delta} S_{\mu_D}(\nu)$$

We can now let $\epsilon$ going to zero, and then $\delta$ going to zero, and then $n$ going to infinity, to conclude (since $S$ is a good rate function). It is not hard to see that $\text{FBB}(\mu_0, \mu_1)$ is closed (see p. 565-566 in [61]).

## Chapter 7

# Voiculescu's non-commutative entropies

One of the motivations to construct free probability theory, and in particular its entropy theory, was to study von Neumann algebras and in particular to prove (or disprove) isomorphisms between them. One of the long standing questions in this matter is to see whether the free group factor $L(F^m)$ with $m$ generators is isomorphic with the free group with $n$ generators when $n \neq m$. Such a problem can be rephrased in free probability terms due to the remark that if $X_1, \cdots, X_m$ (resp. $Y_1, \cdots, Y_m$) are non-commutative variables with law $\tau_X$ and $\tau_Y$ respectively, then $\tau_X = \tau_Y \Rightarrow W^*(X_1, \cdots, X_m) \simeq W^*(Y_1, \cdots, Y_m)$. Therefore,
$$W^*(X_1, \cdots, X_m) \simeq W^*(Y_1, \cdots, Y_m) \Leftrightarrow \tau_X \equiv \tau_Y$$
where $\tau_X \equiv \tau_Y$ when there exists $F \in L^\infty(\tau_X), G \in L^\infty(\tau_Y)$ such that
$$\tau_Y(P) = F_\# \tau_X(P) = \tau_X(P \circ F) \quad \tau_X(P) = G_\# \tau_Y(P) \quad \forall P \in \mathbb{C}\langle X_1, \cdots, X_m \rangle.$$

Here, $L^\infty(\tau_X)$ and $L^\infty(\tau_Y)$ denote the elements of the von Neumann algebras obtained by the Gelfand-Naimark-Segal construction. Therefore, if $\sigma_m$ is the law of $m$ free semi-circular variables $S_1, .., S_m$, $L(F^m) \simeq W^*(S_1, \cdots, S_m)$ and
$$L(F^m) \simeq L(F^n) \Leftrightarrow W^*(S_1, \cdots, S_m) \simeq W^*(S_1, \cdots, S_n) \Leftrightarrow \sigma_m \equiv \sigma_n.$$

The isomorphism problem is therefore equivalent to wonder whether $\sigma_m \equiv \sigma_n$ implies $m = n$ or not.

Note that in the classical setting, two probability measures on $\mathbb{R}^m$ and $\mathbb{R}^n$ are equivalent (in the sense above that they are the push forward of each other) provided they have no atoms and regardless of the choice of $m, n \in \mathbb{N}$. In the non-commutative context, one may expect that such an isomorphism becomes false in view for instance of the work of D. Gaboriau [55] in the context of equivalence relations. However, such a setting carries more structure due to the equivalence relations than free group factors; it concerns the case of algebras with Cartan sub-algebras (which are roughly Abelian sub-algebras) whereas D. Voiculescu [124] proved that every separable $II_1$-factor has no Cartan sub-algebras.

To try to disprove the isomorphism, one would like to construct a function $\delta : \mathcal{M}^{(m)} \to \mathbb{R}$ which is an invariant in the sense that for any $\tau, \widetilde{\tau} \in \mathcal{M}^{(m)}$

$$\tau \equiv \widetilde{\tau} \Rightarrow \delta(\tau) = \delta(\widetilde{\tau}) \qquad (7.0.1)$$

and such that $\delta(\sigma_k) = k$ (where $\sigma_k$, $k \leq m$, can be embedded into $\mathcal{M}^{(m)}$ for $k \leq m$ by taking $m - k$ null operators). D. Voiculescu proposed as a candidate the so-called entropy dimension $\delta$, constructed in the spirit of Minkowski dimension (see (7.1.3) for a definition). It is currently under study whether $\delta$ is an invariant of the von Neumann algebra, that is if it satisfies (7.0.1). We note however that the converse implications is false since N. Brown [27] just produced an example showing that there exists a von-Neumann algebra which is not isomorphic to the free group factor but with same entropy dimension. Eventhough this problem has not yet been solved, some other important questions concerning von Neumann algebras could already be answered by using this approach (see [124] for instance). In fact, free entropy theory is still far from being complete and the understanding of these objects is still too narrow to be applied to its full extent. In this last chapter, we will try to complement this theory by means of large deviations techniques.

We begin by the definitions of the two main entropies introduced by D. Voiculescu, namely the microstates entropy $\chi$ and the microstates-free entropy $\chi^*$. The entropy dimension $\delta$ is defined via the microstates entropy $\chi$ but we shall not study it here.

## 7.1. Definitions

We define here a version of Voiculescu's **microstates entropy** with respect to the Gaussian measure (rather than Lebesgue measure as he initially did : these two points of view are equivalent (see [30]) but the Gaussian point of view is more natural after the previous chapters). I underlined the entropies to make this difference, but otherwise kept the same notations as Voiculescu.

For $\mu \in \mathcal{M}_1^{(m)}$, $n \in \mathbb{N}$, $N \in \mathbb{N}$, $\epsilon > 0$, Voiculescu [122] defines a neighborhood $\Gamma_R(\mu, n, N, \epsilon)$ of the state $\mu$ as the set of matrices $\mathbf{A}_1, .., \mathbf{A}_m$ of $\mathcal{H}_N^m$ such that

$$|\mu(\mathbf{X}_{i_1}..\mathbf{X}_{i_p}) - \text{tr}(\mathbf{A}_{i_1}..\mathbf{A}_{i_p})| < \epsilon$$

for any $1 \leq p \leq n$, $i_1, .., i_p \in \{1, .., m\}^p$ and $||\mathbf{A}_j||_\infty \leq R$. Then, the microstates entropy w.r.t the Gaussian measure is given by

$$\underline{\chi}(\mu) \quad := \quad \sup_{R>0} \inf_{n \in \mathbb{N}} \inf_{\epsilon>0} \limsup_{N \to \infty} N^{-2} \log \mu_N^{\otimes m}[\Gamma_R(\mu, n, N, \epsilon)]$$

This definition of the entropy is done in the spirit of Boltzmann and Shannon. In the classical case where $\mathcal{M}_1^m$ is replaced by $\mathcal{P}(\mathbb{R})$ the entropy with respect to the Gaussian measure $\gamma$ is the relative entropy

$$S(\mu) := \inf_{\epsilon>0} \limsup_{N\to\infty} N^{-1} \log \gamma^{\otimes N}\left(d(\frac{1}{N}\sum_{i=1}^{N}\delta_{x_i}, \mu) \leq \epsilon\right)$$

$$:= \inf_{\epsilon>0, k\in\mathbb{N}} \limsup_{N\to\infty} N^{-1} \log \gamma^{\otimes N}\left(|\frac{1}{N}\sum_{i=1}^{N} f_p(x_i) - \mu(f_p)| \leq \epsilon, p \leq k\right)$$

where the last equality holds if $(f_p)_{p\in\mathbb{N}}$ is a family of uniformly continuous functions, dense in $\mathcal{C}_b^0(\mathbb{R})$. This last way to define the entropy is very close from Voiculescu's definition. In the commutative case, we know by Sanov's theorem that

1. The limsup in the definition of $S$ can be replaced by a liminf, i.e.

$$S(\mu) := \inf_{\epsilon>0} \liminf_{N\to\infty} N^{-1} \log \gamma^{\otimes m}\left(d(\frac{1}{N}\sum_{i=1}^{N}\delta_{x_i}, \mu) \leq \epsilon\right).$$

2. For any $\mu \in \mathcal{P}(\mathbb{R})$ we have the formula

$$S(\mu) = -S^*(\mu)$$

with $S^*$ the relative entropy which is infinite if $\mu$ is not absolutely continuous wrt $\gamma$ and otherwise given by

$$S^*(\mu) = \int \frac{d\mu}{d\gamma}(x) \log \frac{d\mu}{d\gamma}(x) d\gamma(x).$$

These two fundamental results are still lacking in the non-commutative theory except in the case $m=1$ where Voiculescu [123] (see also [10]) has shown that the two limits coincide and that

$$\underline{\chi}(\mu) = \Sigma(\mu) - \frac{1}{2}\int x^2 d\mu(x) + \text{constant}$$

with

$$\Sigma(\mu) = \int\int \log|x-y| d\mu(x) d\mu(y).$$

In the case where $m \geq 2$, Voiculescu [122, 125] proposed an analogue of $S^*$, $\chi^*$, which does not depend at all on the definition of microstates and called for that reason **micro-states-free** entropy. The definition of $\underline{\chi}^*$, is based on the notion of free Fisher information. To generalize the definition of Fisher information to the non-commutative setting, D. Voiculescu noticed that the standard Fisher information can be defined as

$$\phi(\mu) := ||\partial_x^* 1||_{L^2(\mu)}^2$$

with $\partial_x^*$ the adjoint of the derivative $\partial_x$ for the scalar product in $L^2(\mu)$, that is that for every $f \in L^2(\mu)$,

$$\int f\partial_x^* 1 d\mu(x) = \int \partial_x f d\mu(x).$$

When $\mu$ has a density $\rho$ with respect to Lebesgue measure, we simply have $\partial_x^* 1 = -\partial_x \log \rho = -\rho^{-1}\partial_x \rho$. The entropy $S^*$ is related to Fisher information by the formula

$$S^*(\mu) = \frac{1}{2} \int_0^1 \phi(\mu_t^b) dt$$

with $\mu_t^b$ the law at time $t$ of the Brownian bridge between $\delta_0$ and $\mu$.

Such a definition can be naturally extended to the free probability setting as follows. To this end, we begin by describing the definition of the derivative in this setting; let $P \in \mathbb{C}\langle X_1, \ldots, X_m \rangle$ be for instance a monomial function $\prod_{1 \le k \le r}^{\rightarrow} X_{i_k}$. Then, for any operators $\mathbf{X}_1, \cdots, \mathbf{X}_m$ and $\mathbf{Y}_1, \cdots, \mathbf{Y}_m$,

$$P(\mathbf{X} + \epsilon \mathbf{Y}) - P(\mathbf{X}) = \epsilon \sum_{k=1}^r \prod_{1 \le l \le k-1}^{\rightarrow} \mathbf{X}_{i_l} \mathbf{Y}_{i_k} \prod_{k+1 \le l \le r}^{\rightarrow} \mathbf{X}_{i_l} + O(\epsilon^2). \quad (7.1.2)$$

In order to keep track of the place where the operator $\mathbf{Y}$ have to be inserted, the derivative is defined as follows; the derivative $D_{X_i}$ with respect to the $i^{th}$ variable is a linear form from $\mathbb{C}\langle X_1, \ldots, X_m \rangle$ into $\mathbb{C}\langle X_1, \ldots, X_m \rangle \otimes \mathbb{C}\langle X_1, \ldots, X_m \rangle$ satisfying the non-commutative Leibniz rule

$$D_{X_i} PQ = D_{X_i} P \times 1 \otimes Q + P \otimes 1 \times D_{X_i} Q$$

for any $P, Q \in \mathbb{C}\langle \mathbf{X}_1, \ldots, \mathbf{X}_m \rangle$ and

$$D_{X_i} \mathbf{X}_l = 1_{l=i} 1 \otimes 1.$$

One then denotes $\sharp$ the application from $\mathbb{C}\langle X_1, \ldots, X_m \rangle \otimes \mathbb{C}\langle X_1, \ldots, X_m \rangle \times \mathbb{C}\langle X_1, \ldots, X_m \rangle$ into $\mathbb{C}\langle X_1, \ldots, X_m \rangle$ such that $A \otimes B \sharp C = ACB$. It is then not hard to see that (7.1.2) reduces to

$$P(\mathbf{X} + \epsilon \mathbf{Y}) - P(\mathbf{X}) = \epsilon \sum_{j=1}^m D_{X_j} P \sharp \mathbf{Y}_j + O(\epsilon^2).$$

The cyclic derivative $\mathcal{D}_{X_i}$ with respect to the $i$-th variable is given by

$$\mathcal{D}_{X_i} = m \circ D_{X_i}$$

where $m : \mathbb{C}\langle X_1, \cdots, X_n \rangle \otimes \mathbb{C}\langle X_1 \cdots X_n \rangle \to \mathbb{C}\langle X_1, \cdots, X_n \rangle$ is such that $m(P \otimes Q) = QP$. In the case $m = 1$, using the bijection between $\mathbb{C}\langle X \rangle \otimes \mathbb{C}\langle X \rangle$ and $\mathbb{C}\langle X, Y \rangle$, we find that

$$D_X P(x, y) = \frac{P(x) - P(y)}{x - y}.$$

The analogue of $\partial_x^* 1$ in $L^2(\mu)$ is given, for any $\tau \in \mathcal{M}^{(m)}$, as the element in $L^2(\tau)$ such that for any $P \in \mathbb{C}\langle X_1, \ldots, X_m \rangle$,

$$\tau \otimes \tau \left( D_{X_i} P \right) = \tau(P D_{X_i}^* 1).$$

In the case $m = 1$ and so $\tau \in \mathcal{P}(\mathbb{R})$, it is not hard to check (see [123]) that

$$D_X^* 1(x) = 2PV \left( \int (x-y)^{-1} d\tau(y) \right).$$

Free Fisher information is thus given, for $\tau \in \mathcal{M}^{(m)}$, by

$$\Phi^*(\tau) = \sum_{i=1}^{m} \|D_{X_i}^* 1\|_{L^2(\tau)}^2$$

and therefore $\underline{\chi}^*$ is defined by

$$\underline{\chi}^*(\mu) := -\frac{1}{2} \int_0^1 \Phi^*(\mu_t^b) dt$$

with $\mu_t^b$ the distribution of the free Brownian bridge, $\mu_t^b = \mathcal{L}(tX_i + \sqrt{t(1-t)} S_i, 1 \leq i \leq m)$ if $\mu = \mathcal{L}(X_i, 1 \leq i \leq m)$ and $S_i$ are free standard semi-circular variables, free with $(X_i, 1 \leq i \leq m)$.

The conjecture (named 'unification problem' by D. Voiculescu [129]) is

*Conjecture* 7.1: For any $\mu$ such that $\underline{\chi}(\mu) > -\infty$,

$$\underline{\chi}(\mu) = \underline{\chi}^*(\mu).$$

*Remark* 7.2: If the conjecture would hold for any $\mu \in \mathcal{M}^{(m)}$ and not only for $\mu$ with finite microstates entropy, it would provide an affirmative answer to Connes question since it is known that if $\mu$ is the law of $X_1, \cdots, X_m$ and $S_1, \cdots, S_m$ free semicircular variables, free with $\mathbf{X}_1, \cdots, \mathbf{X}_m$ then, for any $\epsilon > 0$ the distribution $\mu \boxplus \sigma_\epsilon$ of $(\mathbf{X}_1 + \epsilon \mathbf{S}_1, \cdots, \mathbf{X}_m + \epsilon \mathbf{S}_m)$ satisfies $\underline{\chi}^*(\mu \boxplus \sigma_\epsilon) > -\infty$, and hence the above equality would imply $\underline{\chi}(\mu \boxplus \sigma_\epsilon) > -\infty$ so that one could find matrices whose empirical distribution approximates $\mu \boxplus \sigma_\epsilon$ and thus $\mu$ since $\epsilon$ can be chosen arbitrarily small.

In [18], we proved that

**Theorem 7.3.** *For any $\mu \in \mathcal{M}^{(m)}$,*

$$\underline{\chi}(\mu) \leq \underline{\chi}^*(\mu).$$

*Moreover, we can define another entropy $\underline{\chi}^{**}$ such that*

$$\underline{\chi}(\mu) \geq \underline{\chi}^{**}(\mu).$$

Typically, $\underline{\chi}^*$ is obtained as an infimum of a rate function on lows of non-commutative processes with given terminal data, whereas $\underline{\chi}^{**}$ is the infimum of the same rate function but on a a priori smaller set.

From this result, we as well obtain bounds on the entropy dimension

**Corollary 7.4.** *Let for $\tau \in \mathcal{M}^{(m)}$*

$$\delta(\tau) := m - \limsup_{\epsilon \to 0} \frac{\chi(\tau \boxplus \sigma_\epsilon)}{\log \epsilon} = m - \limsup_{\epsilon \to 0} \frac{\underline{\chi}(\tau \boxplus \sigma_\epsilon)}{\log \epsilon} \qquad (7.1.3)$$

*where $\tau \boxplus \sigma_\epsilon$ stands for the free convolution by $m$ free semi-circular variables with parameter $\epsilon > 0$. Define accordingly $\delta^*, \delta^{**}$. Then*

$$\delta^{**}(\tau) \leq \delta(\tau) \leq \delta^*(\tau).$$

In a recent work, A. Connes and D. Shlyaktenkho [35] defined another quantity $\Delta$, candidate to be an invariant for von Neumann algebras, by generalizing the notion of $L^2$-homology and $L^2$-Betti numbers for a tracial von Neumann algebra. Such a definition is in particular motivated by the work of D. Gaboriau [55]. They could compare $\Delta$ with $\delta^*$, and therefore, thanks to the above corollary, to $\delta$.

Eventhough $\delta^*, \delta^{**}$ are not simple objects, they can be computed in some cases such as for the law of the DT-operators (see [1]) or in the case of a finitely generated group where I. Mineyev and D. Shlyakhtenko [95] proved that

$$\delta^*(\tau) = \beta_1(G) - \beta_0(G) + 1$$

with the group $L^2$ Betti-numbers $\beta$.

DT-operators have long been an interesting candidate to try to disprove the invariance of $\delta$. A DT-operator $T$ can be constructed as the limit in distribution of upper triangular matrices with i.i.d Gaussian entries (which amounts to consider the law of two self-adjoint non-commutative variables $T + T^*$ and $i(T - T^*)$). If $C$ is a semicircular operator, which is the limit in distribution of the (non Hermitian) matrix with i.i.d Gaussian entries, then $C$ can be written as $T + \widetilde{T}^*$ where $T, \widetilde{T}$ are free copies of $T$. Hence, since $\delta(C) \leq 2$, we can hope that $\delta(T) < 2$. However, based on an heavy computation of moments of these DT-operators due to P. Sniady [111], K. Dykema and U. Haagerup [44] could prove that $T$ generates $L(F^2)$. Hence invariance would be disproved if $\delta(T) < 2$. But in fact, L. Aagaard [1] recently proved that $\delta^*(T) = 2$ which shows at least that $T$ is not a counter-example for the invariance of $\delta^*$ ( and also settle the case for $\delta$ if one believes conjecture 7.1).

We now give the main ideas of the proof of Theorem 7.3.

## 7.2. Large deviation upper bound for the law of the process of the empirical distribution of Hermitian Brownian motions

In [29] and [30], we established with T. Cabanal Duvillard the inequality

*A. Guionnet/Large deviations for random matrices* 154**Theorem 7.5.** *There exists a good rate function $\mathcal{S}$ on $\mathcal{C}([0,1], \mathcal{M}^{(m)})$ such that for any $\mu \in \mathcal{M}^{(m)}$*

$$\underline{\chi}(\mu) \leq -\mathcal{S}(\mu^b) \tag{7.2.4}$$

*and*

$$-\mathcal{S}(\mu^b) \geq \underline{\chi}^*(\mu). \tag{7.2.5}$$

*Inequality (7.2.5) is an equality as soon as $(D^*_{X_1} 1, \ldots, D^*_{X_m} 1)$ is in the cyclic gradient space, i.e*

$$\inf_{F \in \mathcal{C}^1([0,1], \mathcal{CC}^m_{st}(\mathbb{R}))} \{ \int_0^1 \sum_{l=1}^m \mu^b_u (|\mathcal{D}_{X_l} F_u - D^*_{X_l} 1 - \frac{\mathbf{X}_l}{u}|^2) du \} = 0$$

*where the infimum runs over the set $\mathcal{C}^1([0,1], \mathcal{CC}^m_{st}(\mathbb{R}))$ of continuously differentiable functions with values in $\mathcal{CC}^m_{st}(\mathbb{R})$, the restriction of self-adjoint non-commutative functions of $\mathcal{CC}^m_{st}(\mathbb{C})$ defined in (6.2.2). Further, by definition, if $m : \mathbb{C}\langle X_1, \ldots, X_n \rangle \times \mathbb{C}\langle X_1, \ldots, X_n \rangle \to \mathbb{C}\langle X_1, \ldots, X_n \rangle$ is such that $m(P \otimes Q) = QP$, the cyclic derivative is given by $\mathcal{D}_{X_l} = m \circ D_{X_l}$.*

Note here that it is not clear whether $(D^*_{X_1} 1, \ldots, D^*_{X_m} 1)$ should belong to the cyclic gradient space or not in general. This was proved by D. Voiculescu [126] when it is polynomial, and it seems quite natural that this should be the case for states with finite entropy (see a discussion in [30]).

The strategy is here exactly the same than in chapter 4; consider $m$ independent Brownian motion $(\mathbf{H}^N_1, \cdots, \mathbf{H}^N_m)$ and the $\mathcal{M}^{(m)}$-valued process $\hat{\mu}^N_t = \hat{\mu}^N_{\mathbf{H}^N_1(t), \cdots, \mathbf{H}^N_m(t)}$. Then, it can be seen thanks to Ito's calculus that, for any $F \in \mathcal{C}([0,1], \mathcal{CC}^m_{st}(\mathbb{R}))$, if we set

$$\mathcal{S}^{s,t}(\nu, F) = \nu_t(F_t) - \nu_s(F_s) - \int_s^t \nu_u(\partial_u F_u) du$$
$$- \frac{1}{2} \int_s^t \nu_u \otimes \nu_u (\sum_{l=1}^m D_{X_l} \circ \mathcal{D}_{X_l} F_u) du$$

$$\lll F, G \ggg^{s,t}_\nu = \sum_{l=1}^m \int_s^t \nu_u (\mathcal{D}_{X_l} F_u \mathcal{D}^*_{X_l} G_u) du$$

$$\mathcal{S}^{s,t}(\nu) = \sup_{F \in \mathcal{CC}_{st}(\mathbb{R} \times [0,1])} (\mathcal{S}^{s,t}(\nu, F) - \frac{1}{2} \lll F, F \ggg^{s,t}_\nu),$$

then $(\mathcal{S}^{0,t}(\hat{\mu}^N, F), t \geq 0)$ is a martingale with bracket given by $N^{-2} \lll F, F \ggg^{0,t}_{\hat{\mu}^N}$. Hence we are in business and we can prove a large deviation upper bound with good rate function which is infinite if the initial law is not the law of null operators, and otherwise given by $\mathcal{S}^{0,1}$. We can proceed as in Chapter 6 to improve the upper bound when we are considering the deviations of $\hat{\mu}^N_1$ and see that the infimum is taken at a free Brownian bridge. The relation with $\underline{\chi}^*$ comes from the fact that the free Brownian bridge is associated with the field

$$\tau_\mu \left( \frac{X - X_s}{1-s} | X_s \right) = \frac{X_s}{s} - \mathcal{J}^{\mu^b_s}$$

where $\mathcal{J}^\mu = D^*1$ denotes the non-commutative Hilbert transform of $\mu \in \mathcal{M}^{(m)}$ and $\mu_t^b = \tau_\mu \circ \pi_t$ is the time marginal of the free Brownian bridge. This equality is a direct consequence of Corollary 3.9 in [124].

The main problem here to hope to obtain a large deviation lower bound is that the associated Fokker- Planck equations are quite hard to study and we could not find a reasonable condition over the fields to obtain uniqueness, as needed (see Section 4.2.1).

This is the reason why the idea to study large deviation on path space emerged; the lower bound estimates become easier since we shall then deal with uniqueness of strong solutions to these Fokker-Planck equations rather than uniqueness of weak solutions.

### 7.3. Large deviations estimates for the law of the empirical distribution on path space of the Hermitian Brownian motion

In this section, we consider the empirical distribution on path space of independent Hermitian Brownian motions $(\mathbf{H}^1, .., \mathbf{H}^m)$. It is described by the quantities

$$\hat{\sigma}^N(F) = \text{tr}\left(P(\mathbf{H}^{N,i_1}(t_1), \mathbf{H}^{N,i_2}(t_2), \ldots, \mathbf{H}^{N,i_n}(t_n))\right)$$

with $F(x^1, .., x^m) = P(x^{i_1}(t_1), \cdots, x^{i_n}(t_n))$ for any choice of non-commutative test function $P$, any $(t_1, \ldots, t_n) \in [0,1]^n$ and $(i_1, \cdots, i_n) \in \{1, \cdots, m\}$.

The study of the deviations of the law of $\hat{\sigma}^N$ could a priori be performed again thanks to Itô's formula by induction over the number of time marginals. However, this is not a good idea here since we could not obtain the lower bound estimates. Hence, we produced a new idea to generate exponential martingales which is based on the Clark-Ocone formula.

### 7.4. Statement of the results

Let us be more precise. In order to avoid the problems of unboundedness of the operators and still get a nice topology, manageable for free probabilists, we considered in [18] the unitary operators

$$\mathbf{U}_t^{N,l} := \psi(H_t^{N,l}), \qquad \text{with } \psi(x) = (x+4i)(x-4i)^{-1}.$$

If
$$\Omega_t^l := (\mathbf{S}_t^l + 4i)(\mathbf{S}_t^l - 4i)^{-1}$$

with a free Brownian motion $(\mathbf{S}^1, \ldots, \mathbf{S}^m)$, it is not hard to see that

$$\frac{1}{N}Tr((\mathbf{U}_{t_1}^{N,i_1})^{\varepsilon_1} \ldots (\mathbf{U}_{t_n}^{N,i_n})^{\varepsilon_n}) \to_{N\to\infty} \varphi((\Omega_{t_1}^{i_1})^{\varepsilon_1} \ldots (\Omega_{t_n}^{i_n})^{\varepsilon_n})$$

and we shall here study the deviations with respect to this typical behavior. Since the $\mathbf{U}^{N,l}$ are uniformly bounded, polynomial test functions provide a good topology. This amounts to restrict ourselves to a few Stieljes functionals of the Hermitian Brownian motions. However, this is already enough to study Voiculescu's entropy when one considers the deviations toward laws of bounded operators since then the polynomial functions of $(\psi(\mathbf{X}_1), \cdots, \psi(\mathbf{X}_m))$ generates the set of polynomial functions of $(\mathbf{X}_1, \cdots, \mathbf{X}_m)$ and vice-versa (see Lemma 7.7).

Let $\mathcal{F}_{[0,1]}^m$ be the $*$-algebra of the free group generated by $(u_t^i; t \in [0,1], i \in \{1, \ldots, m\})$ (that is the set of polynomial functions generated by $\prod^{\rightarrow}(u_{t_k}^{i_k})_k^\epsilon$, $t_k \in [0,1]$, $i_\kappa \in \{1, \cdots, m\}$, $\epsilon_k = 1, -1$, equipped with the involution $(u_t^i)^* = (u_t^i)^{-1}$) and $\mathcal{F}_{[0,1]}^{m,sa}$ be its self-adjoint elements.

We denote $\mathcal{M}(\mathcal{F}_{[0,1]}^m)$ the set of tracial states on $\mathcal{F}_{[0,1]}^m$ (i.e. the set of linear forms on $\mathcal{F}_{[0,1]}^m$ satisfying the properties of positiveness, total mass equal to one, and traciality and with real restriction to $\mathcal{F}_{[0,1]}^{m,sa}$). We equip $\mathcal{M}(\mathcal{F}_{[0,1]}^m)$ with its weak topology with respect to $\mathcal{F}_{[0,1]}^m$. Let $\mathcal{M}^c(\mathcal{F}_{[0,1]}^m)$ be the subset of $\mathcal{M}(\mathcal{F}_{[0,1]}^m)$ of states such that for any $\varepsilon_1, \ldots, \varepsilon_n \in \{-1, +1\}^n$, and any $i_1, \ldots, i_n \in \{1, \ldots, m\}^n$, the quantity $\tau((u_{t_1}^{i_1})^{\varepsilon_1} \ldots (u_{t_n}^{i_n})^{\varepsilon_n})$ is continuous in the variables $t_1, \ldots, t_n$. Remark that for any process of unitary operators $(\mathbf{U}^1, \ldots, \mathbf{U}^m)$ with values in a $W^*$-probability space $(\mathcal{A}, \varphi)$ we can associate $\tau \in \mathcal{M}(\mathcal{F}_{[0,1]}^m)$ such that

$$\tau(P) = \varphi(P(\mathbf{U})).$$

Reciprocally, GNS construction allows us to associate to any $\tau \in \mathcal{M}(\mathcal{F}_{[0,1]}^m)$ a $W^*$-probability space $(\mathcal{A}, \varphi)$ as above.

In particular, if $\mathbf{S} = (\mathbf{S}^1, \ldots, \mathbf{S}^m)$ is a m-dimensional free Brownian motion, the family $(\frac{\mathbf{S}_t^l + 4i}{\mathbf{S}_t^l - 4i}; t \in [0,1], l \in \{1, \ldots m\})$ defines a state $\sigma \in \mathcal{M}^c(\mathcal{F}_{[0,1]}^m)$. To define our rate function, we need to introduce, for any time $t \in [0,1]$, any process $\mathbf{X} = (\mathbf{X}_t^i; t \in [0,1], i \in \{1, \ldots, m\})$ of self-adjoint operators and a m-dimensional free Brownian motion $\mathbf{S} = (\mathbf{S}^1, \ldots, \mathbf{S}^m)$, $\mathbf{X}$ and $\mathbf{S}$ being free, the process $\mathbf{X}_\cdot^t = \mathbf{X}_{\cdot \wedge t} + \mathbf{S}_{\cdot - t \vee 0}$. If $\tau \in \mathcal{M}^c(\mathcal{F}_{[0,1]}^m)$ is the law of $(\psi(\mathbf{X}_t^i); t \in [0,1], i \in \{1, \ldots, m\})$, we set $\tilde{\tau}^t \in \mathcal{M}^c(\mathcal{F}_{[0,1]}^m)$ the distribution of $(\psi(\mathbf{X}_s^{t,i}); t \in [0,1], i \in \{1, \ldots, m\})$. Finally, we denote $\nabla_t$ the non-commutative Malliavin operator given by

$$\nabla_s^l((u_{t_1}^{i_1})^{\varepsilon_1} \ldots (u_{t_n}^{i_n})^{\varepsilon_n}) = -\sum_{p=1}^n 1_{i_p = l} \frac{\varepsilon_p}{8i}((u_{t_p}^{i_p})^{\varepsilon_p} - 1)(u_{t_{p+1}}^{i_{p+1}})^{\varepsilon_{p+1}} \ldots (u_{t_n}^{i_n})^{\varepsilon_n} \times \\ (u_{t_1}^{i_1})^{\varepsilon_1} \ldots (u_{t_{p-1}}^{i_{p-1}})^{\varepsilon_{p-1}}((u_{t_p}^{i_p})^{\varepsilon_p} - 1)1_{[0,t_p]}(s).$$

Note that this definition is a formal extension of

$$\nabla_s^l(x_{t_1}^{i_1} \ldots x_{t_n}^{i_n}) = \sum_{p=1}^n 1_{i_p = l} x_{t_{p+1}}^{i_{p+1}} \ldots x_{t_n}^{i_n} x_{t_1}^{i_1} \ldots x_{t_{p-1}}^{i_{p-1}} 1_{[0,t_p]}(s)$$

to the (non converging in general) series $u_t^j = -(1/4i)(x_t^j + 4i) \sum_k (x_t^j/4i)^k$, $j \in \{1, \ldots, m\}$. Finally, we denote $\mathcal{B}_t$ the $\sigma$ algebra generated by $\{\mathbf{X}_u^l, u \leq t, 1 \leq l \leq m\}$ and for any $\tau \in \mathcal{M}(\mathcal{F}_{[0,1]}^m)$, $\tau(.|\mathcal{B}_t)$ the conditional expectation knowing $\mathcal{B}_t$, i.e. the projection on $\mathcal{B}_t$ in $L^2(\tau)$. We proved that

**Theorem 7.6.** *Let* $I : \mathcal{M}(\mathcal{F}^m_{[0,1]}) \to \mathbb{R}^+$ *being given by*

$$I(\tau) = \sup_{t\in[0,1]} \sup_{F\in\mathcal{F}^{m,sa}_{[0,1]}} \{\widetilde{\tau}^t(F) - \sigma(F) - \frac{1}{2}\int_0^t \widetilde{\tau}^t\left(\widetilde{\tau}^s(\nabla_s F|\mathcal{B}_s)^2\right) ds\}.$$

*Then*

1. *$I$ is a good rate function.*
2. *Any $\tau \in \{I < \infty\}$ is such that there exists $K^\tau \in L^2(\tau \times ds)$ such that*
   *a)* $\inf_{P\in\mathcal{F}^{m,sa}_{[0,1]}} \int_0^1 \tau\left[|\widetilde{\tau}^s(\nabla_s P|\mathcal{B}_s) - K^\tau_s|^2\right] ds = 0.$
   *b) For any $P \in \mathcal{F}^m_{[0,1]}$, and $t \in [0,1]$, we have*

$$\widetilde{\tau}^t(P) = \widetilde{\tau}^0(P) + \int_0^t \tau\left(\widetilde{\tau}^s(\nabla_s P|\mathcal{B}_s) . K^\tau_s\right) ds.$$

3. *For any closed set $F \subset \mathcal{M}^c(\mathcal{F}^m_{[0,1]})$ we have*

$$\limsup_{N\to\infty} \frac{1}{N^2} \log \mathbb{P}(\hat{\sigma}^N \in F) \leq -\inf_{\tau\in F} I(\tau)$$

4. *Let $\mathcal{M}^{c,\infty}_b(\mathcal{F}^m_{[0,1]})$ be the set of states in $\{I < \infty\}$ such that the infimum in 2.a) is attained. Then, for any open set $O \subset \mathcal{M}(\mathcal{F}^m_{[0,1]})$,*

$$\liminf_{N\to\infty} \frac{1}{N^2} \log \mathbb{P}(\hat{\sigma}^N \in O) \geq - \inf_{\tau\in O\cap \mathcal{M}^{c,\infty}_b(\mathcal{F}^m_{[0,1]})} I(\tau).$$

### 7.4.1. Application to Voiculescu's entropies

To relate Thorem 7.6 to estimates on $\chi$, let us introduce in the spirit of Voiculescu, the entropy $\widetilde{\chi}$ as follows. Let $\Gamma^U_R(\nu, n, N, \epsilon)$ be the set of unitary matrices $V_1, .., V_m$ such that

$$-1 \leq 2^{-1}(V_j + V_j^*) \leq 1 - 2(R^2 + 1)^{-1}$$

for all $j \in \{1, \cdots, m\}$ and

$$|\nu(\mathbf{U}^{\varepsilon_1}_{i_1}..\mathbf{U}^{\varepsilon_p}_{i_p}) - \mathrm{tr}(V^{\varepsilon_1}_{i_1}..V^{\varepsilon_p}_{i_p})| < \epsilon$$

for any $1 \leq p \leq n$, $i_1, .., i_p \in \{1, .., m\}^p$, $\varepsilon_1, \cdots, \varepsilon_p \in \{-1, +1\}^p$. We set

$$\widetilde{\chi}(\nu) := \sup_{R>0} \inf_{n\in\mathbb{N}} \inf_{\epsilon>0} \limsup_{N\to\infty} \frac{1}{N^2} \log P\bigl((\psi(A_1^N), \cdots, \psi(A_m^N)) \in \Gamma^U_R(\nu, n, N, \epsilon)\bigr).$$

Let $\Psi$ be the Möbius function

$$\Psi(\mathbf{X}^1, \ldots, \mathbf{X}^m)_l = \psi(\mathbf{X}^l) = \frac{\mathbf{X}^l + 4i}{\mathbf{X}^l - 4i}.$$

Then, it is not hard to see that

**Lemma 7.7.** *For any $\tau \in \mathcal{M}_1^{(m)}$ which corresponds, by the GNS construction, to the distribution of bounded operators, we have*

$$\underline{\chi}(\tau) = \widetilde{\chi}(\tau \circ \Psi).$$

Moreover, Theorem 7.6 and the contraction principle imply

**Theorem 7.8.** *Let $\pi_1 : \mathcal{M}^{c,\infty}(\mathcal{F}_{[0,1]}^m) \to \mathcal{M}(\mathcal{F}_{\{1\}}^m) = \mathcal{M}_U^m$ be the projection on the algebra generated by $(\mathbf{U}^1(1), \ldots, \mathbf{U}^m(1), \mathbf{U}^1(1)^{-1}, \ldots, \mathbf{U}^m(1)^{-1})$. Then, for any $\tau \in \mathcal{M}^{(m)}$,*

$$\underline{\chi}^{**}(\tau) := -\lim_{\delta \to 0} \inf\{I(\sigma); d(\sigma \circ \pi_1, \tau \circ \Psi) < \delta, \sigma \in \mathcal{M}_b^{c,\infty}(\mathcal{F}_{[0,1]}^m)\} \leq \underline{\chi}(\tau) \quad (7.4.6)$$

*and*

$$\begin{aligned}\underline{\chi}(\tau) &\leq -\lim_{\delta \to \infty} \inf\{I(\sigma); d(\sigma \circ \pi_1, \tau \circ \Psi) < \delta, \sigma \in \mathcal{M}^{c,\infty}(\mathcal{F}_{[0,1]}^m)\} \\ &= -\inf\{I(\sigma), \sigma \circ \pi_1 = \tau \circ \Psi\}.\end{aligned} \quad (7.4.7)$$

The above upper bound can be improved by realizing that the infimum has to be achieved at a free Brownian bridge (generalizing the ideas of Chapter 6, Section 6.7). In fact, if $\mu$ is the distribution of $m$ self-adjoint operators $\{\mathbf{X}_1, \ldots, \mathbf{X}_m\}$ and $\{\mathbf{S}^1, \ldots, \mathbf{S}^m\}$ is a free Brownian motion, free with $\{\mathbf{X}_1, \ldots, \mathbf{X}_m\}$, we denote $\tau_\mu^b$ the distribution of
$\left\{\psi(t\mathbf{X}^l + (1-t)\mathbf{S}_{\frac{t}{1-t}}^l), 1 \leq l \leq m, t \in [0,1]\right\}$. Then

**Theorem 7.9.** *For any $\mu \in \mathcal{M}^{(m)}$,*

$$\underline{\chi}(\mu) \leq -I(\tau_\mu^b) = \underline{\chi}^*(\mu) \quad (7.4.8)$$

### 7.4.2. Proof of Theorem 7.6

We do not prove here that $I$ is a good rate function. The idea of the proof of the large deviations estimates is again based on the construction of exponential martingales; In fact, for $P \in \mathcal{F}_{[0,1]}^m$, it is clear that $M_P^N(t) = E[\hat{\sigma}^N(P)|\mathcal{F}_t] - E[\hat{\sigma}^N(P)]$ is a martingale for the filtration $\mathcal{F}_t$ of the Hermitian Brownian motions. Clark-Ocone formula gives us the bracket of this martingale :

$$< M_P^N >_t = \frac{1}{N^2} \sum_{l=1}^m \int_0^t \mathrm{tr}(E[\nabla_s^l P|\mathcal{F}_s]^2) ds$$

As a consequence, for any $P \in \mathcal{F}_{[0,1]}^m$, and $t \in [0,1]$

$$E[\exp\{N^2(E[\hat{\sigma}^N(P)|\mathcal{F}_t] - E[\hat{\sigma}^N(P)] - \int_0^t \mathrm{tr}(E[\nabla_s P|\mathcal{F}_s]^2) ds)\}] = 1.$$

To deduce the large deviation upper bound from this result, we need to show that

**Proposition 7.10.** *For $P \in \mathcal{F}_{[0,1]}^m$, $\tau \in \mathcal{M}^{c,\infty}(\mathcal{F}_{[0,1]}^m)$, $\epsilon > 0$, $l \in \mathbb{N}$, $L \in \mathbb{R}^+$, define*

$$\begin{aligned} H &:= H(P, L, \epsilon, N, \tau, l) \\ &= \operatorname*{ess\,sup}_{\{d(\hat{\sigma}^N, \tau) < \epsilon; \hat{\sigma}^N \in K_{L\sqrt{\cdot}} \cap \Gamma_L\}} |tr(\mathbb{E}[P(U^N)|\mathcal{H}_t]^l) - \tau(\widetilde{\tau}^t(P|\mathcal{B}_t)^l)| \end{aligned}$$

*then one has, for every $l \in \mathbb{N}$*

$$\sup_{L>0} \limsup_{\epsilon \to 0} \limsup_{N \to \infty} \sup_{\substack{\tau \in K_{L\sqrt{\cdot}} \cap \Gamma_L \\ t \in [0,1]}} H = 0 \qquad (7.4.9)$$

*Here $\{K_{L\sqrt{\cdot}} \cap \Gamma_L\}_{L \in \mathbb{N}}$ are compact subsets of $\mathcal{M}(\mathcal{F}_{[0,1]}^m)$ such that*

$$\limsup_{L \to \infty} \limsup_{N \to \infty} \frac{1}{N^2} \log \mathbb{P}(\hat{\sigma}^N \in (K_{L\sqrt{\cdot}} \cap \Gamma_L)^c) = -\infty.$$

Then, we can apply exactly the same techniques than in Chapter 4 to prove the upper bound.

To obtain the lower bound, we have to obtain uniqueness criteria for equations of the form

$$\widetilde{\tau}^t(P) - \sigma(P) = \int_0^t \tau\left(\widetilde{\tau}^s(\nabla_s P|\mathcal{B}_s)\widetilde{\tau}^s(\nabla_s K|\mathcal{B}_s)\right) ds$$

with fields $K$ as general as possible. We proved in [18], Theorem 6.1, that if $K \in \mathcal{F}_{[0,1]}^m$, the solutions to this equation are strong solutions in the sense that there exists a free Brownian motion $\mathbf{S}$ such $\tau$ is the law of the operator $\mathbf{X}$ satisfying

$$d\mathbf{X}_t = d\mathbf{S}_t + \widetilde{\tau}^t(\nabla_t K|\mathcal{B}_t)(\mathbf{X})dt.$$

But, if $K \in \mathcal{F}_{[0,1]}^m$, it is not hard to see that $\widetilde{\tau}^t(\nabla_t K|\mathcal{B}_t)(\mathbf{X})$ is Lipschitz operator, so that we can see that there exists a unique such operator $\mathbf{X}$, implying the uniqueness of the solution of our free differential equation, and hence the large deviation lower bound.

### 7.5. Discussion and open problems

Note that we have the following heuristic description of $\underline{\chi}^*$ and $\underline{\chi}^{**}$:

$$\underline{\chi}^*(\tau) = -\inf\{\int_0^1 \mu(K_t^2) dt\}$$

where the infimum is taken over all laws $\mu$ of non-commutative processes which are null operators at time 0, operators with law $\tau$ at time one and which are the distributions of 'weak solutions' of

$$d\mathbf{X}_t = d\mathbf{S}_t + K_t(X)dt.$$

$\chi^{**}$ is defined similarly but the infimum is restricted to processes with smooth fields $K$ (actually $K \in \mathcal{F}_{[0,1]}^m$). We then have proved in Theorem 7.8 that

$$\underline{\chi}^{**} \leq \underline{\chi} \leq \underline{\chi}^*$$

and it is legitimate to ask when $\underline{\chi}^{**} = \underline{\chi}^*$. Such a result would show $\underline{\chi} = \underline{\chi}^*$. Note that in the classical case, the relative entropy can actually be described by the above formula by replacing the free Brownian motion by a standard Brownian motion and then all the inequalities become equalities.

This question raises numerous questions :

1. First, inequalities (7.4.6) and (7.4.8) become equalities if $\tau_\mu^b \in \mathcal{M}_b^{c,\infty}(\mathcal{F}_{[0,1]}^m)$ that is if there exists $n$, times $(t_i, 1 \leq i \leq n+1) \in [0,1]^{n+1}$, and polynomial functions $(Q_i, 1 \leq i \leq n)$ and $P$ such that

$$D_{\mu_t^b}^* 1 = \mathcal{J}^{\mu_t^b} = \sum_{i=1}^n 1_{(t_i, t_{i+1}]} Q_i + \widetilde{\tau}_\mu^{b\,t} (\nabla_t P | \mathcal{B}_t).$$

Can we find non trivial $\mu \in \mathcal{M}^{(m)}$ such that this is true?

2. If we follow the ideas of Chapter 4, to improve the lower bound, we would like to regularize the laws by free convolution by free Cauchy variables $\mathbf{C}^\epsilon = (C_1^\epsilon, \cdots, C_m^\epsilon)$ with covariance $\epsilon$. If $\mathbf{X} = (\mathbf{X}_1, \cdots, \mathbf{X}_m)$ is a process satisfying

$$d\mathbf{X}_t = d\mathbf{S}_t + K_t(X_t)dt,$$

for some non-commutative function $K_t$, it is easy to see that $\mathbf{X}^\epsilon = \mathbf{X} + \mathbf{C}^\epsilon$ satisfies the same free Fokker-Planck equation with $K_t^\epsilon(\mathbf{X}_t + \mathbf{C}^\epsilon) = \tau(K_t(\mathbf{X}_t)|\mathbf{X}_t + \mathbf{C}^\epsilon)$. Then, does $K^\epsilon$ is smooth with respect to the operator norm? This is what we proved for one operator in [65]. If this is true in higher dimension, then Connes question is answered positively since by Picard argument

$$d\mathbf{X}_t^\epsilon = d\mathbf{S}_t + K_t^\epsilon(X_t)dt$$

has a unique strong solution and there exists a smooth function $F^\epsilon$ such that for any $t > 0$

$$\mathbf{X}_t^\epsilon = F_t^\epsilon(\mathbf{S}_s, s \leq t).$$

In particular, for any polynomial function $P \in \mathbb{C}\langle X_1, \ldots, X_m\rangle$

$$\mu(P(\mathbf{X} + \mathbf{C}^\epsilon)) = \sigma(P \circ F_1^\epsilon(\mathbf{S}_s, s \leq 1)) = \lim_{N \to \infty} \text{tr}(P \circ F_1^\epsilon(\mathbf{H}_s^N, s \leq 1))$$

where we used in the last line the smoothness of $F_1^\epsilon$ as well as the convergence of the Hermitian Brownian motion towards the free Brownian motion. Hence, since $\epsilon$ is arbitrary, we can approximate $\mu$ by the empirical distribution of the matrices $F_1^\epsilon(\mathbf{H}_s^N, s \leq 1)$, which would answer Connes question positively. As in remark 7.2, the only way to complete the argument without dealing with Connes question would be to be able to prove such a regularization property only for laws with finite entropy, but it is rather unclear how such a condition could enter into the game. This could only be true if the hyperfinite factor would have specific analytical properties.

3. If we think that the only point of interest is what happens at time one, then we can restrict the preceding discussion by showing that if $\mathbf{X} = (\mathbf{X}_1, \cdots, \mathbf{X}_m)$ are non-commutative variables with law $\mu$, and $(\mathcal{J}_i^\mu, 1 \leq i \leq m)$ is the Hilbert transform of $\mu$ and if we let $\mu^\epsilon$ be the law of $\mathbf{X} + \mathbf{C}^\epsilon$, then we would like to show that $(\mathcal{J}_i^{\mu^\epsilon}, 1 \leq i \leq m)$ is smooth for $\epsilon > 0$. In the case $m = 1$, $\mathcal{J}^{\mu^\epsilon}$ is analytic in $\{|\Im(z)| < \epsilon\}$. The generalization to higher dimension is wide open.

4. A related open question posed by D. Voiculescu [129] (in a paragraph entitled Technical problems) could be to try to show that the free convolution acts smoothly on Fisher information in the sense that $t \in \mathbb{R}^+ \to \tau_{X+tS}(|\mathcal{J}_i^{\tau_{X+tS}}|^2)$ is continuous.

5. A different approach to microstates entropy could be to study the generating functions $\Lambda(P)$ given, for $P \in \mathbb{C}\langle X_1, \cdots, X_m \rangle \otimes \mathbb{C}\langle X_1, \cdots, X_m \rangle$, by

$$\lim_{R \to \infty} \limsup_{N \to \infty} \frac{1}{N^2} \log \int_{\|X_N^i\|_\infty \leq R} e^{\mathrm{Tr} \otimes \mathrm{Tr}(P(\mathbf{X}_N^1, \cdots, \mathbf{X}_N^m))} \prod_{1 \leq i \leq m} d\mu_N(\mathbf{X}_N^i)$$

It is easy to see (and written down in [67]) that

$$\underline{\chi}(P) = \inf_{P \in \mathbb{C}\langle X_1, \cdots, X_m \rangle^{\otimes 2}} \{\Lambda(P) - \tau \otimes \tau(P)\}.$$

Reciprocally,

$$\Lambda(P) = \sup_{\tau \in \mathcal{M}^{(m)}} \{\underline{\chi}(\tau) + \tau \otimes \tau(P)\}.$$

Therefore, we see that the understanding of the first order of all matrix models is equivalent to that of $\underline{\chi}$. In particular, the convergence of all of their free energies would allow to replace the limsup in the definition of the microstates entropy by a liminf, which would already be a great achievement in free entropy theory. Note also that in the usual proof of Cramer's theorem for commutative variables, the main point is to show that one can restrict the supremum over the polynomial functions $P \in \mathbb{C}\langle X_1, \cdots, X_m \rangle^{\otimes 2}$ to polynomial functions in $\mathbb{C}\langle X_1, \cdots, X_m \rangle$ (i.e. take linear functions of the empirical distribution). This can not be the case here since this would entail that the microstates entropy is convex which it cannot be according to D. Voiculescu [124] who proved actually that if $\tau \neq \tau' \in \mathcal{M}^{(m)}$ with $m \geq 2$, $\tau$ and $\tau'$ having finite microstates entropy, then $\alpha \tau + (1 - \alpha) \tau'$ have infinite entropy for $\alpha \in (0, 1)$.

**Acknowledgments :** I would like to take this opportunity to thank all my coauthors, I had a great time developping this research program with them. I am particularly indebted toward O. Zeitouni for a carefull reading of preliminary versions of this manuscript. I am also extremely grateful to many people who very kindly helped and encouraged me in my struggle for understanding points in fields that I used to ignore entirely, among whom D. Voiculescu, D. Shlyakhtenko, N. Brown, P. Sniady, C. Villani, D. Serre, Y. Brenier, A. Okounkov, S. Zelditch, V. Kazakov, I. Kostov, B. Eynard. I wish also to thank the scientific committee and the organizers of the XXIX conference on Stochastic processes and Applications for giving me the opportunity to write these notes, as well as to discover the amazingly enjoyable style of Brazilian conferences.